\documentclass[12pt,one column]{IEEEtran}
\usepackage{setspace}
\usepackage{amsmath,mathtools}
\usepackage{amssymb}
\usepackage{amsthm}
\usepackage{bbm}
\usepackage{dsfont}
\usepackage{cite}
\usepackage{color}
\usepackage{epsfig}
\usepackage{epsf}
\usepackage{rotating}
\usepackage{mathrsfs}
\usepackage{textgreek}
\usepackage{upgreek}
\usepackage{epsfig}
\usepackage{graphics}
\usepackage{enumerate}
\usepackage{euscript}
\usepackage{hyperref}
\usepackage[dvipsnames]{xcolor}
\usepackage{accents}
\DeclareMathAccent{\wtilde}{\mathord}{largesymbols}{"65}
\usepackage[T3,T1]{fontenc}
\DeclareSymbolFont{tipa}{T3}{cmr}{m}{n}
\DeclareMathAccent{\invbreve}{\mathalpha}{tipa}{16}
\usepackage[multiple]{footmisc}
\usepackage{anyfontsize}
\usepackage{t1enc}
\usepackage{amssymb}
\usepackage{cite}
\usepackage{color}
\usepackage{epsfig}
\usepackage{epsf}
\usepackage{rotating}
\usepackage{mathrsfs}
\usepackage{epsfig}
\usepackage{graphics}
\usepackage{subfigure}

\newcommand{\vertiii}[1]{{\left\vert\kern-0.25ex\left\vert\kern-0.25ex\left\vert #1 
    \right\vert\kern-0.25ex\right\vert\kern-0.25ex\right\vert}}

\theoremstyle{plain}
\newtheorem{theorem}{Theorem}
\newtheorem{lemma}{Lemma}
\newtheorem{conj}{Conjecture}

\newtheorem{corollary}{Corollary}
\newtheorem{proposition}{Proposition}

\begin{document}
\vskip 1cm


\thispagestyle{empty} \vskip 1cm


\title{{
Refining Concentration Bounds for \\Gaussian Quadratic Chaos\\ with
Applications 
}}
\author{ Kamyar Moshksar\\
\small Columbia College\\ Vancouver, BC, Canada} \maketitle

\begin{abstract} 
This paper studies concentration of measure for Gaussian quadratic chaos in the non-asymptotic regime, where existing bounds are improved and several new bounds are proposed. We begin by slightly tightening the classic Hanson-Wright inequality~(HWI). Specifically, we increase the absolute constant in its formula from the largest known value of~0.125 to at least 0.1457 in the symmetric case. We present a sharper version of Laurent-Massart~inequality~(LMI), which increases the absolute constant in HWI from the largest available value of~$1-\frac{\sqrt{3}}{2}\approx 0.134$ to $\frac{9-\sqrt{17}}{32}\approx 0.152$ in the positive-semidefinite case. Stepping beyond HWI, we develop a sequence of inequalities in the symmetric case indexed by $m=1,2,3,\dots, \infty$ that involve Schatten norms of the underlying matrix. The case $m=1$ recovers HWI. These bounds undergo a phase transition in the sense that if the tail parameter is smaller than a critical threshold $\tau_c$, then $m=1$ is the tightest, whereas if it is larger than $\tau_c$, then $m=\infty$ is the tightest. This leads to a novel bound referred to as the~$m_\infty$~bound. All three HWI, LMI and the~$m_\infty$~inequality rely on Markov's inequality. By avoiding Markov's inequality, we introduce the strong $\chi^2$~inequality in the symmetric case and the weak $\chi^2$~inequality in the positive-semidefinite~case.    To further investigate the $m_\infty$ and $\chi^2$~bounds, we restrict our attention to positive-definite matrices and list all available concentration bounds that depend only on the operator norm, namely, the $m_\infty$~bound, the weak $\chi^2$~bound, relaxed versions of HWI and LMI, and the so-called large-deviations bound. The sharpest among these bounds is always either the $m_\infty$~bound or the weak $\chi^2$~bound. If the matrix dimension $n$ is $2,4$~or~$6$, the weak $\chi^2$~bound is tighter than the $m_\infty$~bound. For even~$n\ge8$, the $m_\infty$~bound is sharper than the weak $\chi^2$~bound if and only if the ratio of the tail parameter over the operator norm lies inside an open interval which expands indefinitely as $n$ grows. By utilizing conditioning, modified versions of $m_\infty$, strong $\chi^2$, and HW inequalities of various orders are proposed, all of which can be significantly tighter than the original inequalities. The sharpness of these modified bounds improves dramatically as the order increases, provided that suitable conditions hold. \textcolor{Black}{By leveraging the classical Ruben-Kotz expansion for the distribution function of a positive-definite quadratic form, we introduce the Ruben-Kotz inequalities~(RKIs) of various truncation levels for general symmetric quadratic forms. At the expense of arbitrarily large computational cost, a key advantage of RKIs over alternative inequalities is that they provide an explicit upper bound on the approximation gap, which vanishes under specific conditions as the truncation level increases.}  The paper concludes with two applications.  First, they are applied to change detection in the covariance matrix of a Gaussian vector. Second, they enhance the performance evaluation of unitary space-time modulation over a wireless MISO channel in block Rayleigh flat-fading when both the transmitter and the receiver operate without channel state information.   
\end{abstract}

\section{Introduction}
\label{sec1}
\subsection{Summary of prior art}
This paper focuses on concentration of measure inequalities for Gaussian quadratic forms in the non-asymptotic regime. A well-known result is the Gaussian version of Hanson-Wright inequality~(HWI)~\cite{HW}. It presents an upper bound on the tail probability for a quadratic form in a random vector of independent standard Gaussian random variables known as order 2 chaos.\footnote{In its most general form, HWI concerns quadratic forms in sub-Gaussian random vectors. } For an integer~$n\geq1$ called the dimension throughout the paper, let $\underline{\boldsymbol{x}}$ be such a random vector of length $n$  and $A$ be an arbitrary $n\times n$ matrix. Denote the difference between $\underline{\boldsymbol{x}}^{\mathsf{T}}A\underline{\boldsymbol{x}}$ and its expected value $\mathbb{E}[\underline{\boldsymbol{x}}^{\mathsf{T}}A\underline{\boldsymbol{x}}]=\mathrm{tr}(A)$ by $\boldsymbol{\Delta}$, i.e., 
\begin{align}
\boldsymbol{\Delta}=\underline{\boldsymbol{x}}^{\mathsf{T}}A\underline{\boldsymbol{x}}-\mathrm{tr}(A).
\end{align}
HWI states that for every tail parameter $t>0$,
 \begin{align}
\label{HW}
\Pr(\boldsymbol{\Delta}> t)\leq \exp\Big(-\kappa\min\Big\{\frac{t^2}{\|A\|_2^2},\frac{t}{\|A\|}\Big\}\Big),
\end{align}
where $\|A\|_2$ is the Hilbert-Schmidt norm of $A$, $\|A\|$ is the operator norm of $A$ and $\kappa$ is an absolute constant that does not depend on $n$, $A$ and $t$.\footnote{Replacing $A$ by $-A$ in (\ref{HW}) results in the same upper bound as in (\ref{HW}) on the left-tail probability $\Pr(\boldsymbol{\Delta}<-t)$.} Throughout the paper, we assume $A$ is symmetric without any loss of generality. This can be easily justified by noting that $\underline{\boldsymbol{x}}^{\mathsf{T}}A\underline{\boldsymbol{x}}=\underline{\boldsymbol{x}}^{\mathsf{T}}\tilde{A}\underline{\boldsymbol{x}}$ where $\tilde{A}=\frac{1}{2}(A+A^{\mathsf{T}})$ is symmetric and $\|\tilde{A}\|_2\leq \|A\|_2$ and $\|\tilde{A}\|\leq \|A\|$. These inequalities follow by applying the triangle inequality for norms and noting that $\|A\|_2=\|A^{\mathsf{T}}\|_2$ and $\|A\|=\|A^{\mathsf{T}}\|$. In the original paper~\cite{HW},  $\|A\|$ is the operator norm of the matrix whose entries are absolute values of entries of $A$. Reference~\cite{HW_RV} shows that the inequality indeed holds as stated in (\ref{HW}) with $\|A\|$ being the operator norm of $A$ itself. An explicit value for $\kappa$ is missing or at least hard to identify in \cite{HW,HW_RV}. Theorem~3.1.4 on page~75 in~\cite{christ} provides the bounds 
\begin{align}
\label{christ}
\Pr(\boldsymbol{\Delta}> t)\leq\left\{\begin{array}{cc}
     \exp(-\frac{t^2}{8\|A\|_2^2}) & 0<t\leq \frac{\|A\|_2^2}{\|A\|}    \\
     \exp(\frac{\|A\|_2^2}{8\|A\|^2}-\frac{t}{4\|A\|}) &   t>\frac{\|A\|_2^2}{\|A\|}
\end{array}\right..
\end{align}
For $ t>\frac{\|A\|_2^2}{\|A\|}$, we have $\frac{\|A\|_2^2}{8\|A\|^2}-\frac{t}{4\|A\|}<-\frac{t}{8\|A\|}$. Therefore, (\ref{christ}) implies HWI with $\kappa=\frac{1}{8}=0.125$. More recently, Giraud~\cite{Giraud} gives the explicit value $\kappa=\frac{1}{8}$ at the end of the proof for Theorem~B.8 on page~304.  This value for $\kappa$ can also be inferred form inequality~(3.28) on page~120 in~\cite{GN} given by 
\begin{align}
\label{gn}
\Pr(\boldsymbol{\Delta}> t)\leq \exp\Big(-\frac{t^2}{4(\|A\|_2^2+\|A\|t)}\Big).
\end{align}
It is stated in \cite{GN} that (\ref{gn}) holds for arbitrary symmetric matrix $A$. However, there is an error in the proof. Inequality $\frac{1}{2}(-\ln(1-2x)-2x)\leq \frac{x^2}{1-2x}$ in the middle of page~120 does not hold for $x<0$.  Nonetheless,~(\ref{gn}) certainly holds when $A$ is positive-semidefinite. Loosening the bound by writing $\|A\|_2^2+\|A\|t\leq 2 \max\{\|A\|_2^2,\|A\|t\}$, we arrive at HWI with $\kappa=\frac{1}{8}$. For a positive-semidefinite matrix $A$,  a key result in the literature on concentration of Gaussian quadratic chaos is Inequality~(4.1) on page~1325 in~\cite{LM} due to Laurent and Massart given by
\begin{align}
\label{LM}
\Pr(\boldsymbol{\Delta}>2\|A\|_2\sqrt{t}+2\|A\|t)\leq \exp(-t),
\end{align}  
for all tail parameters $t>0$. We will refer to (\ref{LM}) as Laurent-Massart inequality~(LMI). It is notable that LMI also appears in Example~2.12 in Textbook~\cite{CI_BLM}. LMI is equivalent to\footnote{See Section~\ref{sec6}.} 
\begin{align}
\label{LM2}
\Pr(\boldsymbol{\Delta}>t)\leq \exp\bigg(-\Big(\frac{\sqrt{\|A\|_2^2+2\|A\|t}-\|A\|_2}{2\|A\|}\Big)^2\bigg).
\end{align}
 We will show in Section~\ref{sec6} that LMI implies HWI at best with $\kappa=1-\frac{\sqrt{3}}{2}\approx 0.134$. 
 
 Throughout this paper, the term ``non-asymptotic regime'' means that the rank of the matrix \(A\), denoted by \(r(A)\), remains finite and does not grow large. Due to the rotational invariance of Gaussian vectors, \(\Pr(\boldsymbol{\Delta} > t)\) depends on \(A\) only through its nonzero eigenvalues. While \(n\) may grow to infinity, the number of these nonzero eigenvalues, which is \(r(A)\), can remain bounded and hence, we are still in the non-asymptotic regime. 

 A natural question is whether the Cumulative Distribution Function (CDF) of \(\underline{\boldsymbol{x}}^{\mathsf{T}}A\underline{\boldsymbol{x}}\), denoted by $F_{\underline{\boldsymbol{x}}^{\mathsf{T}}A\underline{\boldsymbol{x}}}(\cdot)$, can be found in closed form. This problem was studied in several pioneering papers from the 1940s through the 1960s. Specifically, Robbins~\cite{Robbins1}, Gurland~\cite{Gurland} and Ruben~\cite{Ruben}  proposed a power series expansion, a Laguerre series expansion and an expansion in \(\chi ^{2}\)~CDFs, respectively, for positive-definite $A$. A special case of Ruben's result was reported earlier in~\cite{Robbins2} by Robbins and Pitman.  All these  expansions were later unified as special cases of a general expansion by Kotz et al. in~\cite{Kotz1}. Mathai and Provost presented these findings along with further extensions in their comprehensive book~\cite{Mathai}. Assuming the positive-definite matrix $A$ has eigenvalues $\lambda_1,\cdots, \lambda_n>0$, the Ruben-Kotz expansion in $\chi^2$ CDFs is given by 
 \begin{align}
 \label{RK_expansion}
F_{\underline{\boldsymbol{x}}^{\mathsf{T}}A\underline{\boldsymbol{x}}}(z)=\sum_{l=0}^\infty c^{\text{RK}}_{l}F_{n+2l}\Big(\frac{z}{s}\Big),
\end{align}
where $F_{n+2l}(\cdot)=F_{\chi^2_{n+2l}}(\cdot)$ is the CDF of a chi-squared random variable with $n+2l$ degrees of freedom, and $s>0$ is an arbitrary scale factor. The coefficients $c^{\text{RK}}_l=c^{\text{RK}}_l(\lambda_1,\cdots, \lambda_n; s)$ are computed via the recursion 
\begin{align}
\label{RK_recursion}
c_0^{\text{RK}}=\prod_{i=1}^n\Big(\frac{s}{\lambda_i}\Big)^{\frac{1}{2}},\quad c^{\text{RK}}_l=\frac{1}{2l}\sum_{j=0}^{l-1}d^{\text{RK}}_{l-j}c^{\text{RK}}_j,\,\,\, l\geq 1,
\end{align}   
where
\begin{align}
d^{\text{RK}}_l=\sum_{i=1}^n\Big(1-\frac{s}{\lambda_i}\Big)^l,\,\,\, l\geq1.
\end{align}
If $0<s\le\lambda_{\min}(A)= \min_{1\le i\le n}\lambda_i$, the coefficients $c^{\text{RK}}_{l}$ are strictly positive. Under this condition, truncating the series in (\ref{RK_expansion}) at an integer $L\geq 0$ yields the concentration bound
\begin{align}
\label{RK_CB}
\Pr(\boldsymbol{\Delta}>t)=1-F_{\underline{\boldsymbol{x}}^{\mathsf{T}}A\underline{\boldsymbol{x}}}\big(t+\mathrm{tr}(A)\big)\leq 1-\sum_{l=0}^L c^{\text{RK}}_{l}F_{n+2l}\Big(\frac{t+\mathrm{tr}(A)}{s}\Big).
\end{align} 
 A key advantage of this formulation is that the truncation error can be explicitly upper-bounded by 
\begin{align}
\label{error_cap}
  \sum_{l=L+1}^\infty c^{\text{RK}}_{l}F_{n+2l}\Big(\frac{t+\mathrm{tr}(A)}{s}\Big)\leq \sum_{l=L+1}^\infty c^{\text{RK}}_{l}=1-\sum_{l=0}^L c^{\text{RK}}_{l},
\end{align}
where the last step is due to $\sum_{l=0}^\infty c^{\text{RK}}_{l}=1$. 

\textcolor{Black}{For the general symmetric case, Gurland~\cite{Gurland} derived a single-series formula for $F_{\underline{\boldsymbol{x}}^{\mathsf{T}}A\underline{\boldsymbol{x}}}(z)$ under the assumption that the number of positive eigenvalues of $A$ is even when $z>0$ and the number of negative eigenvalues of $A$ is even when $z\le0$. The components in this series are expressed using Laguerre-type polynomials and the incomplete Gamma function.  More recently, Provost and Rudiuk~\cite{Provost} derived highly intricate, explicit formulas for the CDF of $\underline{\boldsymbol{x}}^{\mathsf{T}}A\underline{\boldsymbol{x}}$ as a triple infinite series. The terms of their series incorporate the incomplete Gamma function, as well as Euler's $\psi$ function also known as the digamma function.} 

 Next, we present a detailed account of our main contributions, which include improving existing inequalities, proposing several new inequalities, and exploring their applications to two problems in statistical signal processing and wireless communications.
\subsection{Contributions and main results}
\label{1B}
The paper offers seven main contributions: 
\begin{enumerate}
  \item We slightly modify the proof given in~\cite{Giraud} for HWI. In the course of proof, Reference~\cite{Giraud} uses the bound 
  \begin{align}
-\ln(1-x)\leq x+x^2,\,\,\,  |x|\le\frac{1}{2},
\end{align} 
which ultimately results in $\kappa=0.125$. Instead, we consider
  \begin{align}
  \label{theone}
-\ln(1-x)\leq x+ax^2,\,\,\,  |x|\le b,
\end{align} 
for $a,b>0$. For given $0<b<1$, we determine the smallest $a$ (tightest bound) such that~(\ref{theone})~holds.  As a result, we are able to increase $\kappa$ from $0.125$ to at least $0.1457$ for arbitrary symmetric matrix~$A$. More precisely, we have the following proposition. 
 \begin{proposition}
 \label{prop_1}
  For $0< b<1$  and integer $m\geq1$ define the sequence of functions  
  \begin{align}
  \label{theta_func}
\theta_m(b)=-\frac{1}{b^{m+1}}\Big(\ln(1-b)+b+\frac{b^2}{2}+\frac{b^3}{3}+\cdots+\frac{b^m}{m}\Big).
\end{align}  
If $A$ is symmetric, then HWI holds with 
\begin{align}
\label{HW_constant}
\kappa=\frac{b^*}{4}\approx 0.1457, 
\end{align}
 where $b^*$ is the unique solution for $b$ in the equation\footnote{The equation $2b\theta_1(b)=1$ simplifies to $\ln(1-b)=-\frac{3}{2}b$.  } $2b\theta_1(b)=1$. 
  \end{proposition}  
  \begin{proof}
  See Section~\ref{sec4}.
  \end{proof}
 A special case of Lemma~\ref{lem_1} in Section~\ref{sec4} shows that the smallest value for the parameter~$a$ in (\ref{theone}) is~$\theta_1(b)$. This optimal value for $a$ can also be inferred from Corollary~4 on Page~5 of the independent work cited in~\cite{Gallagher}. 
The functions $\theta_m(b)$ will appear again in Proposition~\ref{prop_4} in below.
\item In the positive-semidefinite case, we present a sharper version of LMI. To prove (\ref{LM}), Reference~\cite{LM} relies on the bound
  \begin{align}
-\ln(1-x)\leq x+\frac{x^2}{2(1-x)},\,\,\, 0\le x< 1.
\end{align}
 Instead, we consider  
\begin{align}
\label{thetwo}
-\ln(1-x)\leq x+\frac{x^2}{2(1-ax)},\,\,\, 0\le x\le b,
\end{align}
for $0<a, b<1$. We show that for every $\frac{2}{3}<a<1$, this inequality holds for $b$ at least as large as~$\frac{3a-2}{a(2a-1)}$. As a result, we derive an improved LMI presented in the next proposition.
\begin{proposition}
\label{prop_2}
Let $A$ be positive-semidefinite and 
\begin{align}
\label{conv}
\alpha=\|A\|,\quad \beta=\|A\|_2^2.
\end{align}
For $\frac{2}{3}<a<1$, let 
\begin{align}
\label{toy}
b=\frac{3a-2}{a(2a-1)},\quad c=\frac{3a-2}{2(1-a)^2}.
\end{align}
 Then 
\begin{align}
\label{en_LM}
\Pr(\boldsymbol{\Delta}>t)\leq e^{\Lambda(t, a)},
\end{align}
where 
\begin{align}
\label{en_exp}
\Lambda(t, a)=\left\{\begin{array}{cc}
    -\frac{t}{2a\alpha}+\frac{\beta}{2a^2\alpha^2}\Big(\sqrt{1+\frac{2a\alpha t}{\beta}}-1\Big)  & \frac{\alpha t}{\beta}\leq c   \\
    -\frac{bt}{2\alpha}+\frac{b^2\beta}{4(1-ab)\alpha^2}  &  \frac{\alpha t}{\beta}> c 
\end{array}\right..
\end{align}
\end{proposition}
\begin{proof}
See Section~\ref{sec5}.
\end{proof}
Letting $a\to1^-$ in (\ref{en_LM}) recovers LMI in (\ref{LM2}). In Subsection~\ref{2A}, we minimize the upper bound on the right-hand side of (\ref{en_LM}) in terms of $a$ and compare the resulting optimized inequality with LMI to demonstrate the improvement. We will see that the optimal value for $a$ is the unique root inside the interval $(\frac{2}{3},1)$ of a quintic polynomial equation\footnote{See Equation~(\ref{quintic}).} and it depends on the matrix~$A$ and the tail parameter~$t$ only through the ratio $\rho=\frac{t\|A\|}{\|A\|_2^2}$. We can approximate this unique root using an analytic formula\footnote{See Equation~(\ref{a_star}).} with an approximation error shown to be less than $0.035$ for all values of $A$~and~$t$. These investigations lead to two variants of LMI named optimal LMI~(OLMI) and augmented LMI~(ALMI). As a consequence of Proposition~\ref{prop_2},  we offer the next result which is a stronger version of Proposition~\ref{prop_1} in the positive-semidefinite case.  
\begin{proposition}
\label{prop_3}
Let $A$ be positive-semidefinite. The inequality in (\ref{en_LM}) implies HWI at best with
\begin{align}
\label{ }
\kappa=\frac{9-\sqrt{17}}{32}\approx 0.1524,
\end{align}
which is achieved for
\begin{align}
a=\frac{7-\sqrt{17}}{4}\approx 0.7192.
\end{align}
 Moreover, LMI in (\ref{LM}) implies HWI at best with $\kappa=1-\frac{\sqrt{3}}{2}\approx 0.134$.
\end{proposition}
\begin{proof}
See Section~\ref{sec6}. 
\end{proof}
  \item We explore beyond HWI for arbitrary symmetric matrix $A$ by considering an extension of~(\ref{theone}), i.e., an inequality of the form 
  \begin{align}
  \label{gen_ineq1}
-\ln(1-x)\leq x+\frac{x^2}{2}+\frac{x^3}{3}+\cdots+\frac{x^m}{m}+a|x|^{m+1},\,\,\, |x|\le b,
\end{align}
for a given integer $m\geq1$ and $a,b>0$. We show\footnote{See Lemma~\ref{lem_1} in Section~\ref{sec4}.} that for every $m\ge1$ and $0<b<1$, the smallest~$a$~(tightest upper bound) for which (\ref{gen_ineq1}) holds is $\theta_m(b)$ as defined in (\ref{theta_func}). 
 Consequently, we obtain a sequence of concentration bounds for Gaussian quadratic chaos that are written in terms of Schatten norms of~$A$.\footnote{See Subsection~I.C for the definition of Schatten $p$-norm.}  
 \begin{proposition}
 \label{prop_4}
 Let $m\geq1$ be an integer, $0<b<1$, $\|A\|_{m+1}$ be Schatten norm of order $m+1$ for a symmetric matrix $A$ and $r(A)$ be the rank of $A$. Then  
 \begin{align}
 \label{HW_general}
\Pr(\boldsymbol{\Delta}>t)\leq e^{\Lambda_m(t, b)},
\end{align} 
 where\footnote{The function $\Lambda(t,a)$ in (\ref{en_exp}) and the functions $\Lambda_m(t,b)$ in (\ref{Lambda_m}) are not related. }
 \begin{align}
 \label{Lambda_m}
\Lambda_m(t, b)=\frac{r(A)}{2}\sum_{k=2}^m \frac{b^k}{k}-\kappa_m(b)\min\bigg\{\frac{t^{1+\frac{1}{m}}}{\|A\|_{m+1}^{1+\frac{1}{m}}}, \frac{t}{\|A\|} \bigg\} 
\end{align}
and 
\begin{align}
\label{kappa_m}
\kappa_m(b)=\frac{m}{2(m+1)}\min\bigg\{\frac{1}{\big((m+1)\theta_m(b)\big)^{\frac{1}{m}}}, b\bigg\}
\end{align}
with $\theta_m(b)$ as defined in (\ref{theta_func}). 
 \end{proposition}
 \begin{proof}
 See Section~\ref{sec7}.
 \end{proof}
Choosing $m=1$ in (\ref{HW_general}) recovers HWI in (\ref{HW})  with $\kappa=\kappa_1(b)$. Schatten norm $\|A\|_{m+1}$ should not be mistaken with the $L_p$-norm of $A$ with $p=m+1$ unless $m=1$ in which case $\|A\|_2$ is Hilbert-Schmidt norm of $A$. The upper bound $e^{\Lambda_m(t, b)}$ in (\ref{HW_general}) is to be minimized over $0<b<1$ and the choice of $m\geq1$. In general, the sequence of upper bounds $\inf_{b\in (0,1)}e^{\Lambda_m(t, b)}$ for $m=1,2,3,\dots$ is not monotone in~$m$ for given $A$ and $t$. However, we conjecture that it undergoes a phase transition as stated~next. 
\begin{conj}
\label{conj_1}
For a symmetric matrix $A$, there exists a critical threshold $\tau_c\ge0$ such that if $t<\tau_c$, then $\inf_{b\in (0,1)}e^{\Lambda_m(t, b)}$ is the smallest for $m=1$ and  if $t>\tau_c$, then $\inf_{b\in (0,1)}e^{\Lambda_m(t, b)}$ is the smallest when $m$ grows to infinity.\footnote{More precisely, $\lim_{m\to\infty}\inf_{b\in (0,1)}e^{\Lambda_m(t, b)}$ exists and it is smaller than $\inf_{b\in (0,1)}e^{\Lambda_m(t, b)}$ for every $m\geq1$.} If every eigenvalue of $A$ has the same absolute value, then $\tau_c=0$. 
\end{conj}
If Conjecture~\ref{conj_1} is true, then one must always either choose $m=1$ or let $m$ grow to infinity. In Subsection~\ref{2B}, we will study the bounds in (\ref{HW_general}) both through a numerical example and analytically where two corollaries to Proposition~\ref{prop_4} are stated.  The first corollary looks into the limiting value of~(\ref{HW_general}) as $m$ grows large. This leads to the simple bound 
  \begin{align}
 \label{infty_m}
\Pr(\boldsymbol{\Delta}>t)\leq \Big(1+\frac{t}{r(A)\|A\|}\Big)^{\frac{r(A)}{2}}e^{-\frac{t}{2\|A\|}},
\end{align}
which we refer to as the $m_\infty$~bound or $m_\infty$~inequality. The right hand side in (\ref{infty_m}) depends on $A$ only through its rank $r(A)$ and its operator norm $\|A\|$. The second corollary presents an estimate (upper bound) for the critical threshold $\tau_c$ mentioned in Conjecture~1. 

\item To further explore the $m_\infty$~bound, we raise the following question in Subsection~\ref{2D}: 

\textit{Let $A$ be an $n\times n$ positive-definite matrix. How much is $\sup_{A\in\mathbb{S}^n_{++}: \|A\|=\nu}\Pr(\boldsymbol{\Delta}>t)$ in terms of $n, \nu, t$?  Here, $\mathbb{S}^n_{++}$ is the cone of $n\times n$ positive-definite matrices. }

We present five upper bounds on $\Pr(\boldsymbol{\Delta}>t)$ that depend on $A$ only through $\|A\|$. One candidate is the~$m_\infty$~bound in (\ref{infty_m}) with $r(A)=n$. The second and third candidates are relaxed versions of HWI and LMI, respectively, which do not depend on $\|A\|_2$ anymore. The relaxed HWI~is 
\begin{align}
\label{relaxed_HW}
\Pr(\boldsymbol{\Delta}>t)\leq e^{-\kappa\min\big\{\frac{t^2}{n\|A\|^2}, \frac{t}{\|A\|}\big\}},\quad \kappa=\frac{9-\sqrt{17}}{32}
\end{align}
and the relaxed LMI is given by
\begin{align}
\label{relaxed_LM}
\Pr(\boldsymbol{\Delta}>t)\leq e^{-\frac{n}{4}\big(\sqrt{1+\frac{2t}{n\|A\|}}-1\big)^2}.
\end{align}
The fourth candidate is based on the chi-squared distribution which we refer to as the weak $\chi^2$~bound. It is given~by 
\begin{align}
\label{cvb}
\Pr(\boldsymbol{\Delta}>t)\leq 1-F_{n}\Big(1+\frac{t}{\|A\|}\Big),
\end{align}
where $F_n(\cdot)=F_{\chi^2_n}(\cdot)$ is the CDF for a chi-squared random variable with $n$ degrees of freedom. The weak $\chi^2$~bound is a relaxed version of the strong $\chi^2$~bound that holds for arbitrary symmetric matrix $A$. It is given~by 
\begin{align}
\label{chi2_bound}
\Pr(\boldsymbol{\Delta}>t)\leq\left\{\begin{array}{cc}
    1-F_{r(A)}(\frac{t+\mathrm{tr}(A)}{\lambda^*_{\max}(A)})  &  \lambda^*_{\max}(A)>0  \\\\
    F_{r(A)}(\frac{t+\mathrm{tr}(A)}{\lambda^*_{\max}(A)})  &  \lambda^*_{\max}(A)<0
\end{array}\right.,
\end{align}
where $r(A)$ is the rank of $A$, $\mathrm{tr}(A)$ is the trace of $A$ and $\lambda^*_{\max}(A)$ is the largest~\textit{nonzero} eigenvalue~of~$A$.  The fifth and final candidate is a further relaxed version of the weak $\chi^2$~bound, known as the large-deviations bound, given by 
\begin{align}
\label{LDB}
\Pr(\boldsymbol{\Delta}>t)<\frac{1}{\sqrt{e}}\Big(\frac{e(1+\frac{t}{\|A\|})}{n}\Big)^{\frac{n}{2}}e^{-\frac{t}{2\|A\|}},\,\,\, \frac{t}{\|A\|}\geq n-1.
\end{align} 
We emphasize that (\ref{LDB}) is valid only for $\frac{t}{\|A\|}\geq n-1$. Derivations of the relaxed HWI,  the relaxed LMI, the strong and weak $\chi^2$ bounds and the large-deviations bound are presented in Subsection~\ref{2D}.   Note that the bounds in (\ref{infty_m}), (\ref{relaxed_HW}), (\ref{relaxed_LM}), (\ref{cvb}) and (\ref{LDB}) depend on $A$ and $t$ only through the ratio $\frac{t}{\|A\|}$. Our first result is the following.  
\begin{proposition}
\label{prop_6}
The $m_\infty$~bound in~(\ref{infty_m}) is sharper than the relaxed HWI in~(\ref{relaxed_HW}), the relaxed LMI in~(\ref{relaxed_LM}) and the large-deviations bound in (\ref{LDB}) for every positive-definite matrix $A$ and $t>0$.
\end{proposition}
\begin{proof}
See Section~\ref{sec9}.
\end{proof}
It follows from Proposition~\ref{prop_6} that among the five candidates we have introduced, the sharpest bound is always either the $m_\infty$~bound or the weak $\chi^2$~bound. Our second result compares these two finalists when the dimension $n$ is even.\footnote{$F_{\chi^2_n}(\cdot)$ only admits a closed-form expression for even values of $n$. If $n$ is odd, one may rely on a lower bound on $F_{\chi^2_n}(\cdot)$ such as $F_{\chi^2_{n+1}}(\cdot)$.} 
\begin{proposition}
\label{prop_7}
 If $n=2,4,6$, then the weak $\chi^2$~bound is sharper than the $m_\infty$~bound regardless of the positive-definite matrix~$A$~and~$t$. If $n$ is an even integer greater than or equal to $8$, there exist positive constants $r_n$ and $r'_n$ such that $r_n<1<r'_n$ and the $m_\infty$~bound is sharper than the weak~$\chi^2$~bound if and only if $r_n<\frac{t}{\|A\|}<r'_n$.
\end{proposition}
\begin{proof}
The proof uses Descartes' rule of signs. See Section~\ref{sec10}.
\end{proof}
Additionally, we prove in Appendix~\ref{A4} that as the even integer $n$ increases, $r_n$ goes to zero and $r'_n$ grows to infinity, i.e., 
\begin{align}
\label{final_claim}
\lim_{n\to\infty}r_n=0,\quad \lim_{n\to\infty}r'_n=\infty.
\end{align}
Thus, the advantage of the $m_\infty$~bound over the weak $\chi^2$~bound is more significant as $n$ grows larger. 

\item We propose modified versions of $m_{\infty }$, strong $\chi ^{2}$, and HW inequalities of various orders. We denote the order of these modified bounds by an integer $k \geq 1$, treating original bounds as order zero. Under suitable conditions, the proposed modifications offer significant improvements over their earlier versions.  They are particularly effective when the symmetric matrix $A$ possesses one or more negative eigenvalues whose absolute values are sizeably larger than the others. However, these improvements are subject to specific viability conditions, which must be verified for the bounds to be applicable. There is a distinct trade-off between bound tightness and computational complexity:
\begin{itemize}
  \item Computational Cost: The modified $m_{\infty }$ and strong $\chi ^{2}$ bounds require more computing time than the modified HW bounds of the same order.
  \item Accuracy vs. Complexity: While increasing the order $k$ typically enhances the effectiveness of the bound, it does so at the expense of higher computational cost.
\end{itemize}
 Nonetheless, a practical advantage of the modified bounds is that they do not involve special functions or numerical integration.  Let us denote the eigenvalues for the symmetric matrix $A$ by $\lambda_1,\dots, \lambda_n$. Let $1\leq k\leq n-1$ and   \begin{align}
\label{poloi_1}
\lambda_1\le\cdots\le\lambda_k<0
\end{align}
be the $k$ largest negative\footnote{By larger negative, we mean a negative number with a larger absolute value.} eigenvalues. Define the parameters  $a_k$ and $b_k$ by 
\begin{align}
\label{a_k_b_k}
a_k=\max_{k+1\leq i\leq n}|\lambda_i|,\quad b_k=\sum_{i=1}^k\lambda_i.
\end{align}
The next proposition presents the $m_\infty$~inequality of order $k$.
\begin{proposition}
\label{prop_17}
Let $r(A)\geq 3$, $\lambda_1,\dots, \lambda_k$ be as in (\ref{poloi_1}) and $r(A)-k$ be even. If 
\begin{align}
\label{poloi_2}
t+b_k\ge0,
\end{align}
then 
\begin{align}
\label{call_my_name}
&\Pr(\boldsymbol{\Delta}>t)\leq e^{-\frac{t+b_k}{2a_k}}\sum_{q=0}^{\frac{r(A)-k}{2}}\frac{c_{k,q}}{(\frac{r(A)-k}{2}-q)!}\Big(1+\frac{t+b_k}{(r(A)-k)a_k}\Big)^{\frac{r(A)-k}{2}-q},
\end{align}
where
\begin{align}
\label{def_c_l}
c_{k,q}=\frac{\big(\frac{r(A)-k}{2}\big)!}{\prod_{i=1}^k\big(1-\frac{\lambda_i}{a_k}\big)^{\frac{1}{2}}} \sum_{\substack{j_1,\dots, j_k\ge0\\ j_1+\cdots+j_k=q}}\prod_{i=1}^k\frac{(2j_i-1)!!}{j_i!}\phi_i^{j_i},
\end{align}
\begin{align}
\label{def_phi_i}
   \phi_i=\frac{\lambda_i}{(r(A)-k)(\lambda_i-a_k)},
\end{align}
$a_k, b_k$ are given in (\ref{a_k_b_k}) and the double factorial notation is given by
\begin{align}
\label{double_factorial}
(2j-1)!!=\left\{\begin{array}{cc}
    1  &  j=0  \\
    1\times 3\times 5\times\cdots\times (2j-1)  &   j\geq 1
\end{array}\right..
\end{align}
\end{proposition}
\begin{proof}
See Section~\ref{sec11}.
\end{proof}
In the special case where $A$ is full-rank,~$n$ is odd and $k=1$, we have $a_1=\max_{2\leq i\leq n}|\lambda_i|$, $b_1=\lambda_1=\lambda_{\min}(A)$ and (\ref{call_my_name}) simplifies~to 
\begin{align}
\label{case_k_1}
\Pr(\boldsymbol{\Delta}>t)\leq \frac{e^{-\frac{t+\lambda_1}{2a_1}}}{\big(1-\frac{\lambda_1}{a_1}\big)^{\frac{1}{2}}}\sum_{q=0}^{\frac{n-1}{2}}(2q-1)!!{\frac{n-1}{2}\choose q}\Big(\frac{\lambda_1}{(n-1)(\lambda_1-a_1)}\Big)^q\Big(1+\frac{t+\lambda_1}{(n-1)a_1}\Big)^{\frac{n-1}{2}-q}. 
\end{align}
In Subsection~\ref{2E}, we observe that for every $k\ge2$, the $m_\infty$~bound of order $k$ is sharper than the $m_\infty$~bound of order $k-2$ for all \textit{sufficiently large} values of $t$. In order to compute the bound in~(\ref{call_my_name}), one needs to find all solutions to $j_1+\cdots+j_k=q$ in nonnegative integers $j_1,\dots, j_k$ for $0\leq q\leq \frac{r(A)-k}{2}$. The number of these solutions is exactly the number of solutions to $j_0+j_1+\cdots+j_k=\frac{r(A)-k}{2}$ in nonnegative integers $j_0,j_1,\dots, j_k$, which is equal to ${\frac{r(A)-k}{2}+k+1-1\choose k+1-1}={\frac{r(A)+k}{2}\choose k}$ by the Stars and Bars Theorem. Generating these solutions can be quite time-consuming, even for moderate values of $r(A)$ and $k$. The numbers $c_{k,q}$ in (\ref{def_c_l}) can be computed efficiently via a recursion similar to the one in~(\ref{RK_recursion}), as given by the next proposition.  
\begin{proposition}
\label{prop6.1}
The numbers $c_{k,q}$ in (\ref{def_c_l}) satisfy the recursion 
\begin{align}
\label{Pill_000}
c_{k,0}=\frac{\big(\frac{r(A)-k}{2}\big)!}{\prod_{i=1}^k\big(1-\frac{\lambda_i}{a_k}\big)^{\frac{1}{2}}},\quad c_{k,q}=\frac{1}{2q}\sum_{j=0}^{q-1}d_{k,q-j}c_{k,j},\,\,\, q\ge1,
\end{align}
where 
\begin{align}
d_{k,q}=\sum_{i=1}^k (2\phi_i)^q,\,\,\, q\geq1
\end{align}
and $\phi_i$ are given in (\ref{def_phi_i}).
\end{proposition}
\begin{proof}
See Section~\ref{sec_123321}.
\end{proof}

Next, we present the strong $\chi^2$~inequality of order $k$. 
\begin{proposition}
\label{prop_9}
Let $r(A)\geq 3$, $\lambda_{\max}(A)>0$, $\lambda_1,\dots, \lambda_k$ be as in (\ref{poloi_1}) and $r(A)-k$ be even. If
\begin{align}
\label{poloi_4}
t+\mathrm{tr}(A)\ge0,
\end{align}
then 
\begin{align}
 \label{modified_chi2_final}
\Pr(\boldsymbol{\Delta}>t)\leq e^{-\frac{t+\mathrm{tr}(A)}{2\lambda_{\max}(A)}}\sum_{m=0}^{\frac{r(A)-k}{2}-1}\frac{1}{2^m}\sum_{q=0}^m\frac{\tilde{c}_{k,q}}{(m-q)!}\Big(\frac{t+\mathrm{tr}(A)}{\lambda_{\max}(A)}\Big)^{m-q},
\end{align}
where 
\begin{align}
\label{tilde_c_kl}
\tilde{c}_{k,q}=\frac{1}{\prod_{i=1}^k\big(1-\frac{\lambda_i}{\lambda_{\max}(A)}\big)^{\frac{1}{2}}} \sum_{\substack{j_1,\dots, j_k\ge0\\ j_1+\cdots+j_k=q}}\prod_{i=1}^k\frac{(2j_i-1)!!}{j_i!}\tilde{\phi}_i^{j_i},
\end{align}
\begin{align}
\label{}
   \tilde{\phi}_i=\frac{\lambda_i}{\lambda_i-\lambda_{\max}(A)}
\end{align}
and the double factorial notation is given in (\ref{double_factorial}). Moreover, the numbers $\tilde{c}_{k,q}$ can be computed via the recursion 
\begin{align}
\label{pill_111}
\tilde{c}_{k,0}=\frac{1}{\prod_{i=1}^k\big(1-\frac{\lambda_i}{\lambda_{\max}(A)}\big)^{\frac{1}{2}}} ,\quad \tilde{c}_{k,q}=\frac{1}{2q}\sum_{j=0}^{q-1}\tilde{d}_{k,q-j}\tilde{c}_{k,j},\,\,\, q\ge1,
\end{align}
where 
\begin{align}
\tilde{d}_{k,q}=\sum_{i=1}^k (2\tilde{\phi}_i)^q,\,\,\, q\geq1.
\end{align}
\end{proposition}
\begin{proof}
See Appendix~\ref{A5}. 
\end{proof}
If $A$ is full-rank, $n$ is odd and $k=1$, the bound in (\ref{modified_chi2_final}) simplifies~to 
\begin{align}
\Pr(\boldsymbol{\Delta}>t)\leq \frac{e^{-\frac{t+\mathrm{tr}(A)}{2\lambda_{\max}(A)}}}{(1-\frac{\lambda_1}{\lambda_{\max}(A)})^{\frac{1}{2}}}\sum_{m=0}^{\frac{n-3}{2}}\frac{1}{2^mm!}\sum_{q=0}^m(2q-1)!!{m\choose q}\Big(\frac{\lambda_1}{\lambda_1-\lambda_{\max}(A)}\Big)^q\Big(\frac{t+\mathrm{tr}(A)}{\lambda_{\max}(A)}\Big)^{m-q}.
\end{align} 
Subsection~\ref{2E} shows that as the order \(k\) increases, the modified strong \(\chi ^{2}\) inequality provides a sharper bound for every deviation parameter \(t>0\). Computing the bound in (\ref{modified_chi2_final}) via (\ref{tilde_c_kl}) requires finding all nonnegative integer solutions to $j_1+\cdots+j_k\leq \frac{r(A)-k}{2}-1$. Finding these ${\frac{r(A)+k}{2}-1\choose k}$ solutions can be computationally expensive. The recursion in (\ref{pill_111}) offers an efficient alternative in computing $\tilde{c}_{k,q}$. 

Our last proposition presents HWI of order $k$. 
\begin{proposition}
\label{prop_8}
Let $\lambda_1,\dots, \lambda_k$ be as in (\ref{poloi_1}). If
\begin{align}
\label{poloi_3}
t+b_k\ge\frac{\sum_{i=k+1}^n \lambda_i^2}{a_k},
\end{align}
then 
\begin{align}
\label{mod_HWI}
\Pr(\boldsymbol{\Delta}>t)\leq e^{-\frac{\kappa(t+b_k)}{a_k}}\prod_{i=1}^k \Big(1-2\kappa\frac{ \lambda_i}{a_k}\Big)^{-\frac{1}{2}},
\end{align}
where $a_k, b_k$ are given in (\ref{a_k_b_k}) and $\kappa$ is given in (\ref{HW_constant}). 
\end{proposition}  
\begin{proof}
See Appendix~\ref{A6}. 
\end{proof}
 In contrast to $m_\infty$ and strong $\chi^2$ inequalities of order $k$, Proposition~\ref{prop_8} does not require $r(A)-k$ be even. Note that the condition in (\ref{poloi_3}) is stronger than the one in (\ref{poloi_2}). Therefore, a modified HWI is viable for larger values of $t$ compared to the modified $m_\infty$~inequality of the same order. It is easy to see that if\footnote{It is understood that $a_0=\|A\|$.} $a_k=a_{k-1}$ for some $k\geq1$, then HWI of order $k$ is \textit{looser} than HWI of order $k-1$ for \textit{every}~$t>0$. If $a_k<a_{k-1}$, then HWI of order $k$ is sharper than HWI of order $k-1$ for \textit{sufficiently large} values of $t$. Moreover, $\frac{\sum_{i=k+1}^n \lambda_i^2}{a_k}-b_k$ is non-decreasing in $k$. Hence, the set of values of $t$ for which the liability condition in (\ref{poloi_3}) holds either remains unchanged or shrinks as the order~$k$ increases. We explain these points further in Subsection~\ref{2E}.   
  
\textcolor{Black}{As mentioned earlier, these modified inequalities are most effective when one or more negative eigenvalues of $A$ is larger in magnitude than the remaining eigenvalues. This pattern occurs in the two applications discussed below, where the positive eigenvalues are always less than 1, whereas some negative eigenvalues have a much larger magnitude.}  

\item \textcolor{Black}{By leveraging the classical Ruben-Kotz expansion for the CDF of a positive-definite quadratic form in~(\ref{RK_expansion}) and allowing for arbitrarily large computational cost, we introduce Ruben-Kotz inequalities~(RKIs) of various truncation levels $L\ge0$ for a symmetric quadratic form. Before we present our main result in this subsection, some notation and definitions are in order. 
\begin{enumerate}
  \item  Let the symmetric matrix $A$ have $u$ negative eigenvalues and $v=r(A)-u$ positive eigenvalues, denoted by $\lambda_1,\dots, \lambda_{u}<0$ and $\lambda_{u+1},\dots, \lambda_{r(A)}>0$, respectively. We define 
\begin{align}
\label{}
 \mathcal{S}^-=\{-\lambda_1,\dots, -\lambda_{u}\},\,\,\, \mathcal{S}^+=\{\lambda_{u+1},,\dots, \lambda_{r(A)}\}
\end{align}
and 
\begin{align}
\label{}
\check{u}=\left\{\begin{array}{cc}
    u  &  \textrm{$u$ is even}  \\
      u-1  &  \textrm{$u$ is odd}    
\end{array}\right.,\quad \hat{v}=\left\{\begin{array}{cc}
    v  &  \textrm{$v$ is even}  \\
      v+1  & \textrm{$v$ is odd}    
\end{array}\right..
\end{align}
Here, $\mathcal{S}^-$ is the set of absolute values of all negative eigenvalues, $\mathcal{S}^+$ is the set  of all positive eigenvalues, $\check{u}$ is the largest even integer less than or equal to $u$ and $\hat{v}$ is the smallest even integer greater than or equal to $v$. 
  \item Let $ c_l^{\text{RK}^-}$ and $ c_l^{\text{RK}^+}$ be the Ruben-Kotz coefficients given by  
\begin{align}
\label{}
   c_l^{\text{RK}^-}=c_l^{\text{RK}}(-\lambda_1,\cdots, -\lambda_{u};s^-),\,\,\,l\geq0
\end{align}
and 
\begin{align}
\label{}
  \quad  c_l^{\text{RK}^+}=c_l^{\text{RK}}(\lambda_{u+1},\cdots, \lambda_{r(A)};s^+),\,\,\,l\ge0,
\end{align}
as defined in (\ref{RK_recursion}), where the scale factors $s^-$ and $s^+$ satisfy 
\begin{align}
\label{s_pm_ineq}
  0< s^-\le\min_{\lambda\in  \mathcal{S}^-} \lambda,\quad  0< s^+\le\min_{\lambda \in \mathcal{S}^+} \lambda.
\end{align}
  \item For an integer $q\geq 0$, we define the numbers $c^-_q$ and $c^+_q$  recursively as 
\begin{align}
\label{c_neg_q}
c^-_{0}=\frac{1}{\prod_{\lambda\in \mathcal{S}^+}\big(1+\frac{\lambda}{s^-}\big)^{\frac{1}{2}}} ,\quad c^-_{q}=\frac{1}{2q}\sum_{j=0}^{q-1}d^-_{q-j}c^-_{j},\,\,\, q\ge1,
\end{align}
where 
\begin{align}
d^-_{q}=\sum_{\lambda\in \mathcal{S}^+} \Big(\frac{2\lambda}{\lambda+s^-}\Big)^q,\,\,\, q\geq1
\end{align}
and 
\begin{align}
\label{c_plus_q}
c^+_{0}=\frac{1}{\prod_{\lambda\in \mathcal{S}^-}\big(1+\frac{\lambda}{s^+}\big)^{\frac{1}{2}}} ,\quad c^+_{q}=\frac{1}{2q}\sum_{j=0}^{q-1}d^+_{q-j}c^+_{j},\,\,\, q\ge1,
\end{align}
where 
\begin{align}
d^+_{q}=\sum_{\lambda\in \mathcal{S}^-} \Big(\frac{2\lambda}{\lambda+s^+}\Big)^q,\,\,\, q\geq1.
\end{align}
\item For $k\ge1$, let 
\begin{align}
\label{KSD}
   \delta_k=\sup_{x\in \mathbb{R}}|F_k(x)-F_{k+1}(x)|
\end{align}
be the Kolmogorov-Smirnov distance between $F_k(\cdot)$ and $F_{k+1}(\cdot)$.
\end{enumerate}
\begin{proposition}
\label{prop_11_last}
If $t+\mathrm{tr}(A)<0$, then 
\begin{align}
\label{RK_bound_minus}
 \Pr(\boldsymbol{\Delta}>t)\le 1-e^{\frac{t+\mathrm{tr}(A)}{2s^-}}\sum_{l=0}^L c_l^{\textup{RK}^-}\sum_{m=0}^{\frac{\check{u}}{2}+l-1}\frac{1}{2^m}\sum_{q=0}^m\frac{c^-_{q}}{(m-q)!}\Big(-\frac{t+\mathrm{tr}(A)}{s^-}\Big)^{m-q}
\end{align}
and the corresponding gap is bounded above by 
\begin{align}
\label{RK_error_minus}
 \mathsf{Gap}^-=  1-\sum_{l=0}^Lc_l^{\textup{RK}^{-}}+\sum_{l=0}^{L}w^{-}_{l}c_l^{\textup{RK}^-}.
\end{align}
If $t+\mathrm{tr}(A)\ge0$, then 
\begin{align}
\label{RK_bound_plus}
 \Pr(\boldsymbol{\Delta}>t)\le 1-\sum_{l=0}^Lc_l^{\textup{RK}^+}+e^{-\frac{t+\mathrm{tr}(A)}{2s^{+}}}\sum_{l=0}^L c_l^{\textup{RK}^+}\sum_{m=0}^{\frac{\hat{v}}{2}+l-1}\frac{1}{2^m}\sum_{q=0}^m\frac{c^+_{q}}{(m-q)!}\Big(\frac{t+\mathrm{tr}(A)}{s^+}\Big)^{m-q}
\end{align}
and the corresponding gap is bounded above by 
\begin{align}
\label{RK_error_plus}
 \mathsf{Gap}^+=  1-\sum_{l=0}^Lc_l^{\textup{RK}^+}+\sum_{l=0}^{L}w^+_{l}c_l^{\textup{RK}^+}.
\end{align}
Here, the weights $w^-_l$ and $w^+_l$ are defined by 
\begin{align}
\label{}
 w^-_l=\left\{\begin{array}{cc}
    0  &  \textrm{$u$ is even}  \\
    \delta_{u+2l-1}  &   \textrm{$u$ is odd}
\end{array}\right.,\quad w^+_l=\left\{\begin{array}{cc}
    0  &  \textrm{$v$ is even}  \\
    \delta_{v+2l}  &   \textrm{$v$ is odd}
\end{array}\right..
\end{align}
\end{proposition}
\begin{proof}
See Section~\ref{secXII}.
\end{proof}
If $t+\mathrm{tr}(A)< 0$,  two sets of coefficients, $(c_l^{\text{RK}^-})_{l=0}^{L}$ and $(c_q^{-})_{q=0}^{\frac{\check{u}}{2}+L-1}$, must be computed recursively. For each $1\leq l\leq L$, calculating $c_l^{\text{RK}^-}$ requires $l-1$ additions and $l+1$ multiplications. Similarly, for each $1\le q\leq \frac{\check{u}}{2}+L-1$, calculating $c^-_q$ requires $q-1$ additions and $q+1$ multiplications. This yields a total workload of additions given by  
\begin{align}
\label{add_wl}
\sum_{l=1}^L (l-1)+\sum_{q=1}^{\frac{\check{u}}{2}+L-1}(q-1)=\frac{(L-1)L}{2}+\frac{1}{2}\Big(\frac{\check{u}}{2}+L-2\Big)\Big(\frac{\check{u}}{2}+L-1\Big),
\end{align}
and multiplications given by
\begin{align}
\label{multiplication_wl}
   \sum_{l=1}^L (l+1)+\sum_{q=1}^{\frac{\check{u}}{2}+L-1}(q+1)=\frac{L(L+3)}{2}+\frac{1}{2}\Big(\frac{\check{u}}{2}+L-1\Big)\Big(\frac{\check{u}}{2}+L+2\Big).
\end{align}
Analogously, if $t+\mathrm{tr}(A)\ge0$, the number of additions and multiplications needed to calculate the coefficients $(c_l^{\text{RK}^+})_{l=0}^{L}$ and $(c_q^{+})_{q=0}^{\frac{\hat{v}}{2}+L-1}$ are given by (\ref{add_wl}) and (\ref{multiplication_wl}), respectively, with $\check{u}$ replaced by $\hat{v}$.  }

 \textcolor{Black}{In Subsection~\ref{2F}, we elaborate on the approximation gaps in (\ref{RK_error_minus}) and (\ref{RK_error_plus}). Specifically, we provide a closed-form formula for the numbers $\delta_k$ in (\ref{KSD}) and summarize known facts about the Ruben-Kotz coefficients. In particular, we explain why $\sum_{l=0}^\infty c_l^{\text{RK}^{-}}=\sum_{l=0}^\infty c_l^{\text{RK}^{+}}=1$. An immediate consequence of this is that if either $t+\mathrm{tr}(A)\le0$ with even $u$ or $t+\mathrm{tr}(A)\ge0$ with even $v$, the approximation gap vanishes as $L\to\infty$. }
 
   \textcolor{Black}{The results in Proposition~\ref{prop_11_last} lead to single-series formulas for the CDF $F_{\underline{\boldsymbol{x}}^{\mathsf{T}}A\underline{\boldsymbol{x}}}(z)$ when either $z<0$ and $u$ is even, or $z\ge0$ and $v$ is even. Provost and Rudiuk~\cite{Provost} also build upon the Ruben-Kotz expansion in~(\ref{RK_expansion}) to develop their general formulas for $F_{\underline{\boldsymbol{x}}^{\mathsf{T}}A\underline{\boldsymbol{x}}}(\cdot)$. In their approach, they first derive double-series formulas for the PDF $p_{\underline{\boldsymbol{x}}^{\mathsf{T}}A\underline{\boldsymbol{x}}}(\cdot)$, which involve Whittaker functions, functions written in terms of the confluent hypergeometric functions. By integrating these PDF formulas, they arrive at their triple-series formulas for the CDF. In contrast, our approach utilizes the Ruben-Kotz expansion in (\ref{RK_expansion}) alongside a conditioning argument to compute $F_{\underline{\boldsymbol{x}}^{\mathsf{T}}A\underline{\boldsymbol{x}}}(\cdot)$ directly, bypassing the need to find the PDF first.    }

\item Applications in statistical signal processing and wireless communications are considered where the modified bounds prove highly effective. A common feature in both applications is that the tail parameter~$t$ depends on the underlying matrix~$A$, i.e., $t$ and $A$ do not vary independently. During the first half of Section~\ref{sec3}, we consider a binary statistical test of \textit{simple} hypotheses on the covariance structure of a Gaussian vector. Specifically, we adopt the Neyman-Pearson paradigm with a significance level of $\alpha$ and an underlying threshold $\gamma$ on the log-likelihood-ratio. Type~I and Type~II error probabilities, denoted by $p_{1|0}$ and $p_{0|1}$, respectively, are identified as tail probabilities of Gaussian quadratic forms.  We denote the largest achievable power of the test by $p_{1|1,\alpha}=\max_{\gamma: p_{1|0}\le\alpha}p_{1|1}$ and the optimal threshold~$\gamma$ that achieves $p_{1|1,\alpha}$~by~$\gamma_{\alpha}$. We derive seven \textit{upper} bounds $\gamma_{\alpha,\mathsf{B}}$  on $\gamma_{\alpha}$, where~$\mathsf{B}$ represents HW, $m_\infty$, strong $\chi^2$, their modified versions and RK. Plugging $\min_{\mathsf{B}}\gamma_{\alpha,\mathsf{B}}$ for $\gamma$ in another set of seven \textit{lower} bounds on $p_{1|1}=1-p_{0|1}$ leads us to seven \textit{lower} bounds on $p_{1|1,\alpha}$. Numerical results reveal that all three modified bounds present considerable improvements over their original versions, with the modified $m_\infty$ and $\chi^2$ bounds being the sharpest among the modified bounds in most cases. \textcolor{Black}{As expected, RKIs can present the tightest bounds, especially in cases where the approximation gap vanishes by increasing the truncation level. However, this comes at a heavy computational cost. For example, our numerical results in Subsection~\ref{3_A} show that for $n=30$, the RK bound exceeds the best modified $m_\infty$ bound for $L=500$, requiring at least hundreds of thousands of arithmetic operations.\footnote{See Fig.~\ref{fig_p11_n}.}  Moreover, the rate of convergence for the series $\sum_{l=0}^\infty c_l^{\textup{RK}^{\pm}}$ is exceptionally low for certain values of $n$. For instance, when $n=31$, even $L=1000$ remains insufficient for the RK bound to surpass the best modified $m_\infty$ bound. } While the task of approximating \(p_{1|{}1,\alpha }\) is traditionally left to numerical integrations via Monte Carlo simulations, to the best of our knowledge, this is the first attempt at deriving such analytical bounds. 

In the second half of Section~\ref{sec3}, we revisit the problem of performance evaluation for unitary space-time modulation over MISO Rayleigh flat-fading channels where neither the transmitter nor the receiver are aware of the channel coefficients.\footnote{We assume there is only one receiver antenna for simplicity. The results can be easily extended to an arbitrary number of receiver antennas. } This problem was first studied in~\cite{Hochwald1}, where pairwise error probabilities under maximum-likelihood decoding were expressed as the tail probabilities of Gaussian quadratic forms. Following a characteristic function approach similar to~\cite{Kay_Detection},  Hochwald and Marzetta~\cite{Hochwald1} derived a simple upper bound on these error probabilities, hereafter referred to as the~HM~bound. Denoting the number of transmitter antennas and the average signal-to-noise ratio at the receiver by $M$ and $\mathsf{SNR}$, respectively, the~HM~bound scales like~$\mathsf{SNR}^{-M}$ as $\mathsf{SNR}$ grows large.   We derive alternative upper bounds using the results introduced earlier in this section. We show that the HW~bound, the $m_\infty$~bound, their viable modifications and the original strong $\chi^2$~bound all saturate as $\mathsf{SNR}$ increases. Consequently, these bounds become ineffective in the high signal-to-noise ratio regime.  Interestingly, the modified strong $\chi^2$~bound of \textit{largest} possible order is always viable, it achieves the same signal-to-noise ratio scaling, $\mathsf{SNR}^{-M}$,  as the~HM~bound and it offers noticeable improvements over the~HM~bound for a wide range of values of $\mathsf{SNR}$, as demonstrated by an example.

  \end{enumerate}
  \subsection{Notation}
  Random variables are denoted by boldfaced letters such as $\boldsymbol{x}$ with expectation $\mathbb{E}[\boldsymbol{x}]$ and variance $\mathrm{Var}(\boldsymbol{x})$. Vectors are denoted by an underline such as $\underline{\boldsymbol{x}}$ with realization $\underline{x}$. \textcolor{Black}{For a random vector $\underline{\boldsymbol{x}}$ and a measurable function $f$, we define the notation 
  \begin{align}
\label{}
  \mathbb{E}_{\underline{x}\sim \underline{\boldsymbol{x}}}[f(\underline{x})]=\mathbb{E}[f(\underline{\boldsymbol{x}})].
\end{align}
This notation will be useful when we perform conditioning to evaluate certain probabilities.} For a continuous random vector~$\underline{\boldsymbol{x}}$, its probability density function (PDF) is denoted by $p_{\underline{\boldsymbol{x}}}(\cdot)$. The cumulative distribution function~(CDF) of a random variable $\boldsymbol{x}$ is denoted by $F_{\boldsymbol{x}}(\cdot)$. \textcolor{Black}{If $\boldsymbol{x}\sim \chi^2_n$ is a chi-squared random variable with $n$ degrees of freedom, we denote its CDF by $F_n(\cdot)$, i.e., 
  \begin{align}
    F_n(x)=F_{\chi^2_n}(x)=\Pr(\chi^2_n\leq x),\,\,\,x\in \mathbb{R}.
\end{align}}A zero-mean Gaussian vector $\underline{\boldsymbol{x}}$ with covariance matrix $C$ is indicated by $\underline{\boldsymbol{x}}\sim N(\underline{0},C)$. The $n\times n$ identity matrix is denoted by~$I_n$. The transpose of a matrix $A$ is denoted by $A^{\mathsf{T}}$ and its rank is denoted by
  \begin{align}
r(A)=\mathrm{rank}(A).
\end{align}
 The (real) eigenvalues of an $n\times n$ symmetric matrix $A$ are denoted by $\lambda_1,\dots, \lambda_n$ and the smallest and  largest eigenvalues of $A$ are denoted by 
 \begin{align}
\lambda_{\min}(A)=\min_{1\leq i\leq n}\lambda_i, \quad \lambda_{\max}(A)=\max_{1\leq i\leq n}\lambda_i.
\end{align}
  The operator norm of $A$ is 
  \begin{align}
\|A\|=\max_{1\leq i\leq n}|\lambda_i|
\end{align} 
and Schatten $p$-norm of $A$ for $p\geq 1$ is defined by  
\begin{align}
\|A\|_p=\big(\sum_{i=1}^n |\lambda_i|^p\big)^{1/p}.
\end{align}
The case $p=2$ reduces to Hilbert-Schmidt norm of $A$. If $A$ is positive-semidefinite, then 
\begin{align}
\|A\|_p=(\mathrm{tr}(A^p))^{1/p},
\end{align} 
where $\mathrm{tr}(\cdot)$ is the trace operator. 
 
\section{Analysis of the results in Subsection~\ref{1B}}
\label{sec2}
 \subsection{Improving LMI} 
 \label{2A}
\begin{figure}
   \centering
    \includegraphics[width=0.45\textwidth]{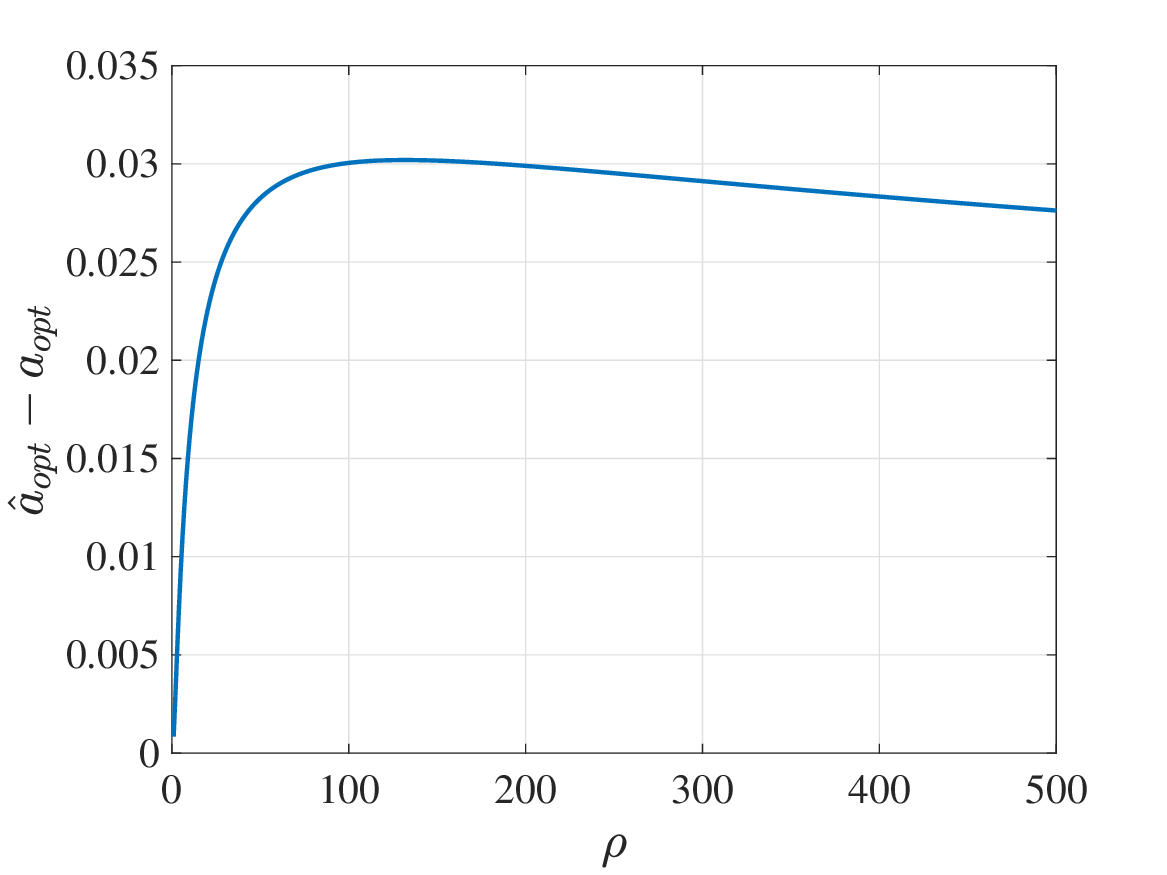}
     \caption{The plot for $\hat{a}_{opt}-a_{opt}$ in terms of $\rho$. Both $a_{opt}$ and $\hat{a}_{opt}$ depend on the matrix $A$ and the tail parameter $t$ through the ratio $\rho$. We observe that $\hat{a}_{opt}-a_{opt}<0.035$ regardless of the value of $\rho$.  }
    \label{pict_777}
\end{figure}

\begin{figure*}[t]
\centering
\subfigure[]{
\includegraphics[scale=0.4]{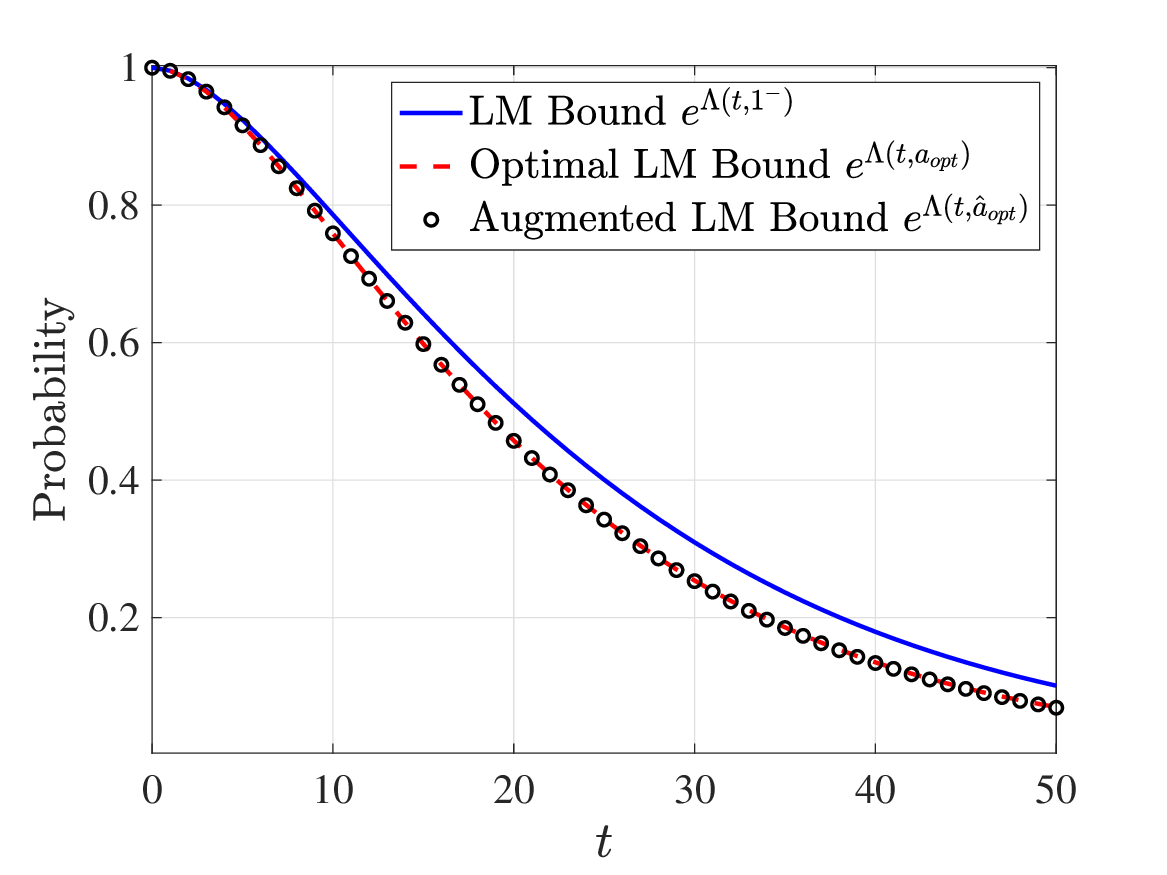}
\label{fig11}
}
\subfigure[]{
\includegraphics[scale=0.4]{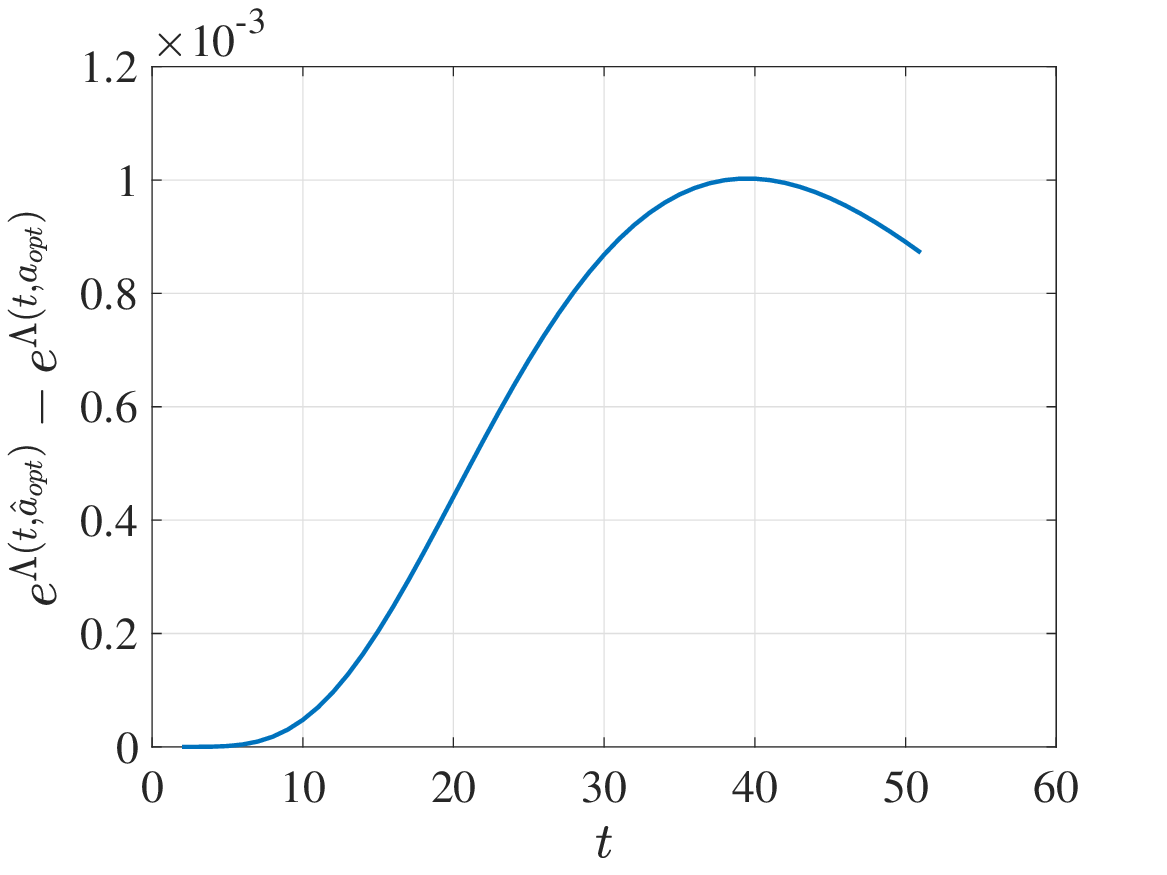}
\label{fig12}
}
\label{fig:subfigureExample}
\caption[Optional caption for list of figures]{ Panel (a) presents plots of the upper bounds in  (\ref{augmented_LM}),  (\ref{optimal_LM}) and (\ref{LM_itself}) in terms of $0\leq t\leq 50$ for a $3\times 3$ positive-semidefinite matrix~$A$ with eigenvalues $0.0206$, $3.6428$ and $6.1735$. The difference between $e^{\Lambda(t,a_{opt})}$ and $e^{\Lambda(t,\hat{a}_{opt})}$ is not noticeable in this plot. Panel~(b) sketches the difference $e^{\Lambda(t, \hat{a}_{opt})}-e^{\Lambda(t, a_{opt})}$.}
\label{fig8}
\end{figure*}

 Consider the improved LMI given in (\ref{en_LM}) in Proposition~\ref{prop_2}. It is easily checked that $\Lambda(t,a)< 0$ for all $t>0$ and $\frac{2}{3}<a<1$ and hence, the bound in (\ref{en_LM}) is nontrivial.  Let us denote the ratio $\frac{\alpha t}{\beta}$ by~$\rho$,~i.e., 
  \begin{align}
  \label{rho_first}
\rho=\frac{\alpha t}{\beta}.
\end{align}
 Recall $c=\frac{3a-2}{2(1-a)^2}$ in (\ref{toy}). The inequalities $\rho>c$ and $\rho\leq c$ given in the formulation of $\Lambda(t,a)$ in (\ref{en_exp}) are solved for the parameter $a$ as $\frac{2}{3}<a<\hat{a}_{opt}$ and $\hat{a}_{opt}\leq a<1$, respectively, where $\hat{a}_{opt}$ is defined by 
 \begin{align}
 \label{a_star}
\hat{a}_{opt}=\frac{4\rho+3-\sqrt{8\rho+9}}{4\rho}.
\end{align}
The reason for choosing this notation becomes clear in a moment. The value of $\hat{a}_{opt}$ always lies in the admissible interval $(\frac{2}{3},1)$ and one can check that 
\begin{align}
\lim_{\rho\to0^+}\hat{a}_{opt}=\frac{2}{3},\quad \,\lim_{\rho\to\infty}\hat{a}_{opt}=1.
\end{align}
Replacing $b$ and $c$ by their values in terms of $a$ and  after some algebra, we find 
\begin{align}
\Lambda(t,a)=\left\{\begin{array}{cc}
    \frac{t}{\alpha}\frac{(3a-2)(2\rho a^2+(3-2\rho)a-2)}{4\rho a^2(2a-1)(1-a)}  & \frac{2}{3}<a<\hat{a}_{opt}  \\
    \frac{t}{\alpha}\big(-\frac{1}{2a}+\frac{1}{2a^2\rho}(\sqrt{1+2a\rho}-1)\big)  & \hat{a}_{opt}\leq a<1
\end{array}\right..
\end{align}
It is easily seen that $\Lambda(t,a)$ is continuous over the whole interval $\frac{2}{3}<a<1$ and it is strictly increasing over the interval $\hat{a}_{opt}\leq a<1$. As such, the absolute minimum value for $\Lambda(t,a)$ with respect to $a$ occurs somewhere over the interval $(\frac{2}{3},\hat{a}_{opt}]$. Looking for the critical numbers of the function $a\mapsto \Lambda(t,a)$ over the interval $(\frac{2}{3},\hat{a}_{opt})$ leads us to the quintic polynomial equation
\begin{align}
\label{quintic}
12\rho a^5+(36-40\rho)a^4+(48\rho-99)a^3+(104-24\rho)a^2+(4\rho-48)a+8=0.
\end{align} 
It is not hard to show that this equation has a unique root denoted by $a_{opt}$ inside the aforementioned interval which must be solved for numerically. Note that both $a_{opt}$ and $\hat{a}_{opt}$ depend on the matrix $A$ and the tail parameter~$t$ only through the ratio $\rho=\frac{t\|A\|}{\|A\|_2^2}$. Fig.~\ref{pict_777} shows the graph for $\hat{a}_{opt}-a_{opt}$ in terms of $0\leq \rho\leq 500$. This numerical observation suggests that $\hat{a}_{opt}-a_{opt}$ never exceeds $0.035$. Having this evidence, we will use $\hat{a}_{opt}$ in place of $a_{opt}$. We will name 
\begin{align}
\label{augmented_LM}
\Pr(\boldsymbol{\Delta}>t)\le e^{\Lambda(t,\hat{a}_{opt})},
\end{align}
the \textit{augmented} LMI and 
\begin{align}
\label{optimal_LM}
\Pr(\boldsymbol{\Delta}>t)\le e^{\Lambda(t,a_{opt})},
\end{align}
the \textit{optimal} LMI. 
As mentioned earlier in Section~\ref{sec1}, the original LMI is recovered from~(\ref{en_LM}) by letting $a\to1^-$ and hence, it is given by  
\begin{align}
\label{LM_itself}
\Pr(\boldsymbol{\Delta}>t)\le e^{\Lambda(t,1^-)}.
\end{align}   
As before, since $\Lambda(t,a)$ is strictly increasing with respect to $a$ over the interval $[\hat{a}_{opt},1)$, we conclude that the augmented LM~bound in (\ref{augmented_LM}) is tighter than the classical LM bound in (\ref{LM_itself}). We demonstrate the improvement through an example. Consider a $3\times 3$ positive-semidefinite matrix $A$ with eigenvalues $0.0206$, $3.6428$ and $6.1735$.\footnote{These numbers were generated randomly.} Panel~(a) in Fig.~\ref{fig8} shows the bounds $e^{\Lambda(t,1^-)}, e^{\Lambda(t, a_{opt})}$ and $e^{\Lambda(t, \hat{a}_{opt})}$ with respect to $0\leq t\leq 50$. The difference between $e^{\Lambda(t,a_{opt})}$ and $e^{\Lambda(t,\hat{a}_{opt})}$ is not noticeable in this plot. Panel~(b) sketches the difference $e^{\Lambda(t, \hat{a}_{opt})}-e^{\Lambda(t, a_{opt})}$.

  \subsection{Beyond HWI}
  \label{2B}
\begin{figure*}[t]
\centering
\subfigure[]{
\includegraphics[scale=0.4]{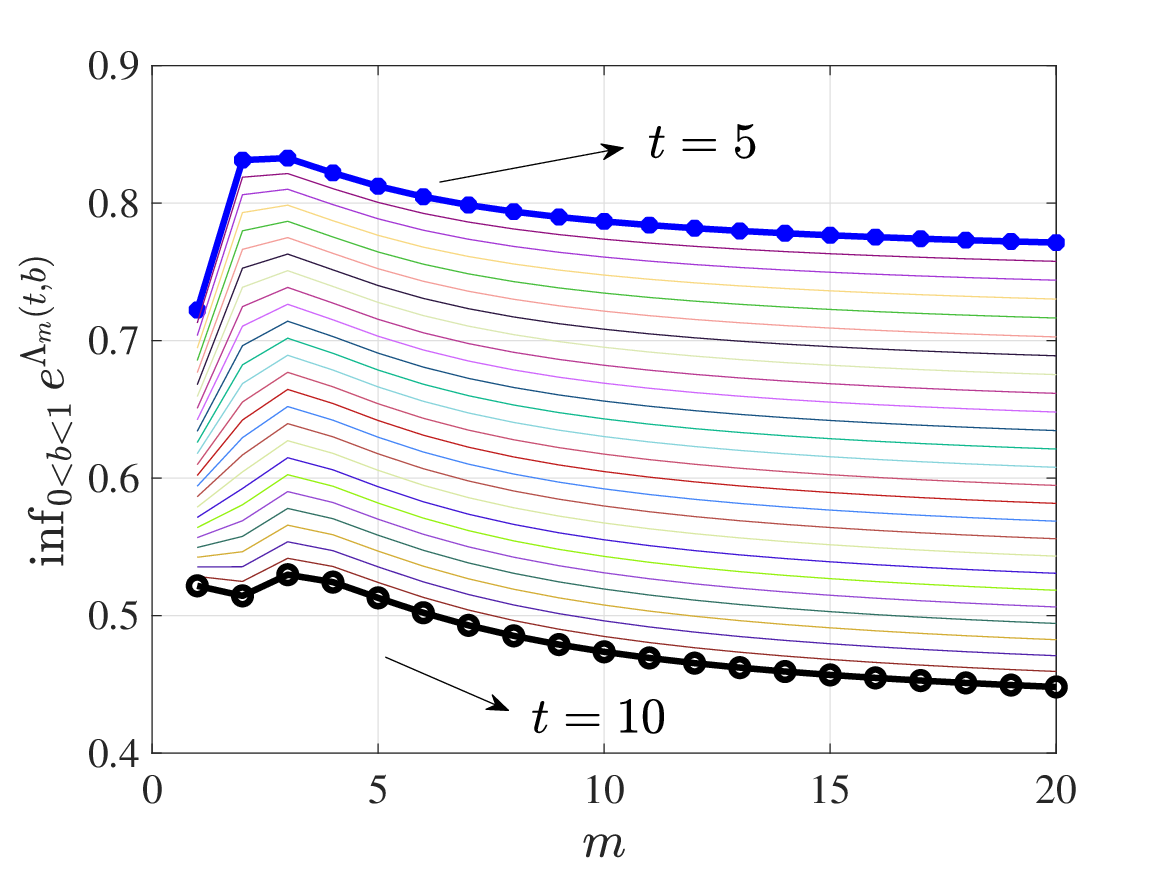}
\label{fig10_1}
}
\subfigure[]{
\includegraphics[scale=0.4]{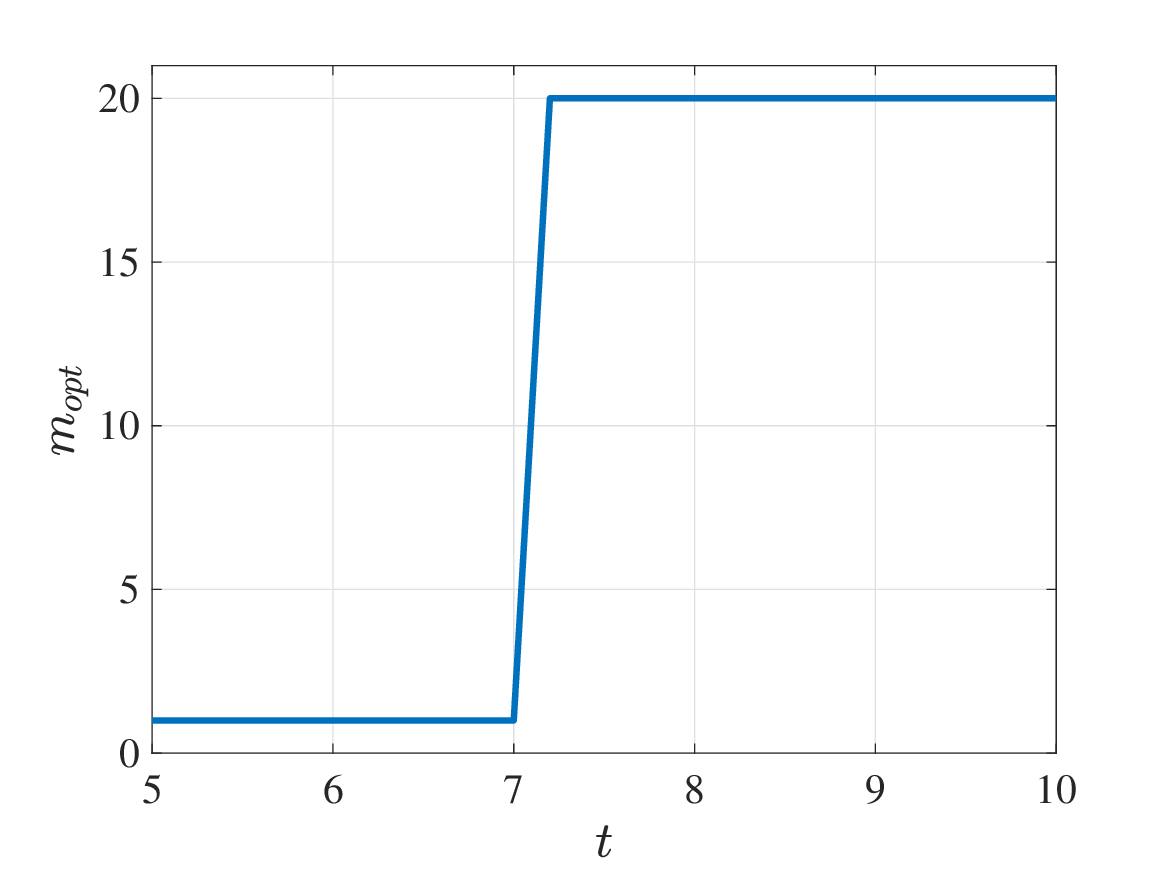}
\label{fig10_2}
}
\label{fig:subfigureExample}
\caption[Optional caption for list of figures]{ Consider a $3\times 3$ symmetric matrix $A$ with eigenvalues $-2.2385$, $-0.6164$ and $0.0880$. Panel (a) presents plots of $\inf _{0<b<1}e^{\Lambda_m(t,b)}$ with respect to $1\leq m\leq 20$ for $5\leq t\leq 10$ in steps of $0.2$. Panel (b) shows the optimum value for $m$ denoted by $m_{opt}$ over the domain $1\leq m\leq 20$ at which $\inf _{0<b<1}e^{\Lambda_m(t,b)}$ is the smallest. If Conjecture~1 is true, then these plots imply that~$\tau_c\le 7.2$.}
\label{fig10}
\end{figure*}

   Next, we investigate the sequence of bounds presented in Proposition~\ref{prop_4}. 
   Two initial remarks are in order: 
  \begin{enumerate}
  \item As we mentioned earlier, the upper bound $e^{\Lambda_m(t,b)}$ must be minimized over $0<b<1$. Such minimization results in meaningful bounds in the sense that $\inf_{0<b<1}e^{\Lambda_m(t,b)}\leq 1$. This is due to the simple observation that $\lim_{b\to0^+}\Lambda_m(t,b)=0$. 
  \item The coefficient $\kappa_m(b)$ in the exponent achieves a global maximum value somewhere inside the interval $0<b<1$. Since $\lim_{b\to0^+}\theta_m(b)=\frac{1}{m+1}$ and $\lim_{b\to1^-}\theta_m(b)=\infty$, we have $\lim_{b\to0^+}\kappa_m(b)=\lim_{b\to1^-}\kappa_m(b)=0$. 
\end{enumerate}
   Let us study the bounds $e^{\Lambda_m(t,b)}$ for $m\ge1$ both numerically and analytically. For our numerical study, consider a $3\times 3$ symmetric matrix $A$ with eigenvalues $-2.2385$, $-0.6164$ and $0.0880$.\footnote{These numbers were generated randomly.} 
   As stated in Conjecture~\ref{conj_1} in Subsection~\ref{1B}, there exists a critical threshold~$\tau_c\ge0$ such that if $t<\tau_c$, then $\inf _{0<b<1}e^{\Lambda_m(t,b)}$  is the smallest when $m=1$ and if $t>\tau_c$, then $\inf _{0<b<1}e^{\Lambda_m(t,b)}$ is the smallest when $m$ grows to infinity. To better observe this phase transition, Fig.~\ref{fig10_1} offers plots of $\inf _{0<b<1}e^{\Lambda_m(t,b)}$  in terms of $1\leq m\leq 20$ for the values $5\leq t\leq 10$ in steps of $0.2$, i.e., $t=5, 5.2, 5.4,\dots,9.8,10$. Fig.~\ref{fig10_2} shows the optimum value for $m$ denoted by $m_{opt}$ over the domain $1\leq m\leq 20$ at which $\inf _{0<b<1}e^{\Lambda_m(t,b)}$ is the smallest. We have $m_{opt}=1$ for $t\leq 7$ and $m_{opt}=20$ for $t\ge7.2$. This tells us that if Conjecture~1 is true, then the critical threshold must satisfy $\tau_c\leq 7.2$. In Corollary~\ref{coro_3} in below, we will present a theoretical upper bound on~$\tau_c$.  

We continue by stating two corollaries to Proposition~\ref{prop_4}.  For our first corollary, we study the behaviour of the bounds in Proposition~\ref{prop_4} as $m$ grows to infinity. This yields what we refer to as the $m_\infty$~inequality.
\begin{corollary}
\label{coro_2}
For a symmetric matrix $A$, 
\begin{align}
\label{inftym}
\Pr(\boldsymbol{\Delta}>t)\leq \liminf_{m\to\infty}\inf_{0<b<1}e^{\Lambda_m(t,b)}\leq \Big(1+\frac{t}{r(A)\|A\|}\Big)^{\frac{r(A)}{2}}e^{-\frac{t}{2\|A\|}}.
\end{align} 
\end{corollary}
\begin{proof}
The first inequality in (\ref{inftym}) is trivial. For the second inequality, we note that\footnote{If $f_n:\mathcal{A}\to\mathbb{R}$ is a sequence of functions on a set $\mathcal{A}$ of real numbers, then $\inf_{a\in \mathcal{A}}f_n(a)\leq f_n(x)$ for every $x\in \mathcal{A}$ and $n\ge1$. Taking limit inferior from both sides, we get $\liminf_{n\to\infty}\inf_{a\in \mathcal{A}}f_n(a)\leq \liminf_{n\to\infty}f_n(x)$ for every $x\in \mathcal{A}$. This says $\liminf_{n\to\infty}\inf_{a\in \mathcal{A}}f_n(a)$ is a lower bound on the set of numbers $\liminf_{n\to\infty}f_n(x)$ for $x\in \mathcal{A}$. Thus, $\liminf_{n\to\infty}\inf_{a\in \mathcal{A}}f_n(a)\leq \inf_{x\in \mathcal{A}}\liminf_{n\to\infty}f_n(x)$ as the infimum of a set is the largest lower bound for that set.  } 
\begin{align}
\label{chick1}
\liminf_{m\to\infty}\inf_{0<b<1}e^{\Lambda_m(t,b)}\leq \inf_{0<b<1}\liminf_{m\to\infty}e^{\Lambda_m(t,b)},
\end{align}
for all $t\ge0$. The limit $\lim_{m\to\infty}\Lambda_m(t,b)$ exists. To see this, first we derive a different formula for $\theta_m(b)$.  Recall the Maclaurin series $-\ln(1-b)=\sum_{i=1}^\infty\frac{b^i}{i}$. Then 
\begin{align}
\theta_m(b)=-\frac{1}{b^{m+1}}\Big(-\sum_{i=m+1}^\infty\frac{b^{i}}{i}\Big)=\sum_{i=m+1}^\infty\frac{b^{i-m-1}}{i}.
\end{align}
Changing the index $i$ to $i+m+1$ results in $\theta_m(b)=\sum_{i=0}^\infty\frac{b^i}{i+m+1}$. We have
\begin{align}
\label{squeeze_bah}
\frac{1}{m+1}=\frac{b^0}{0+m+1}\leq \theta_m(b)\leq \sum_{i=0}^\infty\frac{b^i}{m+1}=\frac{1}{m+1}\sum_{i=0}^\infty b^i=\frac{1}{(1-b)(m+1)}.
\end{align}
Since $\lim_{m\to\infty}(m+1)^{\frac{1}{m}}=1$, it follows by squeezing that 
\begin{align}
\lim_{m\to\infty}(\theta_m(b))^{\frac{1}{m}}=1.
\end{align}
This verifies  
\begin{align}
\label{cherry1}
\lim_{m\to\infty}\kappa_m(b)=\frac{1}{2}\min\{1,b\}=\frac{b}{2}.
\end{align}
It is also evident that 
\begin{align}
\label{cherry2}
\lim_{m\to\infty}\|A\|_{m+1}^{1+\frac{1}{m}}=\|A\|
\end{align}
and
\begin{align}
\label{cherry3}
\lim_{m\to\infty}\sum_{k=2}^m\frac{b^k}{k}=\sum_{k=2}^\infty\frac{b^k}{k}=-\ln(1-b)-b,
\end{align}
where we have used the Maclaurin series for $-\ln(1-b)$ once more. By (\ref{cherry1}), (\ref{cherry2}) and (\ref{cherry3}), 
\begin{align}
\label{coffee}
\lim_{m\to\infty}\Lambda_m(t,b)=-\frac{r(A)}{2}(\ln(1-b)+b)-\frac{bt}{2\|A\|}.
\end{align}
The right-hand side of (\ref{coffee}) is minimized at 
\begin{align}
b=b_{opt}=\frac{t}{t+r(A)\|A\|}.
\end{align}
Plugging this value for $b$ on the right-hand side of (\ref{coffee}), we arrive at 
\begin{align}
\label{chick2}
\inf_{0<b<1}\lim_{m\to\infty}e^{\Lambda_m(t,b)}=\Big(1+\frac{t}{r(A)\|A\|}\Big)^{\frac{r(A)}{2}}e^{-\frac{t}{2\|A\|}}.
\end{align}
By (\ref{chick1}) and (\ref{chick2}), the proof is complete. 
\end{proof}
Using the inequality $(1+x)^{\frac{1}{x}}<e$ for every $x>0$, the $m_\infty$~bound in Corollary~\ref{coro_2} is nontrivial in the sense that it is less than one. 

Our second corollary provides an upper bound on the critical threshold $\tau_c$ in Conjecture~1. 
\begin{corollary}
\label{coro_3}
Let $\hat{\tau}_c$ be the largest positive solution (if it exists) in $t$ to the equation 
\begin{align}
\label{koochil}
r(A)\ln\Big(1+\frac{t}{r(A)\|A\|}\Big)-\frac{t}{\|A\|}+2\kappa\min\Big\{\frac{t^2}{\|A\|_2^2},\frac{t}{\|A\|}\Big\}=0,
\end{align} 

where $\kappa$ is the constant in HWI given in Proposition~\ref{prop_1}. If Conjecture~1 is true, then $\tau_c\leq \hat{\tau}_c$. 
\end{corollary}
\begin{proof}
Define the functions 
\begin{align}
f(t)=\exp\Big(-\kappa\min\Big\{\frac{t^2}{\|A\|_2^2},\frac{t}{\|A\|}\Big\}\Big)
\end{align}
and 
\begin{align}
g(t)=\liminf_{m\to\infty}\inf_{0<b<1}e^{\Lambda_m(t,b)}.
\end{align}
By definition, the critical threshold $\tau_c$ is given by 
\begin{align}
\tau_c=\sup\{t>0: g(t)\geq f(t)\}.
\end{align}
By Corollary~\ref{coro_2},  $g(t)\leq(1+\frac{t}{r(A)\|A\|})^{\frac{r(A)}{2}}e^{-\frac{t}{2\|A\|}}$. It follows that
\begin{align}
\label{quiet_kk}
\tau_c\leq\sup\Big\{t>0: \Big(1+\frac{t}{r(A)\|A\|}\Big)^{\frac{r(A)}{2}}e^{-\frac{t}{2\|A\|}}\ge f(t)\Big\}.
\end{align}
The right hand side of (\ref{quiet_kk}) is the largest solution in $t$ to the equation $\big(1+\frac{t}{r(A)\|A\|}\big)^{\frac{r(A)}{2}}e^{-\frac{t}{2\|A\|}}= f(t)$, which simplifies to~(\ref{koochil}). 
\end{proof}
To illustrate Corollary~\ref{coro_3}, consider the symmetric matrix $A$ in Fig.~\ref{fig10}. The corresponding plots showed that $\tau_c\leq 7.2$. Using Corollary~\ref{coro_3}, we get the sharper bound $\tau_c\leq \hat{\tau}_c\approx 6.147$. 

It is easy to show that if all eigenvalues of $A$ have the same absolute value, then the equation in~(\ref{koochil}) has no positive root. This is the motivation behind the statement included in Conjecture~\ref{conj_1}, which claims that if all eigenvalues of $A$ have the same absolute value, then $\tau_c=0$, implying that $m=\infty$ outperforms $m=1$ for all $t>0$.

\subsection{Upper bounds on $\Pr(\boldsymbol{\Delta}>t)$ that depend on the positive-definite matrix $A$ only through $\|A\|$ }
\label{2D}


 Let $A$ be full-rank, i.e., $r(A)=n$. Then 
\begin{align}
\label{ineq_norm}
\|A\|_2\leq \sqrt{r(A)}\,\|A\|=\sqrt{n}\|A\|.
\end{align}
The exponents in HWI and LMI increase with $\|A\|_2$. This is clear for HWI, which directly yields its relaxed version in (\ref{relaxed_HW}).  As for LMI in (\ref{LM2}), note that $\sqrt{\|A\|_2^2+2\|A\|t}-\|A\|_2$ is positive and decreasing with respect to $\|A\|_2$. Consequently, $-(\sqrt{\|A\|_2^2+2\|A\|t}-\|A\|_2)^2$ increases with $\|A\|_2$. Combining this with (\ref{ineq_norm}) gives  
\begin{align}
\label{relaxed_LM2}
\Pr(\boldsymbol{\Delta}>t)\leq \exp\bigg(-\Big(\frac{\sqrt{n\|A\|^2+2\|A\|t}-\sqrt{n}\|A\|\,}{2\|A\|}\Big)^2\bigg).
\end{align}
Factoring $\sqrt{n}\|A\|$ in $\sqrt{n\|A\|^2+2\|A\|t}-\sqrt{n}\|A\|$ results in the relaxed LMI in (\ref{relaxed_LM}). As stated in Proposition~\ref{prop_6}, the $m_\infty$~inequality is always tighter than both the relaxed HWI and the relaxed LMI. One can also consider the possibility of relaxing the augmented and optimal LM inequalities in (\ref{augmented_LM}) and~(\ref{optimal_LM}),~respectively. This requires that we replace the ratio $\rho=\frac{t\|A\|}{\|A\|_2^2}$ in (\ref{rho_first}) by $\frac{t\|A\|}{(\sqrt{n}\|A\|)^2}=\frac{t}{n\|A\|}$. Numerical results indicate that the $m_\infty$~bound is still tighter that these relaxed inequalities and hence, we do not pursue this further. 

To derive the weak $\chi^2$~bound in (\ref{cvb}), we first present an alternative upper bound on $\Pr(\boldsymbol{\Delta}>t)$ for an arbitrary symmetric matrix $A$ without relying on Markov's inequality. By (\ref{innooo}) in Section~\ref{sec4}, 
\begin{align}
\boldsymbol{x}^{\mathsf{T}}A\boldsymbol{x}\leq \lambda^*_{\max}(A)\chi^2_{r(A)},
\end{align}
where $\lambda^*_{\max}(A)$ is the largest \textit{nonzero} eigenvalue of $A$ and $r(A)$ is the rank of $A$ and $\chi^2_{r(A)}$ is a chi-squared random variable with $r(A)$ degrees of freedom. We consider two cases.  If $\lambda^*_{\max}(A)>0$, then 
  \begin{align}
\Pr(\boldsymbol{\Delta}>t)&\leq  \Pr(\lambda^*_{\max}(A)\chi^2_{r(A)}-\mathrm{tr}(A)>t)=\Pr\Big(\chi^2_{r(A)}>\frac{t+\mathrm{tr}(A)}{\lambda^*_{\max}(A)}\Big)=1-F_{r(A)}\Big(\frac{t+\mathrm{tr}(A)}{\lambda^*_{\max}(A)}\Big).
\end{align}
If $\lambda^*_{\max}(A)<0$, then 
  \begin{align}
\Pr(\boldsymbol{\Delta}>t)\leq  \Pr\Big(\chi^2_{r(A)}<\frac{t+\mathrm{tr}(A)}{\lambda^*_{\max}(A)}\Big)=F_{r(A)}\Big(\frac{t+\mathrm{tr}(A)}{\lambda^*_{\max}(A)}\Big).
\end{align}
In one line, we have the strong $\chi^2$~bound in (\ref{chi2_bound}). If $A$ is positive-definite, then $r(A)=n$,  $\lambda^*_{\max}(A)=\|A\|$, $\mathrm{tr}(A)\geq \lambda^*_{\max}(A)$ and $\frac{t+\mathrm{tr}(A)}{\lambda^*_{\max}(A)}\geq \frac{t+\|A\|}{\|A\|}=1+\frac{t}{\|A\|}$. This leads to the weak $\chi^2$~bound in (\ref{cvb}).  
 It is well-known that
\begin{align}
\label{exact_chi2_CDF}
F_{n}(x)=\frac{\gamma(\frac{n}{2}, \frac{x}{2})}{\Gamma(\frac{n}{2})},\,\,\,x\ge0
\end{align}
where $\gamma(a,z)=\int_0^zv^{a-1}e^{-v}dv$ is the lower incomplete Gamma function and $\Gamma(a)=\int_0^\infty v^{a-1}e^{-v}dv$ 
is the Gamma function. For simplicity and to gain insight, let us assume that $n$ is even for the remainder of this subsection. Then $F_{n}(\cdot)$ admits a closed-form expression given by
\begin{align}
\label{even_n}
F_{n}(x)=1-e^{-\frac{x}{2}}\sum_{i=0}^{\frac{n}{2}-1}\frac{x^i}{2^ii!},\,\,\,x\ge0.
\end{align}
By (\ref{cvb}) and (\ref{even_n}), the weak $\chi^2$~bound becomes 
\begin{align}
\label{cvb_cvb}
\Pr(\boldsymbol{\Delta}>t)\leq  e^{-\frac{1}{2}(1+\frac{t}{\|A\|})}\sum_{i=0}^{\frac{n}{2}-1}\frac{1}{2^ii!}\Big(1+\frac{t}{\|A\|}\Big)^i=e^{-\frac{t}{2\|A\|}}\bigg(\frac{1}{\sqrt{e}} \sum_{i=0}^{\frac{n}{2}-1}\frac{1}{2^ii!}\Big(1+\frac{t}{\|A\|}\Big)^i\bigg).
\end{align}
 Both the $m_\infty$~bound in~(\ref{inftym}) and the weak $\chi^2$~bound in (\ref{cvb_cvb}) depend on $A$ and $t$ through the ratio 
 \begin{align}
\label{yek}
r=\frac{t}{\|A\|}. 
\end{align}
 We can rewrite the $m_\infty$~bound as
\begin{align}
\label{minfty_pol}
\Pr(\boldsymbol{\Delta}>t)\leq P^{(1)}_n(r)e^{-\frac{r}{2}}
\end{align}
and the weak $\chi^2$~bound as 
\begin{align}
\Pr(\boldsymbol{\Delta}>t)\leq P_n^{(2)}(r)e^{-\frac{r}{2}},
\end{align}
where the polynomials $P_n^{(1)}(r)$ and $P^{(2)}_n(r)$ are defined by 
\begin{align}
\label{poly_1}
P_n^{(1)}(r)=\Big(1+\frac{r}{n}\Big)^{\frac{n}{2}}
\end{align}
and 
\begin{align}
\label{poly_2}
P_n^{(2)}(r)= \frac{1}{\sqrt{e}}\sum_{i=0}^{\frac{n}{2}-1}\frac{(1+r)^i}{2^ii!},
\end{align}
respectively. Any comparison between the two bounds boils down to comparing the two polynomials in~(\ref{poly_1}) and (\ref{poly_2}). The first claim in Proposition~\ref{prop_7} is that $P^{(2)}_n(r)<P_n^{(1)}(r)$ for every~$r\geq 0$ when $n=2,4,6$. The second claim is that if $n\ge8$ is even, then $P_n^{(1)}(r)<P_n^{(2)}(r)$ if and only if $r$ lies in the open interval $(r_n, r'_n)$. The numbers $r_n, r'_n$ are, in fact, the only two positive roots of the polynomial $P_n^{(1)}(r)-P_n^{(2)}(r)$, as verified in the proof of Proposition~\ref{prop_7} in Section~\ref{sec10}. 

We close this subsection by recalling the well-known large-deviations lower bound on $F_{n}(\cdot)$ given by
\begin{align}
\label{ld}
F_{n}(x)\ge1- \mathrm{exp}\Big(-\frac{x}{2}+\frac{n}{2}\Big(1+\ln\frac{x}{n}\Big)\Big),\,\,\, x\geq n.
\end{align}
This bound is a consequence of Markov's inequality and a proof can be found in~\cite{Ghosh}.  We emphasize that~(\ref{ld}) only holds for $x\ge n$. Then one can further loosen the weak $\chi^2$~bound in (\ref{cvb}) as 
\begin{align}
\label{ld_chi2}
\Pr(\boldsymbol{\Delta}>t)&\leq  1-F_{n}(1+r)\leq  \mathrm{exp}\Big(-\frac{1+r}{2}+\frac{n}{2}\Big(1+\ln\frac{1+r}{n}\Big)\Big)=\frac{1}{\sqrt{e}}\Big(\frac{e(1+r)}{n}\Big)^{\frac{n}{2}}e^{-\frac{r}{2}},\,\,\, 1+r\geq n.
\end{align}
This is the large-deviations bound introduced earlier in (\ref{LDB}). It is always looser than the $m_\infty$~bound for every $n\ge2$, as shown at the end of the proof of Proposition~\ref{prop_6} in Section~\ref{sec9}.

\subsection{The $m_\infty$~bound can be tighter than the strong $\chi^2$~bound in the symmetric case and the LM bound in the positive-semidefinite case} 
\label{2_E}
We begin by presenting sufficient conditions under which the $m_\infty$~bound outperforms the \textit{strong} $\chi^2$~bound. Let $n$ be an even integer, and let $A$ be a full-rank symmetric matrix (i.e., $r(A)=n$) that satisfies  
\begin{align}
\label{ex111}
\mathrm{tr}(A)=0,\quad \lambda^*_{\max}(A)=\|A\|>0.
\end{align}
For example, a $6\times 6$ symmetric matrix $A$ with eigenvalues $-2, -1, -1, -1, 2, 3$ satisfies the conditions in (\ref{ex111}), while the matrix $-A$ does not.  The strong $\chi^2$~bound in (\ref{chi2_bound}) simplifies to 
\begin{align}
\label{BF_1}
\Pr(\boldsymbol{\Delta}>t)<1-F_{n}(r)=e^{-\frac{r}{2}}\sum_{i=0}^{\frac{n}{2}-1}\frac{r^i}{2^ii!},
\end{align}
where the last step is due to (\ref{even_n}) and $r$ is given in (\ref{yek}). The $m_\infty$~bound can be written as 
\begin{align}
\label{BF_2}
\Pr(\boldsymbol{\Delta}>t)<\Big(1+\frac{r}{n}\Big)^{\frac{n}{2}}e^{-\frac{r}{2}}=e^{-\frac{r}{2}}\sum_{i=0}^{\frac{n}{2}}\frac{{\frac{n}{2}\choose i}r^i}{n^i}.
\end{align}
 It is easy to see that for $n=2,4$, the $m_\infty$~bound is always looser than the strong $\chi^2$~bound. If $n\geq 6$,  the $m_\infty$~bound is tighter than the strong $\chi^2$~bound if and only if 
\begin{align}
\label{r_dp}
0<r<r''_n,
\end{align}
 where the sequence $r''_n$ satisfies 
 \begin{align}
 \label{excitement}
\lim_{n\to\infty}r''_n=\infty.
\end{align}
We next provide the justification for these claims.  Let $P_n(\cdot)$ be the difference between the two polynomials on the right-hand sides in~(\ref{BF_1}) and (\ref{BF_2}). As $\frac{{\frac{n}{2}\choose i}}{n^i}=\frac{1}{2^ii!}$ for $i=0,1$, we can write 
\begin{align}
\label{not_des}
P_n(r)=\sum_{i=0}^{\frac{n}{2}}\frac{{\frac{n}{2}\choose i}r^i}{n^i}-\sum_{i=0}^{\frac{n}{2}-1}\frac{r^i}{2^ii!}=\frac{r^{\frac{n}{2}}}{n^{\frac{n}{2}}}+\sum_{i=2}^{\frac{n}{2}-1}\Big(\frac{{\frac{n}{2}\choose i}}{n^i}-\frac{1}{2^ii!}\Big)r^i.
\end{align}
The sum on the right-hand side is not void due to $\frac{n}{2}-1\geq \frac{6}{2}-1=2$. Since  $\frac{{\frac{n}{2}\choose i}}{n^i}<\frac{1}{2^ii!}$ for $2\leq i\leq \frac{n}{2}-1$, we conclude that $P_n(\cdot)$ has a positive leading coefficient $\frac{1}{n^{\frac{n}{2}}}$ and all its other remaining coefficients are negative. It is now an easy exercise to show that there exists a (unique) constant $r''_n>0$ such that $P_n(r)<0$ for $0<r<r''_n$ and $P_n(r)>0$ for $r>r''_n$. Next, we prove (\ref{excitement}). If we can show that $P_n(n)<0$ for sufficiently large $n$, it follows that $r''_n>n$ for all such values, which verifies~(\ref{excitement}).  We have 
\begin{align}
\label{st_111}
P_n(n)=2^{\frac{n}{2}}-\sum_{i=0}^{\frac{n}{2}-1}\frac{n^i}{2^ii!}=2^{\frac{n}{2}}-e^{\frac{n}{2}}\sum_{i=0}^{\frac{n}{2}-1}\frac{e^{-\frac{n}{2}}(\frac{n}{2})^i}{i!}\stackrel{}{=} 2^{\frac{n}{2}}-e^{\frac{n}{2}}\mathrm{Pr}(\boldsymbol{x}_n\leq  n/2-1),
\end{align}
where  $\boldsymbol{x}_n\sim \mathrm{Poisson}(\frac{n}{2})$ is a Poisson random variable with parameter $\frac{n}{2}$. Define 
\begin{align}
\boldsymbol{z}_n=\sqrt{\frac{2}{n}}\Big(\boldsymbol{x}_n-\frac{n}{2}\Big)
\end{align}
and let us write   
\begin{align}
\label{st_222}
\Pr(\boldsymbol{x}_n\leq n/2-1)=\Pr(\boldsymbol{x}_n\leq n/2)-\Pr(\boldsymbol{x}_n=n/2)=\Pr(\boldsymbol{z}_n\leq 0)-\Pr(\boldsymbol{x}_n=n/2).
\end{align}
Since $\frac{n}{2}$ is an integer, $\boldsymbol{x}_n$ can be viewed as a sum of $\frac{n}{2}$ independent $\mathrm{Poisson}(1)$ random variables. Because $\mathrm{E}[\boldsymbol{x}_n]=\mathrm{Var}(\boldsymbol{x_n})=\frac{n}{2}$, invoking the Central Limit Theorem shows that the standardized sequence~$\boldsymbol{z}_n$ converges in distribution to a standard normal random variable. Therefore,  
\begin{align}
\label{st_333}
\lim_{n\to\infty}\Pr(\boldsymbol{z}_n\leq 0)=\Phi(0)=0.5,
\end{align}
where $\Phi(\cdot)$ is the standard normal CDF.  Moreover,
\begin{align}
\label{st_444}
\lim_{n\to\infty}\Pr(\boldsymbol{x}_n= n/2)=\lim_{n\to\infty}\frac{e^{-\frac{n}{2}}(\frac{n}{2})^{\frac{n}{2}}}{(\frac{n}{2})!}=\lim_{n\to\infty}\frac{e^{-\frac{n}{2}}(\frac{n}{2})^{\frac{n}{2}}}{\sqrt{\pi n}e^{-\frac{n}{2}}(\frac{n}{2})^{\frac{n}{2}}}=\lim_{n\to\infty}\frac{1}{\sqrt{\pi n}}=0,
\end{align}
where we have used Stirling approximation formula. By~(\ref{st_222}), (\ref{st_333}) and (\ref{st_444}), $\lim_{n\to\infty}\mathrm{Pr}(\boldsymbol{x}_n\leq  \frac{n}{2}-1)=\frac{1}{2}$ and hence, it follows from~(\ref{st_111}) that $\lim_{n\to\infty}P_n(n)=\lim_{n\to\infty} (2^{\frac{n}{2}}-\frac{1}{2}e^{\frac{n}{2}})=-\infty$, i.e., $P_n(n)<0$ for all sufficiently large $n$. 

We close this subsection by observing that for a given positive-semidefinite matrix $A$,  the $m_\infty$~inequality is tighter than LMI for sufficiently large values of the deviation parameter $t$. Recall from Subsection~\ref{2D} that the LM~bound, as given by the right-hand side of (\ref{LM2}), increases with $\|A\|_2$. Since $\|A_2\|\geq \|A\|$, the expression
\begin{align}
\label{LB_LM}
\exp\bigg(-\Big(\frac{\sqrt{\|A\|^2+2\|A\|t}-\|A\|}{2\|A\|}\Big)^2\bigg)=e^{\frac{1}{2}(\sqrt{1+2r}-1)}e^{-\frac{r}{2}}
\end{align}   
is less than or equal to the LM~bound. Consequently, the $m_\infty$~bound $(1+\frac{r}{n})^{\frac{n}{2}}e^{-\frac{r}{2}}$ is tighter than the LM~bound whenever it is smaller than (\ref{LB_LM}). This condition reduces to $(1+\frac{r}{n})^{\frac{n}{2}}<e^{\frac{1}{2}(\sqrt{1+2r}-1)}$, or equivalently, $\sqrt{1+2r}-1>n\ln(1+\frac{r}{n})$, which clearly holds for all sufficiently large $r=\frac{t}{\|A\|}$.

\subsection{Modified $m_\infty$, strong $\chi^2$, and HW inequalities}
\label{2E}
For notational simplicity, we let $r(A)=n$. It is easy to see that the strong $\chi^2$~bound of order $k$ can only become sharper as $k$ increases, for any $t>0$. Let $\lambda_1\le\cdots\le\lambda_k\le\lambda_{k+1}<0$ and~$\boldsymbol{\xi}_1,\dots,\boldsymbol{\xi}_k,\boldsymbol{\xi}_{k+1}$ be standard normal random variables. By~(\ref{chi2_proof}) in Appendix~\ref{A5}, the strong $\chi^2$~bound of order $k$ is equal to $1-\mathbb{E}[F_{n-k}(\boldsymbol{\varrho}_k)]$ where $\boldsymbol{\varrho}_k=\frac{t+\mathrm{tr}(A)-\sum_{i=1}^{k}\lambda_i\boldsymbol{\xi}^2_i}{\lambda_{\max}(A)}$. We have 
\begin{align}
\label{CDF_fun}
F_{n-k}(\boldsymbol{\varrho}_k)\stackrel{(a)}{\le}F_{n-k}(\boldsymbol{\varrho}_{k+1})\stackrel{(b)}{\le}F_{n-(k+1)}(\boldsymbol{\varrho}_{k+1}),
\end{align}
where $(a)$ is due to $\boldsymbol{\varrho}_k\leq \boldsymbol{\varrho}_{k+1}$ and the fact that $F_{n-k}(\cdot)$ is non-decreasing and $(b)$ is due to the fact that if $n_1<n_2$, then $F_{n_2}(x)\leq F_{n_1}(x)$ for all $x\in \mathbb{R}$. Applying $\mathbb{E}[\cdot]$ to both sides in (\ref{CDF_fun}) verifies our claim. 

This consistent monotonic behaviour does not hold for the modified $m_\infty$ and HW bounds. Focusing first on the modified $m_\infty$ bounds, the upper bound in (\ref{call_my_name}) is the product of the exponential term $e^{-\frac{t+b_k}{2a_k}}$ and a polynomial in $t$ whose degree, $\frac{n-k}{2}$, decreases as the order $k$ increases. Since 
\begin{align}
\label{parag_above}
a_{k}=\max_{k+1\leq i\leq n}|\lambda_i|\geq \max_{k+2\leq i\leq n}|\lambda_i|=a_{k+1},
\end{align}
 the exponential term can only decay faster with respect to $t$ for larger $k$. These observations verify that the modified $m_\infty$~inequality becomes sharper as the order increases, provided that $t$ is sufficiently large. 
 
 Next, we examine the modified HW bounds. If $a_k>a_{k+1}$, then the exponential term $e^{-\frac{\kappa(t+b_k)}{a_k}}$ in~(\ref{mod_HWI}), which is the only term dependent on $t$, decays faster with respect to $t$ at order $k+1$ than at order $k$. If $a_k=a_{k+1}$, the ratio of the HW~bound of order $k+1$ over the HW~bound of order $k$ is given by 
 \begin{align}
 \frac{e^{-\frac{\kappa(t+b_{k+1 })}{a_k}}\prod_{i=1}^{k+1} (1-2\kappa\frac{\lambda_i}{a_k})^{-\frac{1}{2}}}{e^{-\frac{\kappa(t+b_{k})}{a_k}}\prod_{i=1}^{k} (1-2\kappa\frac{\lambda_i}{a_k})^{-\frac{1}{2}}}=e^{\kappa\frac{-\lambda_{k+1}}{a_k}} \Big(1-2\kappa\frac{\lambda_{k+1}}{a_k}\Big)^{-\frac{1}{2}}\stackrel{}{>}\frac{1-\kappa\frac{\lambda_{k+1}}{a_k}}{(1-2\kappa\frac{\lambda_{k+1}}{a_k})^{\frac{1}{2}}}>1,
\end{align} 
where we have applied the elementary inequalities $e^x>1+x>\sqrt{1+2x}$ for $x>0$. Hence, the HW~bound of order $k+1$ becomes looser than the HW~bound of order $k$ for every $t>0$.

We now analyze the domains of the modified bounds as functions of the deviation parameter~$t$. The viability condition for the modified $\chi^2$ bounds in (\ref{poloi_4}) does not depend on the order~$k$. Hence, the domain is the identical interval $[-\mathrm{tr}(A),\infty)$ for all orders.  

For the modified $m_\infty$ bound of order $k$, the viability condition in (\ref{poloi_3}) yields the domain $[-b_k,\infty)$. Because $-b_k$ increases with $k$, this domain shrinks as the order $k$ increases.  

Similarly, the viability condition in (\ref{poloi_4}) for the modified HW bound of order $k$ results in the domain $\big[\frac{\sum_{i=k+1}^n\lambda_i^2}{a_k}-b_k,\infty\big)$. The lower bound $\frac{\sum_{i=k+1}^n\lambda_i^2}{a_k}-b_k$ is non-decreasing in $k$. To verify this property, let $\lambda_1\leq \cdots\leq \lambda_k\leq \lambda_{k+1}<0$ be the $k+1$ smallest eigenvalues. It follows that 
\begin{align}
\label{order_compare}
& \frac{\sum_{i=k+1}^n\lambda_i^2}{a_k}-b_k\le\frac{\sum_{i=k+2}^n\lambda_i^2}{a_{k+1}}-b_{k+1}\notag\\
 \iff \hskip0.5cm& b_{k+1}-b_k\le\frac{\sum_{i=k+2}^n\lambda_i^2}{a_{k+1}}-\frac{\sum_{i=k+1}^n\lambda_i^2}{a_k}\notag\\
  \iff \hskip0.5cm& \lambda_{k+1}\le\frac{1}{a_ka_{k+1}}\Big(a_k\sum_{i=k+2}^n\lambda_i^2-a_{k+1}\sum_{i=k+1}^n\lambda_i^2\Big)\notag\\
    \iff \hskip0.5cm& \lambda_{k+1}\le\frac{1}{a_ka_{k+1}}\Big(-a_{k+1}\lambda_{k+1}^2+(a_k-a_{k+1})\sum_{i=k+2}^n\lambda_i^2\Big)\notag\\
     \iff \hskip0.5cm& \lambda_{k+1}+\frac{\lambda_{k+1}^2}{a_k}\le\frac{a_k-a_{k+1}}{a_ka_{k+1}}\sum_{i=k+2}^n\lambda_i^2.
 \end{align} 
 The final inequality in (\ref{order_compare}) always holds. Indeed, the right-hand side satisfies $\frac{a_k-a_{k+1}}{a_ka_{k+1}}\sum_{i=k+2}^n\lambda_i^2\geq0$ due to (\ref{parag_above}), whereas the left-hand side satisfies $\lambda_{k+1}+\frac{\lambda_{k+1}^2}{a_k}=\frac{\lambda_{k+1}}{a_k}(a_k-|\lambda_{k+1}|)\leq 0$ because $|\lambda_{k+1}|\leq a_k$. Consequently,  the domain of the modified HW bound of order $k$ either remains unchanged or shrinks as the order $k$ increases.
 
 It is also immediate that each of the best viable modified $m_\infty$, $\chi^2$ and HW bounds is a right-continuous function with respect to $t$. This property follows from the simple observation that if a collection of functions $f_i:[a_i,\infty)\to\mathbb{R}$ is continuous on its respective domain for each $i=1,\dots, m$, then their pointwise minimum $\min_{1\leq i\leq m}f_i$ is right-continuous on its domain $[\min_{1\leq i\leq m}a_i,\infty)$. This property will be utilized in Subsection~\ref{3_A}.


We end this subsection with a numerical example. 

\begin{figure*}[t]
\centering
\subfigure[]{
\includegraphics[scale=0.35]{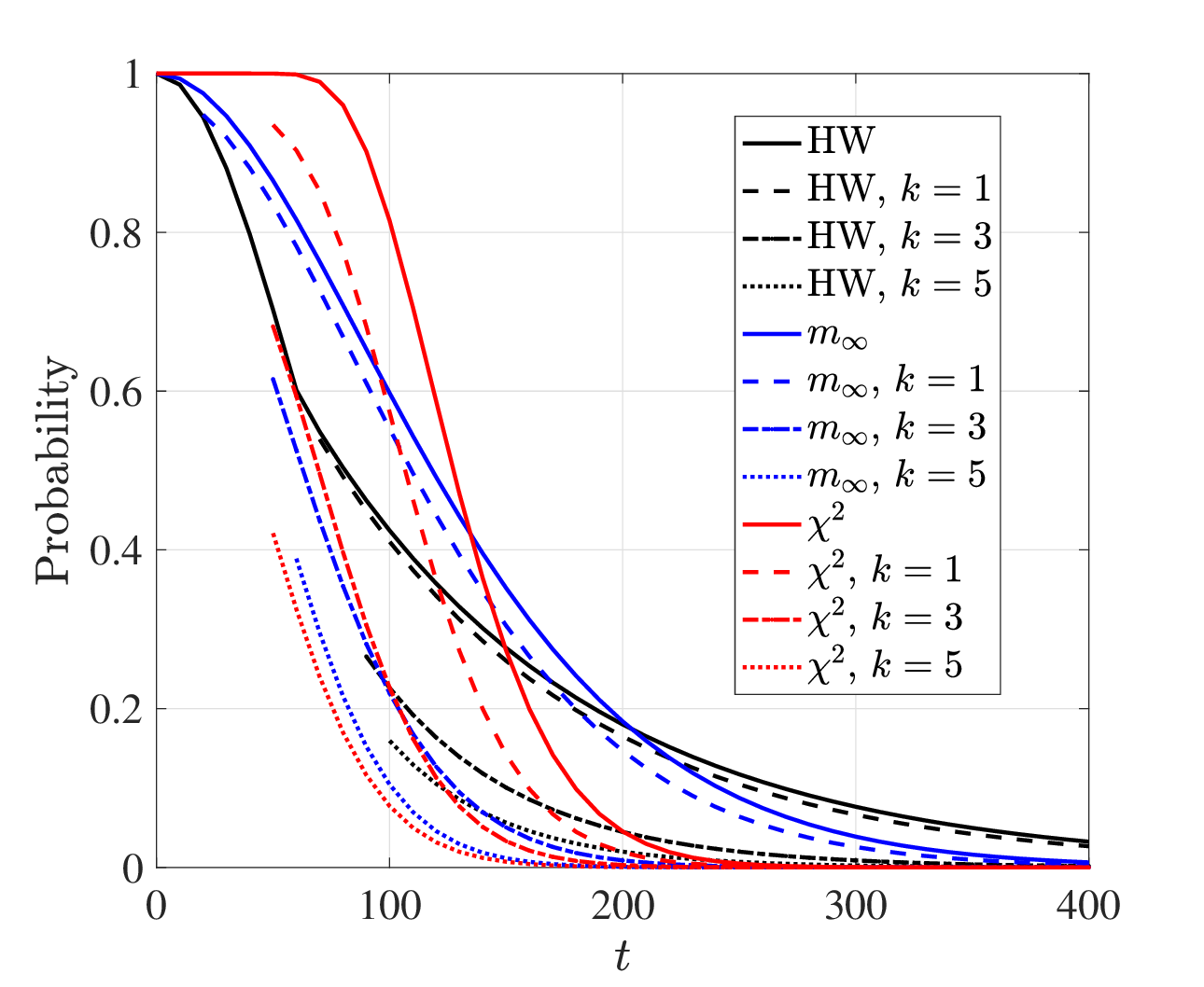}
\label{mod_bounds_figa}
}
\subfigure[]{
\includegraphics[scale=0.35]{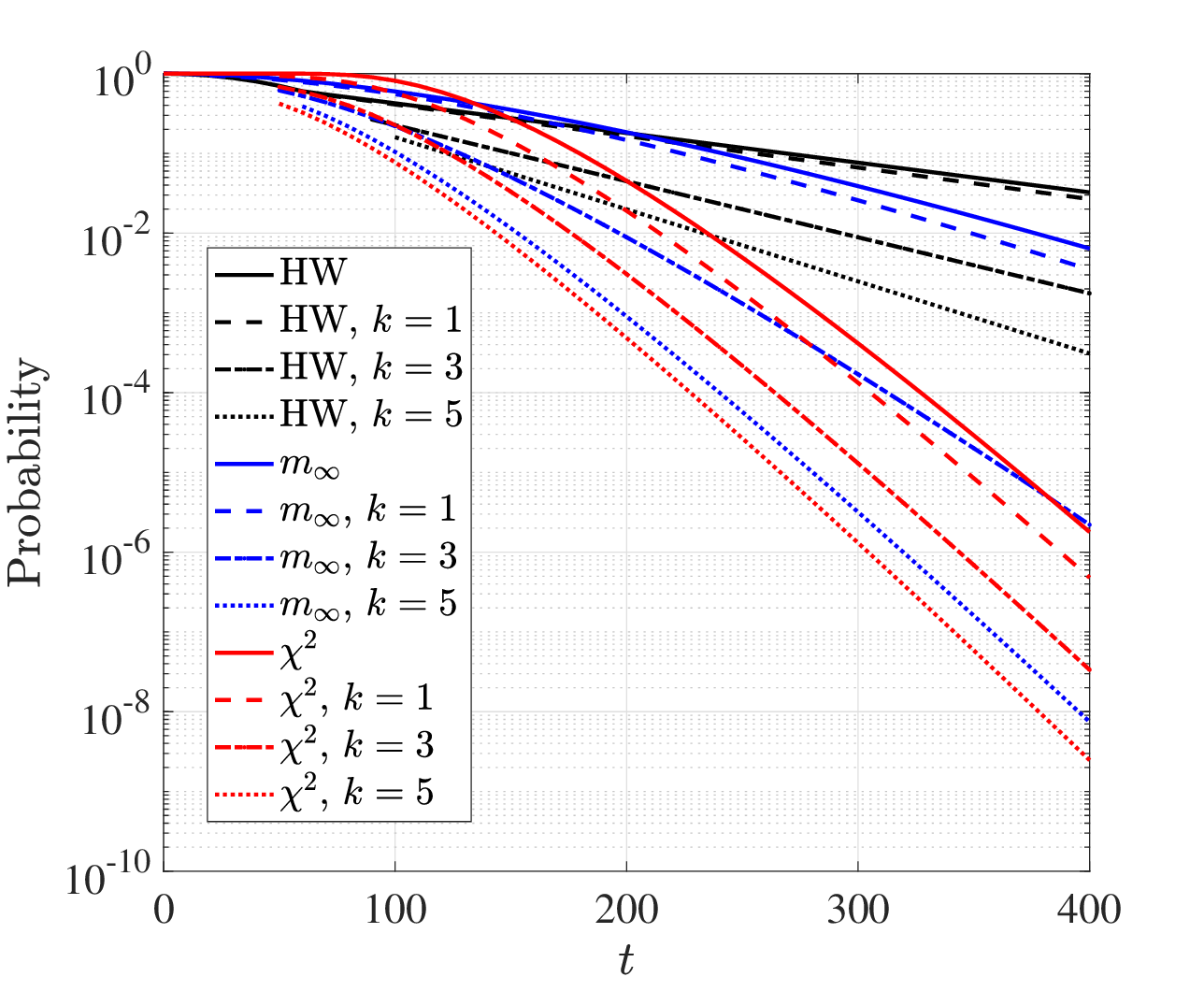}
\label{mod_bounds_figb}
}
\label{fig:subfigureExample}
\caption[Optional caption for list of figures]{Consider a $13\times 13$ symmetric matrix~$A$ with eigenvalues -17, -16, -11, -9, -6, -6, -6, 2, 3, 4, 7, 7, 7.  Fig.~\ref{mod_bounds_figa} plots the original HW, $m_\infty$ and strong $\chi^2$ bounds, along with their modified versions of orders $k=1,3,5$ for $0\le t\leq 400$ in steps of $10$. Fig.~\ref{mod_bounds_figb} displays the same curves using a logarithmic scale for the vertical axis.}
\label{mod_bounds_fig}
\end{figure*}

\textbf{Example}- \textcolor{Black}{Consider a $13\times 13$ symmetric matrix~$A$ with eigenvalues -17, -16, -11, -9, -6, -6, -6, 2, 3, 4, 7, 7, 7.  Fig.~\ref{mod_bounds_figa} plots the original HW, $m_\infty$ and strong $\chi^2$ bounds, along with their modified versions of orders $k=1,3,5$ for $0\le t\leq 400$ in steps of $10$. Fig.~\ref{mod_bounds_figb} displays the same curves using a logarithmic scale for the vertical axis.~$\diamond$}  

\subsection{ RK inequalities  }
\label{2F}
\textcolor{Black}{In this subsection, we study the approximation gaps in (\ref{RK_error_minus}) and (\ref{RK_error_plus}). First, we examine the numbers $\delta_k$ in (\ref{KSD}). An explicit formula for these values is given by 
\begin{align}
\label{closed_form_KSD}
 \delta_k=F_k(x_k)-F_{k+1}(x_k),
\end{align}
where
\begin{align}
\label{KSD_x_k}
 x_k=2\bigg(\frac{\Gamma(\frac{k+1}{2})}{\Gamma(\frac{k}{2})}\bigg)^2.
\end{align}
To see this, let $f_k(x)=F_{k}(x)-F_{k+1}(x)$ for $x\in \mathbb{R}$. We have $f_k(x)=0$ for all $x\le 0$. For every $x>0$, 
\begin{align}
\label{}
 f_k'(x)&=F'_{k}(x)-F'_{k+1}(x)\notag\\&=p_{\chi^2_k}(x)-p_{\chi^2_{k+1}}(x)\notag\\&=\frac{x^{\frac{k}{2}-1}e^{-\frac{x}{2}}}{2^{\frac{k}{2}}\Gamma(\frac{k}{2})}-\frac{x^{\frac{k+1}{2}-1}e^{-\frac{x}{2}}}{2^{\frac{k+1}{2}}\Gamma(\frac{k+1}{2})}\notag\\
 &=\frac{x^{\frac{k}{2}-1}e^{-\frac{x}{2}}}{2^{\frac{k}{2}}}\Big(\frac{1}{\Gamma(\frac{k}{2})}-\frac{\sqrt{x}}{\sqrt{2}\,\Gamma(\frac{k+1}{2})}\Big).
\end{align}
Then $f_k'(x)>0$ for $0<x<x_k$ and $f_k'(x)<0$ for $x>x_k$ where $x_k$ is given in (\ref{KSD_x_k}). It follows from the first derivative test that $\max_{x\in \mathbb{R}}f_k(x)=f_k(x_k)$, as promised.} 

\textcolor{Black}{Let us recall Gautschi's inequality~\cite{Qi} on the ratio of Gamma functions, which states that
\begin{align}
\label{}
x^{1-a} < \frac{\Gamma(x+1)}{\Gamma(x+a)}<(x+1)^{1-a},
\end{align}
for every $x\ge0$ and $0<a<1$. Setting $x=\frac{k-1}{2}$ and $a=\frac{1}{2}$ yields
\begin{align}
\label{}
\Big(\frac{k-1}{2}\Big)^{\frac{1}{2}}<\frac{\Gamma(\frac{k+1}{2})}{\Gamma(\frac{k}{2})}<\Big(\frac{k+1}{2}\Big)^{\frac{1}{2}}.
\end{align}
Consequently, the numbers $x_k$ in (\ref{KSD_x_k}) satisfy 
\begin{align}
\label{bounds_x_k}
  k-1<x_k<k+1,\,\,\,k\ge1.
\end{align}
Using (\ref{bounds_x_k}) together with the Berry-Esseen Central Limit Theorem~\cite{Durrett}, one can show that $\lim_{k\to\infty}\delta_k=0$. We conjecture that $\delta_k$ is monotonically decreasing with respect to $k$. Hence, 
\begin{align}
\label{}
\forall k\ge1,\,\,\,\delta_k\leq \delta_1=F_1(x_1)-F_2(x_1)=\frac{\gamma(\frac{1}{2},\frac{1}{\pi})}{\Gamma(\frac{1}{2})}-\frac{\gamma(1,\frac{1}{\pi})}{\Gamma(1)}\approx 0.3024,
\end{align}   
where the penultimate step is due to (\ref{exact_chi2_CDF}) and $x_1=2\big(\frac{\Gamma(1)}{\Gamma(\frac{1}{2})}\big)^2=\frac{2}{\pi}$.}

\textcolor{Black}{We continue by summarizing several known facts about the Ruben-Kotz expansion in (\ref{RK_expansion}). The coefficients $c^{\text{RK}}_l=c^{\text{RK}}_l(\lambda_1,\cdots, \lambda_n; s)$, defined through the recursion in (\ref{RK_recursion}), do not admit closed-form expressions unless the numbers $\lambda_1,\dots, \lambda_n$ appear in pairs. Kotz et al.~\cite{Kotz1} provide the alternative key formulation 
\begin{align}
\label{Key_Kotz}
 c^{\text{RK}}_l=c^{\text{RK}}_0\frac{\mathbb{E}[\boldsymbol{Q}^l]}{2^ll!},\,\,\,l\ge1.
\end{align}
Here, the quadratic form $\boldsymbol{Q}$ is given by  
\begin{align}
\label{}
  \boldsymbol{Q}=\sum_{i=1}^n\Big(1-\frac{s}{\lambda_i}\Big)\boldsymbol{\xi}_i^2,
\end{align}
where $\boldsymbol{\xi}_1,\dots, \boldsymbol{\xi}_n$ are independent standard normal random variables. Using~(\ref{Key_Kotz}), Kotz et al.~\cite{Kotz1} derive the bounds 
\begin{align}
\label{Kotz_Bound}
 | c^{\text{RK}}_l|\leq c^{\text{RK}}_0 {\frac{n}{2}+l-1\choose l}(\max_{1\le i\le n}\rho_i)^l,\,\,\,l\ge1,
\end{align} 
where ${\frac{n}{2}+l-1\choose l}=\frac{\Gamma(\frac{n}{2}+l)}{\Gamma(\frac{n}{2})l!}$ and $\rho_i=|1-\frac{s}{\lambda_i}|$ for $1\le i\le n$. One can also verify (\ref{Kotz_Bound}) by applying strong induction to the recursion in~(\ref{RK_recursion}). }

\textcolor{Black}{An immediate consequence of (\ref{Kotz_Bound}) is that if $\max_{1\le i\le n}\rho_i<1$, or equivalently, $0<s<2\lambda_{\min}(A)$, the sequence of Ruben-Kotz coefficients $(c^{\text{RK}}_l)_{l\ge0}$ is absolutely summable. In fact,
\begin{align}
\label{}
  \sum_{l=0}^\infty|c^{\text{RK}}_l|\leq c^{\text{RK}}_0\sum_{l=0}^\infty  {\frac{n}{2}+l-1\choose l}(\max_{1\le i\le n}\rho_i)^l=c^{\text{RK}}_0\big(1-\max_{1\leq i\leq n}\rho_i\big)^{-\frac{n}{2}}<\infty,
\end{align} 
where we have used the binomial series $(1-x)^{-a}=\sum_{l=0}^\infty{a+l-1\choose l}x^l$ for $a\in \mathbb{R}$ and $|x|<1$. Due to the absolute summability of $(c^{\text{RK}}_l)_{l\ge0}$, invoking the Weierstrass M-test~\cite{Rudin} shows that the Ruben-Kotz series in~(\ref{RK_expansion}) converges uniformly on $\mathbb{R}$. This uniform convergence warrants an exchange of limits, yielding 
\begin{align}
\label{}
 1=\lim_{z\to\infty}F_{\underline{\boldsymbol{x}}^{\mathsf{T}}A\underline{\boldsymbol{x}}}(z)=\lim_{z\to\infty}\sum_{l=0}^\infty c^{\text{RK}}_{l}F_{n+2l}\Big(\frac{z}{s}\Big)=\sum_{l=0}^\infty c^{\text{RK}}_{l}\lim_{z\to\infty}F_{n+2l}\Big(\frac{z}{s}\Big)=\sum_{l=0}^\infty c^{\text{RK}}_{l},
\end{align}
which proves $\sum_{l=0}^\infty c^{\text{RK}}_{l}=1$. More recently, Ropokis et al.~\cite{Ropokis} derived significantly tighter bounds than those in (\ref{Kotz_Bound}).}

\textcolor{Black}{ Returning now to the setup in Proposition~\ref{prop_11_last}, we note that since $s^-$ and $s^+$ satisfy the conditions in~(\ref{s_pm_ineq}), it holds that $\sum_{l=0}^\infty c^{\text{RK}^\pm}_l=1$.  If either $t+\mathrm{tr}(A)<0$ with even $u$, or $t+\mathrm{tr}(A)\ge0$ with even $v$, then the approximation gaps $\mathsf{Gap}^{\pm}=1-\sum_{l=0}^{\infty}c_l^{\text{RK}^\pm}$, given in~(\ref{RK_error_minus}) and (\ref{RK_error_plus}),  tend to zero as $L$ grows large, yielding the exact value for $\Pr(\boldsymbol{\Delta}>t)$. In other words, RK inequalities are tight as $L\to\infty$.  On the other hand, if either $t+\mathrm{tr}(A)<0$ with odd $u$, or $t+\mathrm{tr}(A)\ge0$ with odd $v$, the gaps $\mathsf{Gap}^{\pm}=1-\sum_{l=0}^L c_l^{\text{RK}^\pm}+\sum_{l=0}^{L}w^{\pm}_{l}c_l^{\text{RK}^\pm}$ saturate at $\sum_{l=0}^\infty w_l^{\pm}c_l^{\text{RK}^\pm}$ as $L$ grows large, where $w_l^{-}=\delta_{u+2l-1}$ and $w_l^{+}=\delta_{v+2l}$. Accepting that $(\delta_k)_{k\ge1}$ is monotonically decreasing, it follows that $\sum_{l=0}^\infty w_l^{\pm}c_l^{\text{RK}^\pm}\leq w_0^{\pm}$. Hence, regardless of the tail parameter $t$, the limiting value of the error never exceeds $\max\{\delta_{u-1},\delta_v\}$.   }

\begin{figure*}[t]
\centering
\subfigure[]{
\includegraphics[scale=0.35]{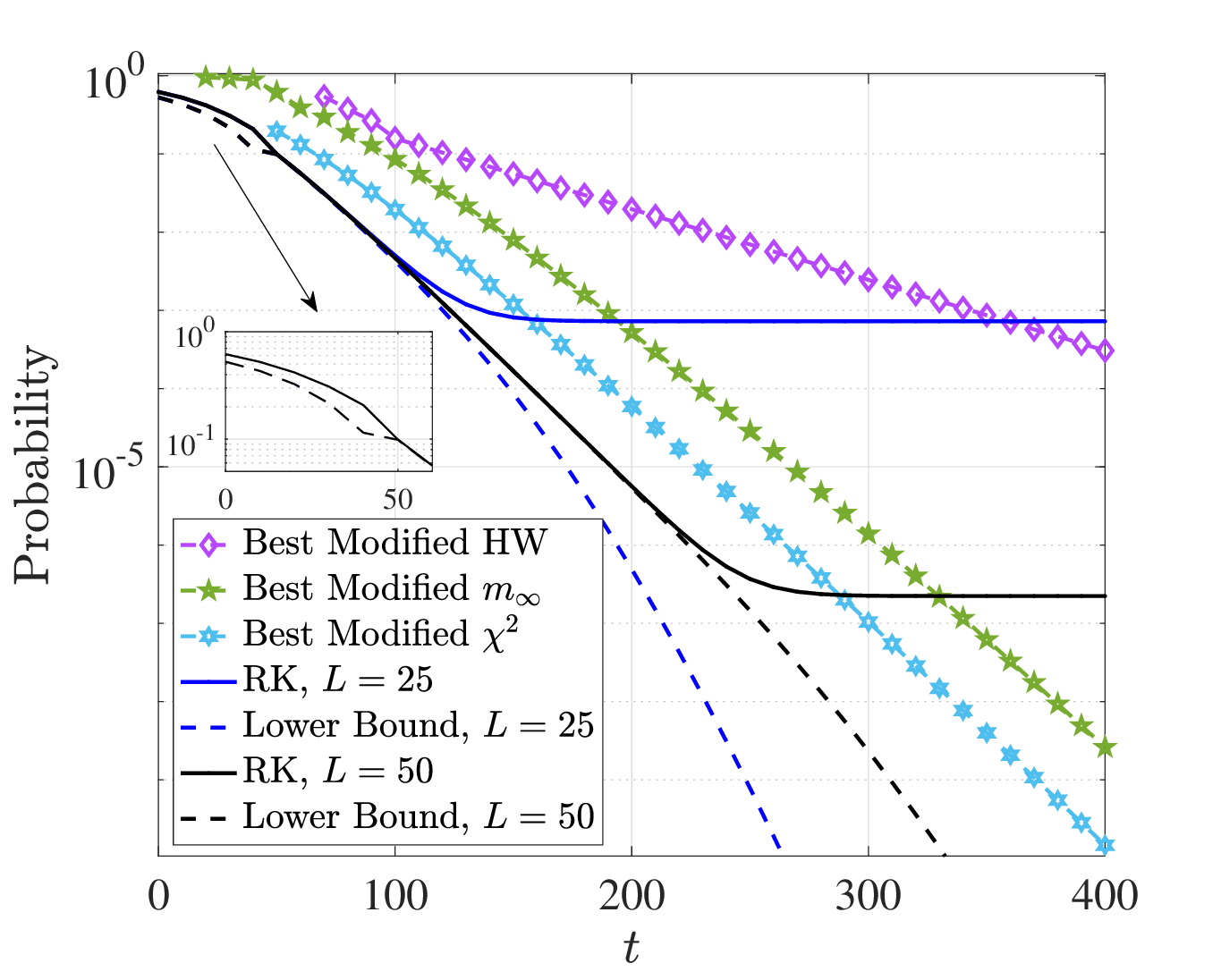}
\label{RK_figa}
}
\subfigure[]{
\includegraphics[scale=0.35]{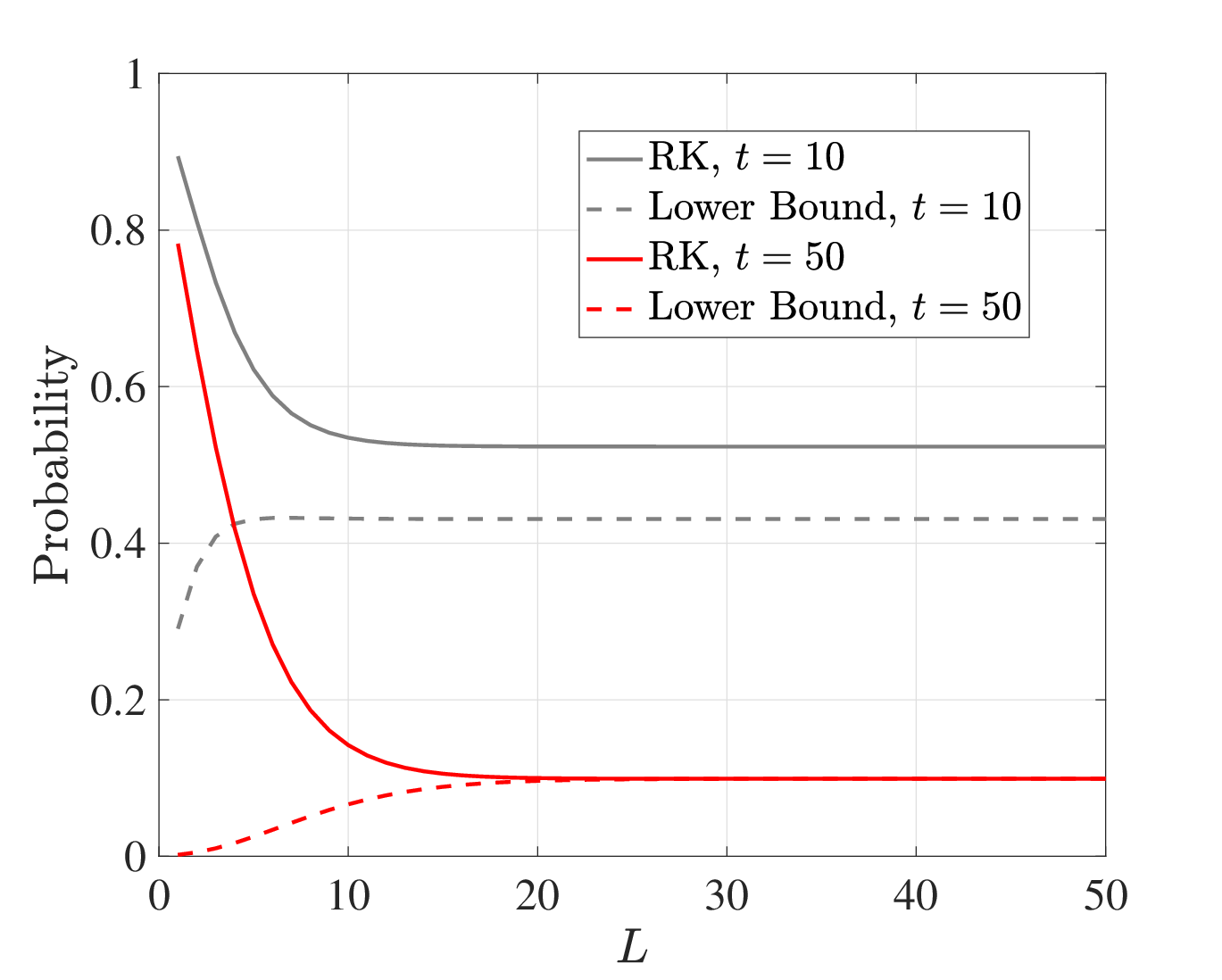}
\label{RK_figb}
}
\label{fig:subfigureExample}
\caption[Optional caption for list of figures]{Consider the $13\times 13$ matrix $A$ in Fig.~\ref{mod_bounds_fig}. Panel~(a) presents the plots for the RK bounds alongside their corresponding lower bounds (the RK bounds minus the gaps) versus $0\leq t\leq 400$ for $L=25$ and $L=50$. Panel~(b) plots the RK bounds and their corresponding lower bounds versus $1\leq L\leq 50$ for $t=10$ and $t=50$. Since $\mathrm{tr}(A)=-41$ and $u=7$ is odd, the approximation gap $\mathsf{Gap}^-$ for $t<41$ saturates at $\sum_{l=0}^\infty \delta_{6+2l}c_l^{\text{RK}^-}\le \delta_6\approx 0.1177$. On the other hand, since $v=6$ is even, the approximation gap $\mathsf{Gap}^+$ for $t\geq41$ vanishes asymptotically as $L$ grows large.      }
\label{RK_fig}
\end{figure*}

\textcolor{Black}{\textbf{Example}- Consider the $13\times 13$ matrix $A$ in Fig.~\ref{mod_bounds_fig}. Here, $u=7$, $v=6$ and $\mathrm{tr}(A)=-41$. The RK bound in (\ref{RK_bound_minus}) applies for $t< 41$. Since $u$ is odd, the bound on the approximation gap in (\ref{RK_error_minus}) does not vanish as $L$ increases. It converges to $\lim_{L\to\infty}\mathsf{Gap}^-=\sum_{l=0}^\infty \delta_{6+2l}c_l^{\text{RK}^-}\leq \delta_6\approx 0.1177$.  The RK bound in (\ref{RK_bound_plus}) applies for $t\geq 41$. Because $v$ is even, the bound on the approximation gap in (\ref{RK_error_plus}) vanishes asymptotically, i.e., $\lim_{L\to\infty}\mathsf{Gap}^+=0$.}

 \textcolor{Black}{Fig.~\ref{RK_figa} plots the RK bounds alongside their corresponding lower bounds (upper bounds minus gaps) for $0\leq t\leq 400$ in steps of $10$, evaluated for $L=25$ and $L=50$ separately, with the vertical axis on a logarithmic scale.  For comparison,  we have also included the curves for the best modified HW, $m_\infty$ and strong $\chi^2$ bounds. Finally, Fig.~\ref{RK_figb} plots the RK bounds and the corresponding lower bounds against $1\leq L\leq 50$ for $t=10$ and $t=50$ separately, where the vertical axis is on a linear scale. For $t=10$,   
\begin{align}
\label{}
 \mathsf{Gap}^-\big|_{L=25}\approx  \mathsf{Gap}^-\big|_{L=50}\approx 0.0925,
\end{align}
while for $t=50$, 
\begin{align}
\label{}
 \mathsf{Gap}^+\big|_{L=25}\approx 7.25\times 10^{-4},\quad \mathsf{Gap}^+\big|_{L=50}\approx 2.23\times 10^{-7}.
\end{align}
}
 

\section{Applications}
\label{sec3}
This section evaluates the performance of the bounds from Section I across two applications in statistical signal processing and wireless communications. Unlike the framework in Section II, where the tail parameter $t$ and the matrix $A$ are treated as independent, these applications feature a dependency between $t$ and $A$. The primary takeaway is that while the standard $m_{\infty }$, strong $\chi ^{2}$, and HW bounds may underperform, their modified versions in Propositions~\ref{prop_17},~\ref{prop_9} and~\ref{prop_8} can provide substantial performance gains.
  \subsection{Statistical Signal Processing}
  \label{3_A}
 Let $\boldsymbol{y}_1,\boldsymbol{y}_2,\dots, \boldsymbol{y}_n$ be samples from a zero-mean Gaussian process. We consider the statistical test of \textit{simple} hypotheses concerning the covariance matrix of $\underline{\boldsymbol{y}}=(\boldsymbol{y}_1,\dots, \boldsymbol{y}_n)^T$ given by 
 \begin{align}
 \label{H0_H1}H_0: \mathrm{Cov}(\underline{\boldsymbol{y}})=C_{0},\quad H_1: \mathrm{Cov}(\underline{\boldsymbol{y}})=C_{1}, 
 \end{align} 
 where $C_{0}$ and $C_{1}$ are given positive-definite matrices. In practice, such testing is often performed in an online fashion as new data becomes available by using variants of Wald's Sequential Probability Ratio Test~(SPRT), where the decision rule is designed to minimize the expected number of samples subject to given thresholds on the Type~I and Type~II error probabilities~\cite{Wald, Blostein}. Here, however, we assume the sample size $n$ is fixed in advance and hence adopt the Neyman-Pearson paradigm. The log-likelihood-ratio simplifies as 
\begin{align}
\label{llr}
\textsf{LLR}(\underline{y})=\ln\frac{p_{\underline{\boldsymbol{y}}}(\underline{y};H_1)}{p_{\underline{\boldsymbol{y}}}(\underline{y};H_0)}=\frac{1}{2}\underline{y}^{\mathsf{T}}(C^{-1}_0-C_1^{-1})\underline{y}+\frac{1}{2}\ln\det(C_1^{-1}C_0).
\end{align}
At a significance level of $\alpha$, the optimal decision rule votes in favour of $H_1$ if and only if 
  \begin{align}
\textsf{LLR}(\underline{y})>\gamma,
\end{align} 
where $\gamma\in \mathbb{R}$ is a threshold. Let $p_{1|0}$ and $p_{0|1}$ be the Type~I and Type~II error probabilities, i.e., $p_{1|0}$ is the probability of wrongfully voting in favour of $H_1$ and $p_{0|1}$ is the probability of wrongfully voting in favour of $H_0$. The threshold $\gamma$ is then chosen such that $p_{0|1}$ is minimized subject to $p_{1|0}\le\alpha$.  We have
\begin{align}
p_{1|0}=\Pr(\textsf{LLR}(\underline{\boldsymbol{y}})>\gamma|H_0)=\Pr\big(\underline{\boldsymbol{y}}^{\mathsf{T}}(C^{-1}_0-C_1^{-1})\underline{\boldsymbol{y}}>2\gamma-\ln\det(C_1^{-1}C_0)|H_0\big).
\end{align}
Let us represent $\underline{\boldsymbol{y}}$ as 
\begin{align}
\label{whitening}
\underline{\boldsymbol{y}}=\left\{\begin{array}{cc}
    C_0^{1/2}\underline{\boldsymbol{x}}  & \textrm{under $H_0$}   \\
    C_1^{1/2}\underline{\boldsymbol{x}}  &   \textrm{under $H_1$}  
\end{array}\right.,\quad \underline{\boldsymbol{x}}\sim N(\underline{0},I_n). 
\end{align}
Then 
\begin{align}
\label{yek_mah}
p_{1|0} &=\Pr\big((C_0^{1/2}\underline{\boldsymbol{x}})^{\mathsf{T}}(C^{-1}_0-C_1^{-1})C_0^{1/2}\underline{\boldsymbol{x}}>2\gamma-\ln\det(C_1^{-1}C_0)\big)\notag\\
&=\Pr\big(\underline{\boldsymbol{x}}^{\mathsf{T}}(I_n-C_0^{1/2}C_1^{-1}C_0^{1/2})\underline{\boldsymbol{x}}>2\gamma-\ln\det(C_1^{-1}C_0)\big).
\end{align}
Define 
\begin{align}
\label{A10}
A_{1|0}=I_n-C_0^{1/2}C_1^{-1}C_0^{1/2}.
\end{align} 
Subtracting $\mathrm{tr}(A_{1|0})$ from both sides of the inequality in (\ref{yek_mah}), 
\begin{align}
\label{TYPE_I}
p_{1|0}=\Pr\big(\underline{\boldsymbol{x}}^{\mathsf{T}}A_{1|0}\underline{\boldsymbol{x}}-\mathrm{tr}(A_{1|0})>c_{1|0}+2\gamma\big),
\end{align}
where the constant $c_{1|0}$ is given by 
\begin{align}
\label{c_1_0}
c_{1|0}=-\ln\det(C_1^{-1}C_0)-\mathrm{tr}(A_{1|0}).
\end{align}
We will denote the deviation parameter in (\ref{TYPE_I}) by $t_{1|0}$, i.e., 
\begin{align}
t_{1|0}=c_{1|0}+2\gamma.
\end{align}
It is easy to see that $c_{1|0}>0$. Let $\mu_1,\dots, \mu_n$ be the eigenvalues for the matrix~$C_1^{-1}C_0$. Then $\mu_i$ are all positive and we have  
\begin{align}
\label{l_back}
c_{1|0} &=-\ln\det(C_1^{-1}C_0)-\mathrm{tr}(I_n-C_0^{1/2}C_1^{-1}C_0^{1/2})\notag\\&=-\ln\det(C_1^{-1}C_0)-\mathrm{tr}(I_n-C_1^{-1}C^{1/2}_0C^{1/2}_0)\notag\\
&=-\ln\det(C_1^{-1}C_0)-\mathrm{tr}(I_n-C_1^{-1}C_0)\notag\\&=\sum_{i=1}^n(\mu_i-1-\ln \mu_i)>0,
\end{align}
where the second step is due to $\mathrm{tr}(AB)=\mathrm{tr}(BA)$ for matrices $A, B$ and the last step is due to $\ln x<x-1$ for all positive $x\neq 1$ and the fact that at least one of $\mu_1,\dots, \mu_n$ is not equal to $1$, as otherwise $C_0=C_1$. Note that the matrix $C_1^{-1}C_0$ is not necessarily positive-definite, however, it is diagonalizable. If all its eigenvalues are $1$, then $C_1^{-1}C_0=I_n$ and hence, $C_0=C_1$.
 
Similarly, 
\begin{align}
\label{TYPE_II}
p_{0|1}=\Pr\big(\underline{\boldsymbol{x}}^{\mathsf{T}}A_{0|1}\underline{\boldsymbol{x}}-\mathrm{tr}(A_{0|1})>t_{0|1}\big),
\end{align}
where 
\begin{align}
\label{A01}
A_{0|1}=I_n-C_1^{1/2}C_0^{-1}C_1^{1/2},
\end{align}
\begin{align}
\label{fifa1}
t_{0|1}=c_{0|1}-2\gamma
\end{align}
and 
\begin{align}
\label{fifa2}
c_{0|1}=-\ln\det(C_0^{-1}C_1)-\mathrm{tr}(A_{0|1})=\sum_{i=1}^n\Big(\frac{1}{\mu_i}-1-\ln\frac{1}{ \mu_i}\Big)>0.
\end{align}
Both probabilities $p_{1|0}$ and $p_{0|1}$ depend on $C_0$ and $C_1$ through the eigenvalues $\mu_1,\dots,\mu_n$ of $C_1^{-1}C_0$. To see this, let us write $\underline{\boldsymbol{x}}^{\mathsf{T}}A_{1|0}\underline{\boldsymbol{x}}=\sum_{i=1}^n \lambda_i \boldsymbol{\xi}_i^2$ where $\lambda_i$ are the eigenvalues of $A_{1|0}$ and $\boldsymbol{\xi}_i$ are independent standard normal random variables.Denoting the eigenvalues of $C_0^{1/2}C_1^{-1}C_0^{1/2}$ by $\tilde{\lambda}_i$, it is clear that $\lambda_i=1-\tilde{\lambda}_i$. Since for every two invertible $n\times n$ matrices $A, B$, the eigenvalues for $AB$ and $BA$ are identical, the numbers $\tilde{\lambda}_i$ are exactly the eigenvalues for $C_1^{-1}C_0^{1/2}C_0^{1/2}=C_1^{-1}C_0$, i.e., $\tilde{\lambda}_i=\mu_i$. Hence, $\lambda_i=1-\mu_i$, proving our claim. Similarly, the eigenvalues of $A_{0|1}$ are $1-\frac{1}{\mu_i}$. 

\begin{figure*}[t]
\centering
\subfigure[]{
\includegraphics[scale=0.4]{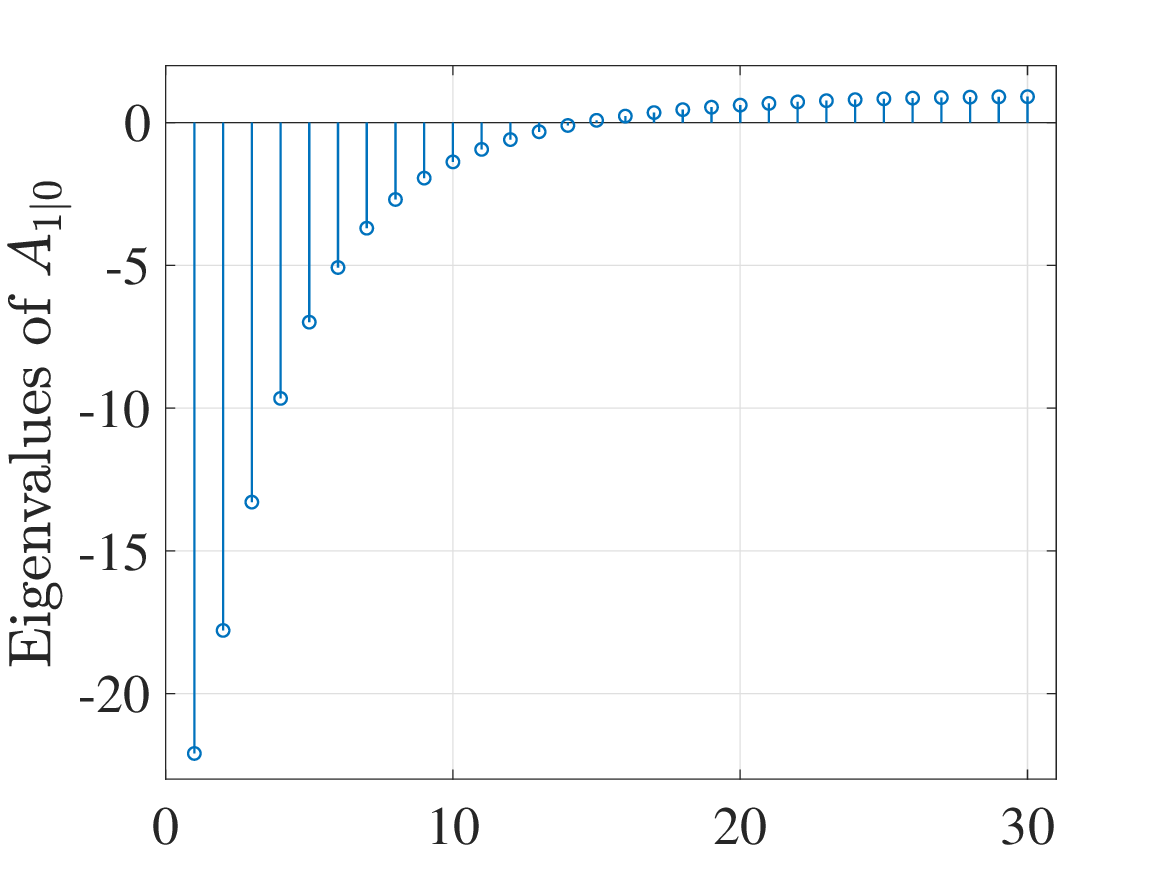}
\label{fig_eig1}
}
\subfigure[]{
\includegraphics[scale=0.4]{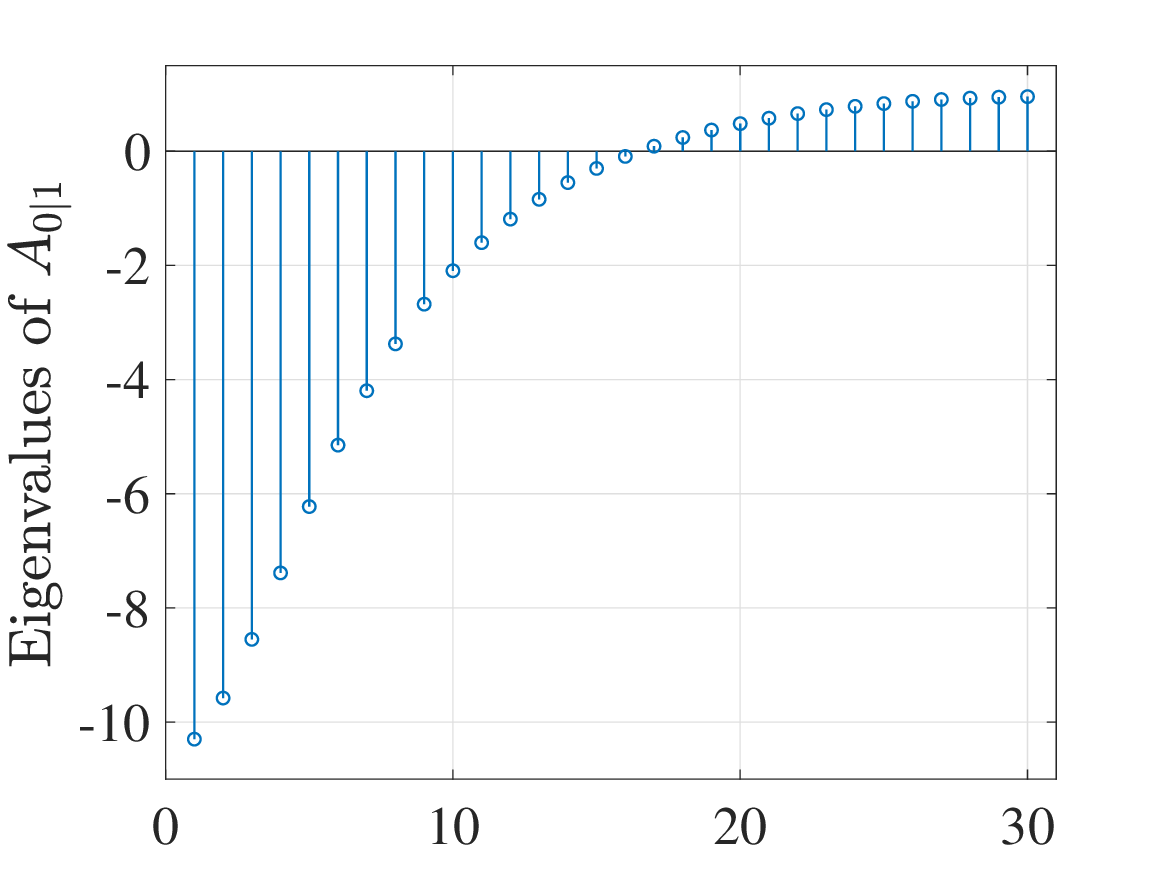}
\label{fig_eig2}
}
\label{fig:subfigureExample}
\caption[Optional caption for list of figures]{Let $\theta_0=-0.5$, $\theta_1=0.7$ and $n=30$. Panel~(a) plots the sorted eigenvalues $1-\mu_i$ of $A_{1|0}$, while Panel~(b) shows the sorted eigenvalues $1-\frac{1}{\mu_i}$ of $A_{0|1}$. Both matrices exhibit several negative eigenvalues that are significantly larger in magnitude than the others, particularly in the case of \(A_{1|0}\).   }
\label{fig_eig}
\end{figure*}

\label{PL_1} We list a few remarks: 
\begin{enumerate}
\item The eigenvalues of $A_{1|0}$ and $A_{0|1}$ are $1-\mu_i$ and $1-\frac{1}{\mu_i}$, respectively, for $i=1,\dots, n$, as established in the previous paragraph.   
\item If $t_{1|0}\ge0$, i.e., 
\begin{align}
\label{normal_bounds}
\gamma\geq -\frac{1}{2}c_{1|0},
\end{align}
 then one can readily apply HWI and the $m_\infty$~inequality to produce two upper bounds on $p_{1|0}$ in~(\ref{TYPE_I}).  Similarly, if $t_{0|1}\ge0$, i.e., 
\begin{align}
\label{normal_bounds_1}
\gamma\le\frac{1}{2}c_{0|1},
\end{align}
 then one can apply HWI and the $m_\infty$~inequality to produce two upper bounds on $p_{0|1}$ in~(\ref{TYPE_II}). A third upper bound which does not require positivity of $t_{1|0}$ or $t_{0|1}$ is  the strong $\chi^2$~bound. The symmetric matrices $A_{1|0}$ and $A_{0|1}$ are not necessarily positive-semidefinite and hence, LMI and its variants in Section~\ref{sec2} may not be applicable. 
   \item Recall that condition (\ref{poloi_4}) in Proposition~\ref{prop_9} must hold for the modified strong $\chi^2$~bound of any order to be viable. We have $t_{1|0}+\mathrm{tr}(A_{1|0})=2\gamma+c_{1|0}+\mathrm{tr}(A_{1|0})=2\gamma-\ln\prod_{i=1}^n\mu_i$ and $t_{0|1}+\mathrm{tr}(A_{0|1})=-2\gamma+c_{0|1}+\mathrm{tr}(A_{0|1})=-2\gamma+\ln\prod_{i=1}^n\mu_i$. If $2\gamma\neq \ln \prod_{i=1}^n\mu_i$, one of $t_{1|0}+\mathrm{tr}(A_{1|0})$ and $t_{0|1}+\mathrm{tr}(A_{0|1})$ is negative and therefore, the modified $\chi^2$ bound for one of $p_{1|0}$ and $p_{0|1}$ is not viable. 
  \item 
In scenarios where all \(\mu_1,\dots, \mu_n\) are distinct and include both positive and negative values, as in the example below, the common practice is to compute approximate values for \(p_{1|{}0}\) and \(p_{0|{}1}\) via Monte Carlo simulations. In the special case where \(\mu_1,\dots, \mu_n\) appear in pairs, one can compute \(p_{1|{}0}\) and \(p_{0|{}1}\) in closed form via a characteristic function approach as in~\cite{Kay_Detection}. Therefore, for the more general non-paired scenarios, it is highly desirable to derive theoretical bounds on~\(p_{1|{}0}\)~and~\(p_{0|{}1}\) whose computations do not require numerical integration.
\end{enumerate}
\textbf{Example}-  Let $\boldsymbol{y}_1,\dots, \boldsymbol{y}_n$ be $n$ observations from a moving average process of order $1$ with dynamics  
\begin{align}
\label{}
   \boldsymbol{y}_i=\boldsymbol{\varepsilon}_i+\theta \boldsymbol{\varepsilon}_{i-1},\,\,\,i=1,\cdots, n,
\end{align}
where $\boldsymbol{\varepsilon}_0, \boldsymbol{\varepsilon}_1, \dots, \boldsymbol{\varepsilon}_n$ are independent $N(0,\sigma^2)$ random variables and $-1<\theta<1$ is the moving average model parameter. Let us consider the hypothesis test 
\begin{align}
\label{}
  H_0: \theta=\theta_0,\quad H_1: \theta=\theta_1,
\end{align} for given constants $-1<\theta_0,\theta_1<1$. 
Equivalently, we have the test of hypothesis in (\ref{H0_H1}) with the covariance matrices 
\begin{eqnarray}
C_0=\sigma^2(1+\theta_0^2)I_n+\sigma^2\theta_0J_n,\quad C_1=\sigma^2(1+\theta_1^2)I_n+\sigma^2\theta_1J_n,
\end{eqnarray}
where $J_n$ is an $n\times n$ matrix whose $(i,j)^{\textrm{th}}$ entry is given by 
\begin{align}
\label{}
   [J_n]_{i,j}=\left\{\begin{array}{cc}
    1  &  |i-j|=1  \\
    0  &   \mathrm{otherwise}
\end{array}\right..
\end{align}
It is well-known that symmetric tridiagonal Toeplitz matrices $aI_n+bJ_n$ for $a,b\in \mathbb{R}$ are simultaneously diagonalizable. Their eigenvalues are given by $a+2b\cos(\frac{i\pi}{n+1})$ for $1\leq i\leq n$. Consequently, the eigenvalues $\mu_i$ of $C_1^{-1}C_0$ can be expressed as   
\begin{eqnarray}
\label{toronto_2}
\mu_i=\frac{1+\theta_0^2+2\theta_0\cos(\frac{i\pi}{n+1})}{1+\theta_1^2+2\theta_1\cos(\frac{i\pi}{n+1})},\,\,\, i=1,\cdots, n.
\end{eqnarray}
Since $\mu_i$ do not depend on $\sigma^2$, it follows that $p_{1|0}$ and $p_{0|1}$ depend only on $\theta_0, \theta_1$ and $n$. The eigenvalues of $A_{1|0}$ are $1-\mu_i=\frac{2(\theta_1-\theta_0)(\frac{\theta_0+\theta_1}{2}+\cos(\frac{i\pi}{n+1}))}{1+\theta_1^2+2\theta_1\cos(\frac{i\pi}{n+1})}$. Given that $-1<\frac{\theta_0+\theta_1}{2}<1$, for sufficiently large $n$, the expression $\frac{\theta_0+\theta_1}{2}+\cos(\frac{i\pi}{n+1})$ changes sign as $i$ ranges from $1$ to $n$. Therefore, the matrix $A_{1|0}$ possesses both  positive and negative eigenvalues. A similar property holds for the matrix $A_{0|1}$. For $\theta_0=-0.5$, $\theta_1=0.7$ and $n=30$, Fig.~\ref{fig_eig1} plots the sorted eigenvalues \(1-\mu_i\) for \(A_{1|{}0}\), while Fig.~\ref{fig_eig2} shows the sorted eigenvalues \(1-\frac{1}{\mu_i}\) for \(A_{0|{}1}\). All positive eigenvalues lie within the interval \((0,1)\). Furthermore, both matrices have negative eigenvalues with larger magnitudes, especially \(A_{1|{}0}\). The modified bounds should therefore outperform the original ones, provided the viability conditions hold.

Recall that the power of the test, denoted by $p_{1|1}(\gamma)$, is defined by $1-p_{0|1}(\gamma)$, where we have made the dependence on the threshold $\gamma$ explicit. Let $p_{1|1,\alpha}$ be the largest achievable power of the test at significance level $\alpha$, and let $\gamma_{\alpha}$ be the optimal threshold that achieves $p_{1|1,\alpha}$, i.e., 
\begin{align}
\label{luis_1}
p_{1|1,\alpha}=\max_{\gamma: \,p_{1|0}(\gamma)\le \alpha} p_{1|1}(\gamma),\quad \gamma_{\alpha}=\arg\max_{\gamma: \,p_{1|0}(\gamma)\le \alpha} p_{1|1}(\gamma).
\end{align}
The only available upper bounds on $p_{1|0}(\gamma)$ are HW (with $\kappa=0.1457$), $m_\infty$, strong $\chi^2$, their modified versions and RK. We denote these seven upper bounds by $p_{1|0, \mathsf{B}}(\gamma)$, where `$\mathsf{B}$' is one of HW, $m_\infty$, strong $\chi^2$, best modified HW, best modified $m_\infty$, best modified strong $\chi^2$ and RK bounds. Each function $p_{1|0, \mathsf{B}}(\gamma)$ has its own domain $\mathscr{D}_{1|0,\mathsf{B}}$. For example, as noted in (\ref{normal_bounds}), $\mathscr{D}_{1|0, \text{HW}}=\mathscr{D}_{1|0, m_\infty}=[-\frac{1}{2}c_{1|0},\infty)$, while $\mathscr{D}_{1|0, \chi^2}=\mathscr{D}_{1|0, \text{RK}}=\mathbb{R}$. For the modified bounds, the domains include the values of $\gamma$ satisfying the respective viability conditions detailed in Propositions~\ref{prop_17}, \ref{prop_9} and \ref{prop_8}. Each of these domains also forms a closed ray extending to $\infty$.  

Since $p_{1|0,\mathsf{B}}(\gamma)$ is a decreasing and right-continuous function with respect to $\gamma$, as explained in Subsection~\ref{2E},  there exists a threshold $\gamma_{\alpha,\mathsf{B}}\in \mathbb{R}$ such that  
\begin{align}
\label{set_interval}
  \{\gamma\in \mathscr{D}_{1|0, \mathsf{B}}:\quad p_{1|0, \mathsf{B}}(\gamma)\leq \alpha\}=[\gamma_{\alpha,\mathsf{B}},\infty).
\end{align}
Each $\gamma_{\alpha,\mathsf{B}}$ serves as an upper bound on~$\gamma_{\alpha}$. This can be justified as follows. By (\ref{set_interval}), $p_{1|0, \mathsf{B}}(\gamma_{\alpha, \mathsf{B}})\leq \alpha$. Since $p_{1|0}(\gamma)\le p_{1|0, \mathsf{B}}(\gamma)$ for all $\gamma\in \mathbb{R}$, we have $p_{1|0}(\gamma_{\alpha,\mathsf{B}})\leq \alpha$. Combining this with the definition of $\gamma_{\alpha}$ in (\ref{luis_1}), it follows that $p_{1|1}(\gamma_{\alpha, \mathsf{B}})\leq p_{1|1}(\gamma_{\alpha})$. This inequality implies that $\gamma_{\alpha}\leq \gamma_{\alpha,\mathsf{B}}$ because $p_{1|1}(\gamma)$ is a non-increasing function of $\gamma$. We denote our best available upper bound on $\gamma_{\alpha}$ by $\hat{\gamma}_{\alpha}$ defined~as 
\begin{align}
\label{hat_gamma}
\hat{\gamma}_{\alpha}=\min_{\mathsf{B}}\gamma_{\alpha,\mathsf{B}}.
\end{align}   
To compute each $\gamma_{\alpha,\mathsf{B}}$, one needs to solve the inequality $p_{1|0,\mathsf{B}}(\gamma)\leq \alpha$ for $\gamma$. This can be done efficiently as follows. Let $\mathscr{G}_1\subseteq \mathscr{G}_2\subseteq \cdots$ be a nested sequence of grids of real numbers such that $\bigcup_{l=1}^\infty\mathscr{G}_l$ is dense in~$\mathbb{R}$. For example, let $\mathscr{G}_l=\{m2^{-l}: m\in \mathbb{Z}\}$ for $l\ge1$. For each $l$, find the smallest $\gamma\in \mathscr{G}_l$ (denoted by $\gamma^{(l)}_{\alpha,\mathsf{B}}$) that satisfies $p_{1|0,\mathsf{B}}(\gamma)\leq \alpha$. The resulting sequence $(\gamma^{(l)}_{\alpha,\mathsf{B}})_{l\ge1}$ is decreasing and converges to $\gamma_{\alpha,\mathsf{B}}$. In fact, one can use each $\gamma^{(l)}_{\alpha,\mathsf{B}}$ in place of $\gamma_{\alpha,\mathsf{B}}$ in~(\ref{hat_gamma}), which results in only a slightly looser upper bound on $\gamma_\alpha$.  Alternatively, one can compute $\gamma_{\alpha,\text{HW}}$ in closed form as 
\begin{align}
\label{gamma_HW}
\gamma_{\alpha,\text{HW}}=\frac{1}{2}\max\Big\{ \Big(-\frac{\ln\alpha}{\kappa}\Big)^{1/2}\|A_{1|0}\|_2, -\frac{\ln\alpha}{\kappa}\|A_{1|0}\|\Big\}-\frac{1}{2}c_{1|0},\quad \kappa=0.1457.
\end{align}
The expressions for $\gamma_{\alpha,m_\infty}$ and $\gamma_{\alpha,\chi^2}$  involve the Lambert W function and the inverse CDF of the $\chi^2$ distribution, respectively, and are given by  
\begin{align}
\label{gamma_Minf}
 \gamma_{\alpha, m_\infty}=-\frac{n}{2}\|A_{1|0}\|\big(1+W_{-1}(-\alpha^{\frac{2}{n}}/e)\big)-\frac{1}{2}c_{1|0}
\end{align}
and
\begin{align}
\label{gamma_chi2}
\gamma_{\alpha, \chi^2}=\frac{1}{2}\lambda_{\max}(A_{1|0})F_{n}^{-1}(1-\alpha)-\frac{1}{2}\mathrm{tr}(A_{1|0})-\frac{1}{2}c_{1|0},
\end{align}
where $c_{1|0}$ is defined in (\ref{c_1_0}) and $F_n^{-1}(\cdot)$ is the inverse CDF function for a $\chi^2_n$ random variable. Both (\ref{gamma_Minf}) and (\ref{gamma_chi2}) assume that $r(A_{1|0})=n$, which holds true if $\cos(\frac{i\pi}{n+1})\neq -\frac{\theta_0+\theta_1}{2}$ for all $1\le i\leq n$. Otherwise, $n$ is replaced by $n-1$, as there is at most one choice of $1\le i\le n$ for which $\cos(\frac{i\pi}{n+1})=-\frac{\theta_0+\theta_1}{2}$ holds. 


As in the case of $p_{1|0}(\gamma)$, the only available upper bounds on $p_{0|1}(\gamma)$, and hence the only available lower bounds on the power of the test $p_{1|1}(\gamma)$, are HW ($\kappa=0.1457$), $m_\infty$, strong $\chi^2$, their modified versions and RK. We denote these seven lower bounds by $p_{1|1, \mathsf{B}}(\gamma)$, where `$\mathsf{B}$' is one of HW, $m_\infty$, strong $\chi^2$, best modified HW, best modified $m_\infty$, best modified strong $\chi^2$ and RK bounds. Substituting $\hat{\gamma}_{\alpha}$ (given in (\ref{hat_gamma})) for~$\gamma$ yields seven lower bounds $p_{1|1, \mathsf{B}}(\hat{\gamma}_{\alpha})$ on the largest achievable power $p_{1|1,\alpha}=p_{1|1}(\gamma_\alpha)$ defined in~(\ref{luis_1}). This substitution is valid only when $\hat{\gamma}_{\alpha}$ belongs to the domain of $p_{1|1,\mathsf{B}}(\gamma)=1-p_{0|1,\mathsf{B}}(\gamma)$, denoted by~$\mathscr{D}_{0|1,\mathsf{B}}$. For example, as noted in (\ref{normal_bounds_1}), $\mathscr{D}_{0|1,\text{HW}}=\mathscr{D}_{0|1, m_\infty}=(-\infty, \frac{1}{2}c_{0|1}]$. Thus, $p_{1|1,\text{HW}}(\hat{\gamma}_\alpha)$ and $p_{1|1,m_\infty}(\hat{\gamma}_\alpha)$ are valid if and only if $\hat{\gamma}_{\alpha}\leq \frac{1}{2}c_{0|1}$.  
 
  \begin{figure*}[t]
\centering
\subfigure[]{
\includegraphics[scale=0.35]{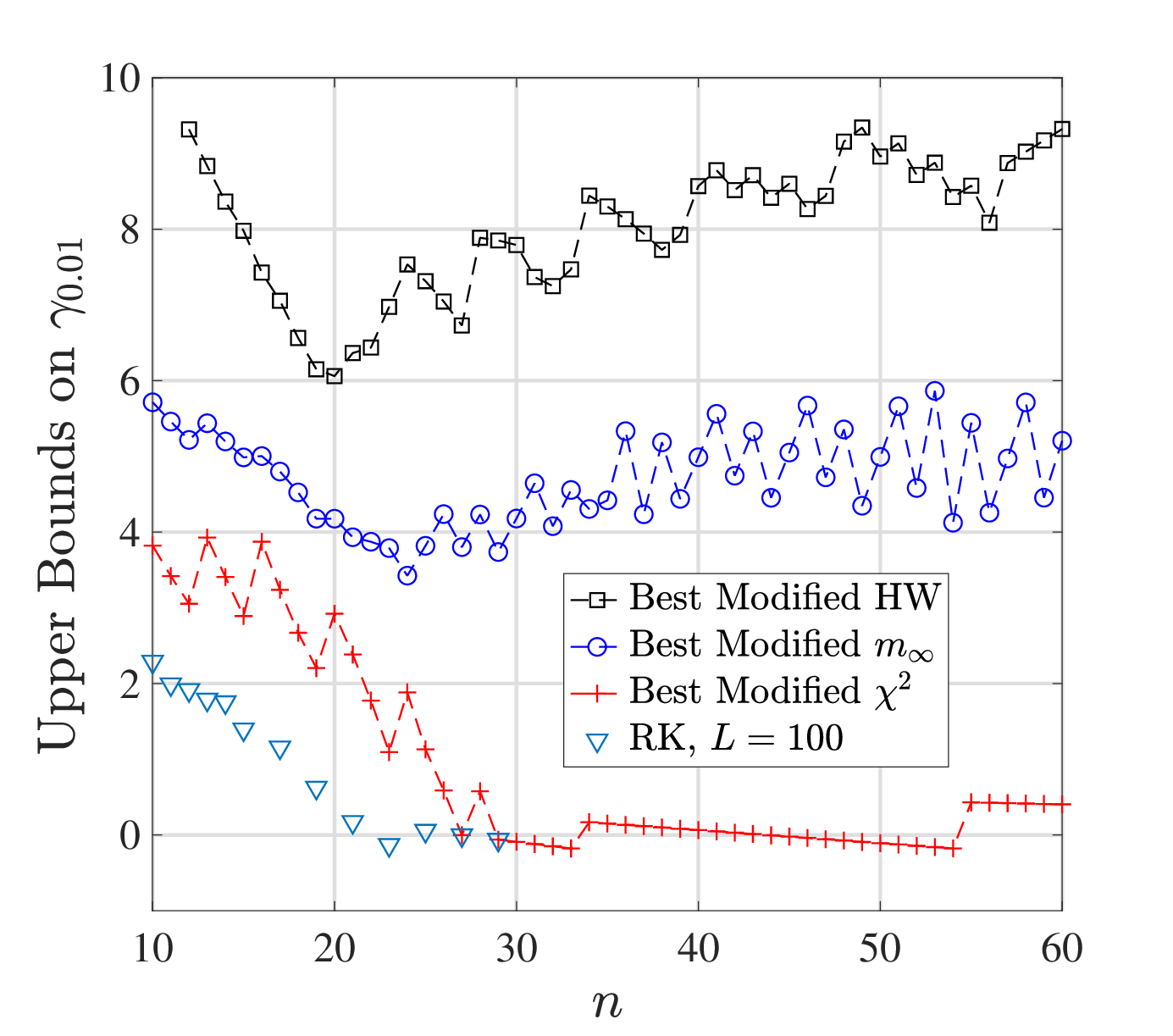}
\label{fig_gamma_n}
}
\subfigure[]{
\includegraphics[scale=0.35]{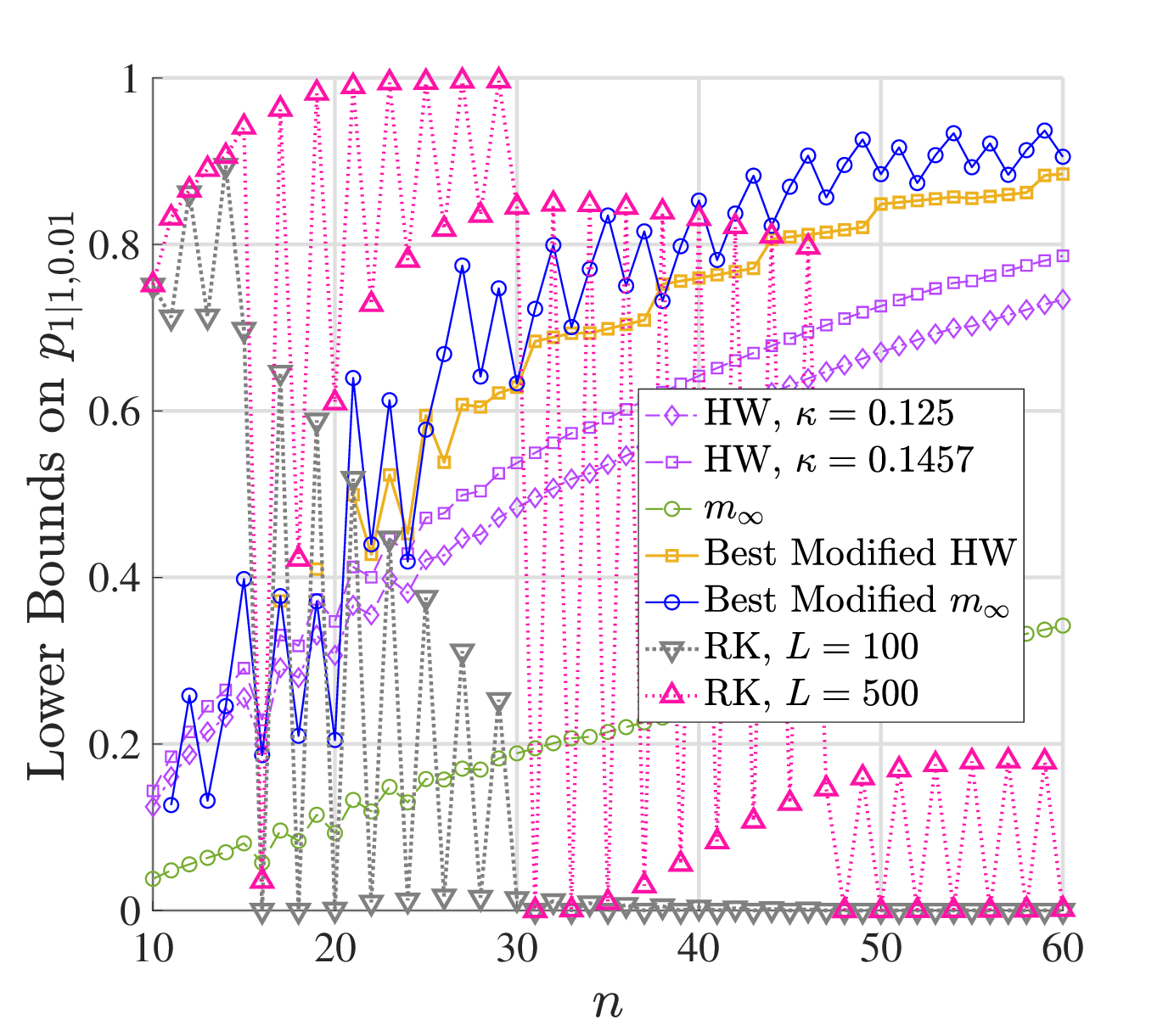}
\label{fig_p11_n}
}
\label{fig:subfigureExample}
\caption[Optional caption for list of figures]{Let $\alpha=0.01$, $\theta_0=-0.5$ and $\theta_1=0.7$. Panel~(a) presents plots of various upper bounds on $\gamma_{0.01}$ with respect to $10\le n\le60$, where we have only considered the RK bound of truncation level $L=100$. Panel~(b) plots the lower bounds on $p_{1|1 ,0.01}$ against $n$, where we have only considered the RK bounds of truncation levels $L=100$ and $500$.   }
\label{Results_TDL3}
\end{figure*}

 \textcolor{Black}{ Let $\alpha=0.01$, $\theta_0=-0.5$ and $\theta_1=0.7$. Fig.~\ref{fig_gamma_n} displays various upper bounds on $\gamma_{0.01}$ against $10\le n\le 60$, considering only the RK bound with a truncation level of $L=100$. Depending on the value of $n$, the tightest upper bound is either the best modified $\chi^2$ bound or the RK bound. Because the values for $\gamma_{0.01, \text{HW}}$, $\gamma_{0.01,m_\infty}$ and $\gamma_{0.01,\chi^2}$ are considerably larger, they are omitted from this plot. Next, Fig.~\ref{fig_p11_n} plots the lower bounds on $p_{1|1 ,0.01}$ against $n$,  considering only the RK bounds of truncation levels $L=100$ and $500$.  The strong $\chi^2$ bound is excluded due to its poor performance, as its values remain very close to zero. Moreover, none of the modified strong $\chi^2$ bounds are viable. For the RK bounds, we observe a downward oscillating trend. This downward trajectory is expected, as a larger truncation level $L$ is typically required at higher values of $n$ to maintain the same approximation gap. Notably, at certain values of $n$ such as 16 and 31, there is negligible difference between the RK bounds for $L=100$ and $L=500$. This behaviour stems from an extremely slow convergence rate of the series $\sum_{l=0}^\infty c_l^{\text{RK}^\pm}$ at those specific points. For every $39\le n\le60$, the best modified $m_\infty$ bound remains the tightest among all considered bounds.~$\diamond$     }

\begin{figure*}[t]
\centering
\subfigure[]{
\includegraphics[scale=0.3]{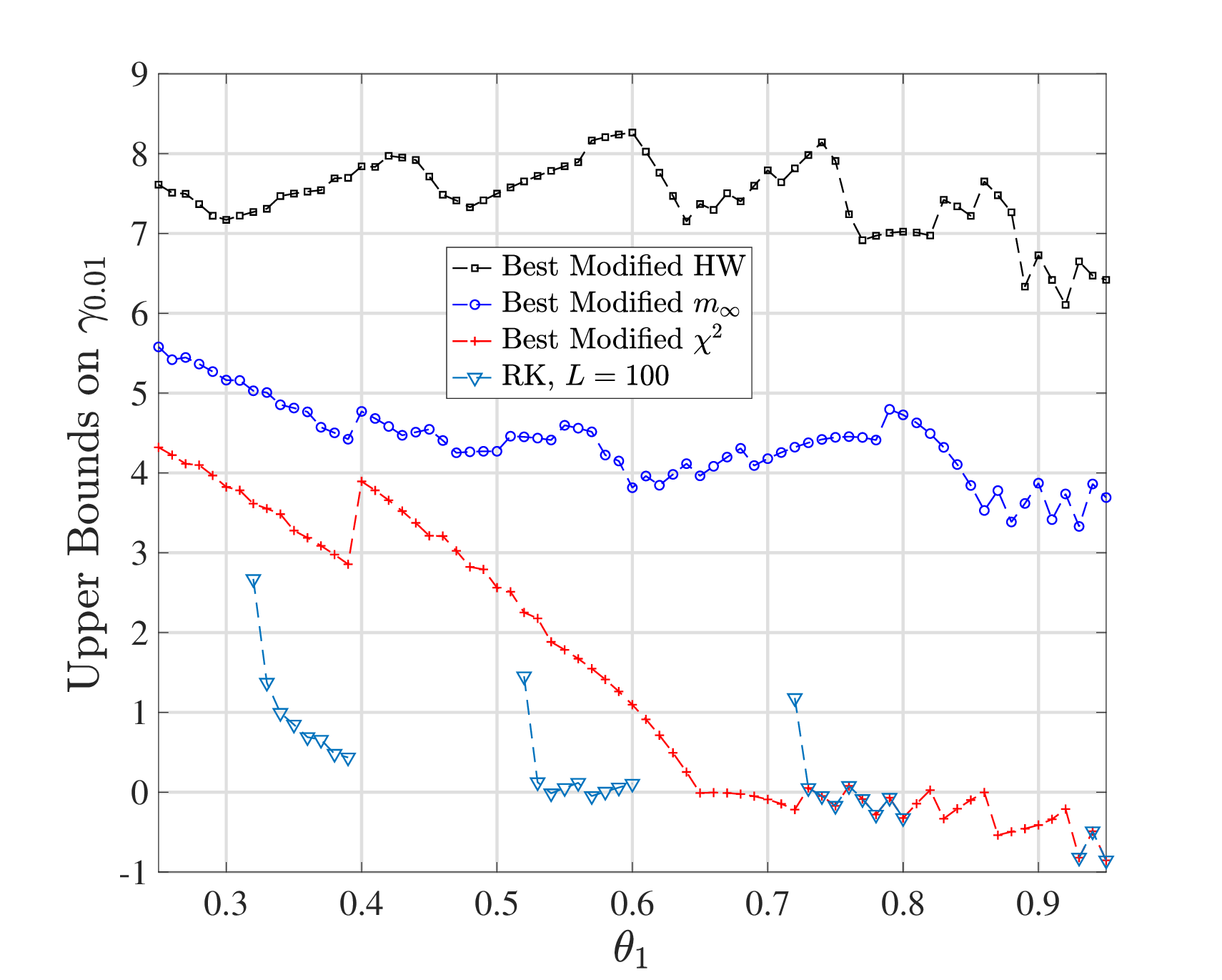}
\label{fig_gamma_theta}
}
\subfigure[]{
\includegraphics[scale=0.3]{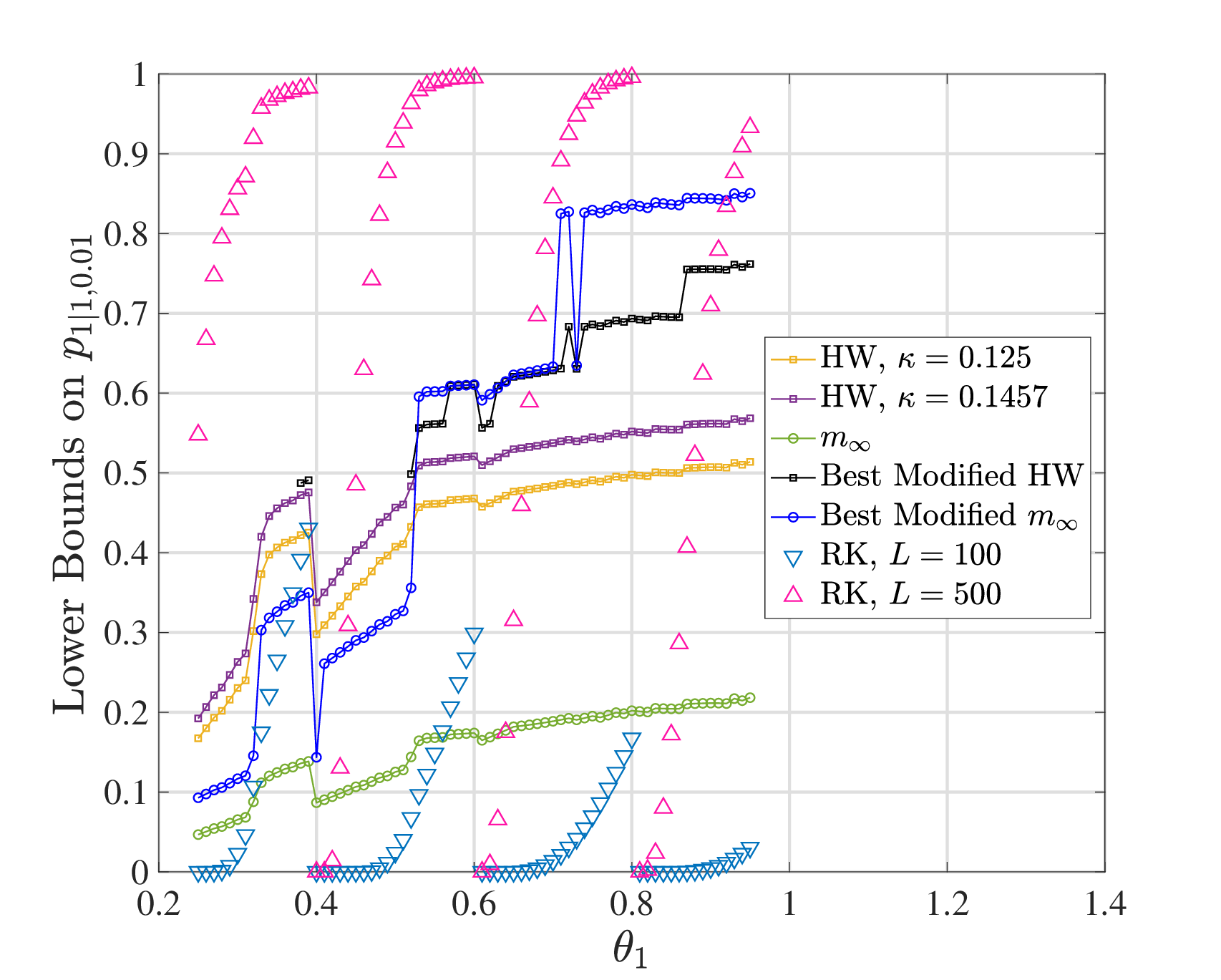}
\label{fig_p11_theta}
}
\label{fig:subfigureExample}
\caption[Optional caption for list of figures]{Let $\alpha=0.01$, $n=30$ and $\theta_0=-0.5$. Panel~(a) shows various upper bounds on $\gamma_{0.01}$ against $0.25\le \theta_1\le 0.95$ in steps of 0.01. 
 }
\label{Results_TDL3}
\end{figure*}

  Let $\alpha=0.01$, $n=30$ and $\theta_0=-0.5$. Fig.~\ref{fig_gamma_theta} shows various upper bounds on $\gamma_{0.01}$ against $0.25\le \theta_1\leq 0.95$ in steps of 0.01,  where we have only considered the RK bound of truncation level $L=100$. Once again, the best modified $\chi^2$ bound and the RK bound yield the tightest upper bound  depending on the value of $\theta_1$. Fig.~\ref{fig_p11_theta} displays the lower bounds on $p_{1|1 ,0.01}$ as a function of $\theta_1$,  where we have only considered the RK bounds of truncation levels $L=100, 500$.~$\diamond$

 
\subsection{Wireless Communications}
\label{3_C}
Let us consider MISO wireless communications in Rayleigh flat-fading with $M$ transmitter antennas and one receiver antenna. Restriction to one receiver antenna is only for simplicity of presentation. The results can be easily extended to any number of receiver antennas.  The channel coefficients $\boldsymbol{h}_{i}\sim CN(0,1)$ between the $i^{\textrm{th}}$ transmitter antenna and the receiver antenna for $i=1,\dots, M$ are independent standard circularly symmetric complex Gaussian random variables that stay constant over a coherence period of $T$ symbol intervals. These coefficients then change independently and stay constant over the next coherence period of the same length.  The channel equation over one such period is given by 
\begin{align}
\label{mimo1}
\underline{\boldsymbol{y}}=\zeta\boldsymbol{S}\underline{\boldsymbol{h}}+\underline{\boldsymbol{z}},\quad \zeta=\Big(\frac{T\mathsf{SNR}}{M}\Big)^{1/2},
\end{align}
where $\boldsymbol{S}$ is the $T\times M$ transmitted matrix, $\underline{\boldsymbol{h}}=[\boldsymbol{h}_1\,\,\cdots\quad \boldsymbol{h}_M]^{\mathsf{T}}\sim CN(\underline{0},I_M)$ is the channel vector which is \textit{unknown} to both the transmitter and the receiver, $\underline{\boldsymbol{z}}\sim CN(\underline{0},I_T)$ is the additive noise vector at the receiver, $\underline{\boldsymbol{y}}$ is the received vector of length~$T$ 
and $\mathsf{SNR}$ is the expected signal-to-noise ratio at the receiver antenna. Motivated by the optimal capacity-achieving  distribution for the random matrix $\boldsymbol{S}$, the authors of~\cite{Hochwald1} introduced unitary space-time modulation. In this approach, information bits at the transmitter are mapped to $T\times M$ unitary matrices $\Phi$ that satisfy $\Phi^{\mathsf{H}}\Phi=I_M$, where $A^{\mathsf{H}}$ denotes the conjugate transpose of matrix $A$. For simplicity, we only address binary modulation, where the  bits $0$ and $1$  are mapped to~$T\times M$ unitary matrices $\Phi_0$~and~$\Phi_1$, respectively. The optimal maximum-likelihood decoder takes on the simple form 
\begin{align}
\hat{\Phi}=\arg\max_{\Phi\in \{\Phi_0, \Phi_1\}}\underline{y}^{\mathsf{H}}\Phi\Phi^{\mathsf{H}}\underline{y}.
\end{align}     
The probability of decoding error can be written~as 
\begin{align}
p_e=\frac{1}{2}(\Pr(\mathrm{error}|\Phi_0)+\Pr(\mathrm{error}|\Phi_1)),
\end{align}  
where $\Pr(\mathrm{error}|\Phi_i)$ is the probability of decoding error assuming $\Phi_i$ is transmitted for $i=0,1$. It turns out that $\Pr(\mathrm{error}|\Phi_0)=\Pr(\mathrm{error}|\Phi_1)$ and hence,
\begin{align}
\label{soodokoo1}
p_e=\Pr(\mathrm{error}|\Phi_0)=\Pr\big(\underline{\boldsymbol{y}}^{\mathsf{H}}(\Phi_1\Phi_1^{\mathsf{H}}-\Phi_0\Phi_0^{\mathsf{H}})\underline{\boldsymbol{y}}>0| \Phi_0\big).
\end{align}
 A key observation made in \cite{Hochwald1} is that $p_e$ depends on $\Phi_0$ and $\Phi_1$ only through the singular values $0\leq d_1\leq d_2\leq \cdots\leq d_M$ of the $M\times M$ matrix $\Phi_1^{\mathsf{H}}\Phi_0$. 
Define the diagonal matrix 
\begin{align}
D=\mathrm{diag}(d_1,\cdots, d_M).
\end{align}
If $0$ is the transmitted bit, it follows that  
\begin{align}
\label{soodokoo2}
\underline{\boldsymbol{y}}^{\mathsf{H}}(\Phi_1\Phi_1^{\mathsf{H}}-\Phi_0\Phi_0^{\mathsf{H}})\underline{\boldsymbol{y}}=\underline{\boldsymbol{w}}^{\mathsf{H}}J\underline{\boldsymbol{w}},
\end{align}
where $\underline{\boldsymbol{w}}\sim CN(\underline{0},C)$ is a vector of length $2M$. Here, the covariance matrix $C$ is given by 
\begin{align}
\label{c_mat}
C=\mathbb{E}[\underline{\boldsymbol{w}}\,\underline{\boldsymbol{w}}^{\mathsf{H}} ]=\begin{bmatrix}
    (\zeta^2+1)I_M  & (\zeta^2+1)D   \\
     (\zeta^2+1)D &  \zeta^2D^2+I_M
\end{bmatrix}
\end{align}
and the $2M\times 2M$ matrix $J$ is defined by 
\begin{align}
\label{j_mat}
J=\begin{bmatrix}
    -I_M &  O_M  \\
    O_M  &  I_M
\end{bmatrix},
\end{align}
where $O_M$ denotes the $M\times M$ zero matrix. By (\ref{soodokoo1}) and (\ref{soodokoo2}), 
\begin{align}
\label{astronomical1}
p_e=\Pr(\underline{\boldsymbol{w}}^{\mathsf{H}}J\underline{\boldsymbol{w}}>0).
\end{align}
In~\cite{Hochwald1},  $p_e$ is computed using a characteristic function approach, which leads to the simple upper bound 
\begin{align}
\label{H_M}
p_e\leq \frac{1}{2}\prod_{i=1}^M \frac{1}{1+\frac{\zeta^4(1-d_i^2)}{4(1+\zeta^2)}}.
\end{align}
We will refer to this bound as Hochwald-Marzetta bound or the~HM~bound for short. 
Instead, we choose to work directly with the quadratic form in (\ref{astronomical1}) to write $p_e$ in a form that is suitable for applying the bounds in Section~\ref{sec1}. Representing $\underline{\boldsymbol{w}}$ as 
\begin{align}
\underline{\boldsymbol{w}}=C^{1/2}\underline{\tilde{\boldsymbol{x}}},\quad \underline{\tilde{\boldsymbol{x}}}\sim CN(\underline{0},I_{2M}),
\end{align}
we can write
\begin{align}
\label{pnm_1}
p_e=\Pr(\underline{\tilde{\boldsymbol{x}}}^{\mathsf{H}}B\underline{\tilde{\boldsymbol{x}}}>0),
\end{align}
where 
\begin{align}
B=C^{1/2}JC^{1/2}.
\end{align}
The matrix $B$ is real. If $\underline{\tilde{\boldsymbol{x}}}_r$ and $\underline{\tilde{\boldsymbol{x}}}_i$ are the real and imaginary parts of $\underline{\tilde{\boldsymbol{x}}}$, then 
\begin{align}
\label{pnm_2}
\underline{\tilde{\boldsymbol{x}}}^{\mathsf{H}}B\underline{\tilde{\boldsymbol{x}}}=\underline{\tilde{\boldsymbol{x}}}_r^{\mathsf{T}}B\underline{\tilde{\boldsymbol{x}}}_r+\underline{\tilde{\boldsymbol{x}}}_i^{\mathsf{T}}B\underline{\tilde{\boldsymbol{x}}}_i=\frac{1}{2}\underline{\boldsymbol{x}}^{\mathsf{T}}A\underline{\boldsymbol{x}},
\end{align}
where $A$ and $\underline{\boldsymbol{x}}$ are given by\footnote{Since $\underline{\tilde{\boldsymbol{x}}}_r,\underline{\tilde{\boldsymbol{x}}}_i\sim N(\underline{0}, \frac{1}{2}I_{2M})$ are independent, it follows that $\underline{\boldsymbol{x}}\sim N(\underline{0},I_{4M})$. } 
\begin{align}
\label{A_B}
A=\begin{bmatrix}
    B  & O_{2M}   \\
    O_{2M}  & B 
\end{bmatrix},\quad \underline{\boldsymbol{x}}=\sqrt{2}[\underline{\tilde{\boldsymbol{x}}}_r^{\mathsf{T}}\quad \underline{\tilde{\boldsymbol{x}}}_i^{\mathsf{T}}]^{\mathsf{T}}\sim N(\underline{0},I_{4M}).
\end{align}
By (\ref{pnm_1}) and (\ref{pnm_2}),  
\begin{align}
\label{WCom}
p_e=\Pr(\underline{\boldsymbol{x}}^{\mathsf{T}}A\underline{\boldsymbol{x}}>0).
\end{align}
Subtracting $\mathrm{tr}(A)$ from both sides of the inequality in (\ref{WCom}), 
\begin{align}
p_e=\Pr\big(\underline{\boldsymbol{x}}^{\mathsf{T}}A\underline{\boldsymbol{x}}-\mathrm{tr}(A)>-\mathrm{tr}(A)\big).
\end{align}
Now, $p_e$ is in a form that is suitable for applying the bounds in Section~\ref{sec1}. The deviation parameter $t$ is given by 
\begin{align}
t=-\mathrm{tr}(A).
\end{align}
 It is easy to see that $t>0$. In fact, 
\begin{align}
\label{pnm_3}
\mathrm{tr}(A)=2\mathrm{tr}(B)=2\mathrm{tr}(C^{1/2}JC^{1/2})=2\mathrm{tr}(C^{1/2}C^{1/2}J)=2\mathrm{tr}(CJ).
\end{align} 
By (\ref{c_mat}) and (\ref{j_mat}), 
\begin{align}
CJ=\begin{bmatrix}
    -(\zeta^2+1)I_M  & (\zeta^2+1)D   \\
     -(\zeta^2+1)D &  \zeta^2D^2+I_M
\end{bmatrix}.
\end{align}
Therefore, 
\begin{align}
\label{pnm_4}
\mathrm{tr}(CJ)=\mathrm{tr}\big(\zeta^2D^2+I_M-(\zeta^2+1)I_M\big)=\zeta^2\mathrm{tr}(D^2-I_M)=\zeta^2\sum_{i=1}^M(d_i^2-1).
\end{align} 
It is shown in \cite{Hochwald1} that $0\leq d_i\leq 1$ for all $i=1,\dots, M$. If $d_1=\cdots=d_M=1$, then one can easily verify that $A=0$. This leads to $\underline{\boldsymbol{y}}^{\mathsf{H}}(\Phi_1\Phi_1^{\mathsf{H}}-\Phi_0\Phi_0^{\mathsf{H}})\underline{\boldsymbol{y}}=\frac{1}{2}\underline{\boldsymbol{x}}^{\mathsf{T}}A\underline{\boldsymbol{x}}=0$ and hence, the receiver can not decide which one of $0$ or $1$ is transmitted. To avoid this tie, there must exist at least one $1\leq i\leq M$ such that $d_i< 1$. Using this fact together with (\ref{pnm_3}) and (\ref{pnm_4}), we conclude that $\mathrm{tr}(A)<0$ and hence, $t=-\mathrm{tr}(A)>0$. The matrix $A$ is symmetric, however, it is not positive-semidefinite as we will verify momentarily. Thus, the only available upper bounds on $p_e$ are HW, $m_\infty$, strong $\chi^2$ and their modified versions.  To compute these bounds, one needs to find the eigenvalues of $A$. By (\ref{A_B}), eigenvalues of $A$ are those of~$B$ with double multiplicity. Eigenvalues of $B=C^{1/2}JC^{1/2}$ are exactly those of $C^{1/2}C^{1/2}J=CJ$.  It is a straightforward observation that for a block matrix $\small \begin{bmatrix}
    M_{1,1}  &  M_{1,2}  \\
    M_{2,1}  &  M_{2,2}
\end{bmatrix}$, where $M_{i,j}$ are diagonal matrices of the same size, the determinant is given by $\det (M_{1,1}M_{2,2}-M_{1,2}M_{2,1})$. Using this formula and applying simple algebra, the characteristic function of~$B=C^{1/2}JC^{1/2}$ is given by
\begin{align}
\det(\psi I_{2M}-B)=\det(\psi I_{2M}-CJ)=\prod_{i=1}^M\big(\psi^2+\zeta^2(1-d_i^2)\psi-(\zeta^2+1)(1-d_i^2)\big)
\end{align}
and the eigenvalues of $B$ are given by 
\begin{align}
\label{poolgh}
\psi_{i}^{\pm}=\left\{\begin{array}{cc}
    0  &  d_i=1  \\
    \frac{1}{2}(1-d_i^2)\zeta^2\big(-1\pm (1+\frac{4(\zeta^2+1)}{(1-d_i^2)\zeta^4})^{1/2}\big)  &   0\leq d_i<1
\end{array}\right.,\quad i=1,\dots, M.
\end{align}
We know that there exists $1\leq i\leq M$ such that $d_i<1$. Then $\psi_{i}^{-}<0$ and $\psi_{i}^{+}>0$. This verifies that $A$ is not positive-semidefinite. It is also easy to see that if $d_i<1$, then $\psi_i^-$ scales linearly with $\zeta^2$ (and hence with $\mathsf{SNR}$), while $0<\psi_i^+<1$ regardless of $\zeta>0$. In other words, 
\begin{align}
\label{pio1}
|\psi_i^-|=\Theta(\mathsf{SNR})\,\,\, \text{as $\mathsf{SNR}\to\infty$} 
\end{align}
and 
\begin{align}
\label{pio2}
\psi_i^+=\mathcal{O}(1)\,\,\, \text{as $\mathsf{SNR}\to\infty$}.
\end{align}
 Without loss of generality, we assume $0\leq d_1\leq d_2\leq \cdots\le d_{M}<1$. Otherwise, one can replace $M$ in the following by $M'=|\{1\le i\leq M: d_i<1\}|$. We have  $\mathrm{rank}(A)=4M$.  The HW~bound and the $m_\infty$~bound are given~by  $p_e\leq \exp(-\kappa\min\{\frac{(\mathrm{tr}(A))^2}{\|A\|_2^2},-\frac{\mathrm{tr}(A)}{\|A\|}\})$
and $p_e\leq \big(1-\frac{\mathrm{tr}(A)}{4M\|A\|}\big)^{2M}e^{\frac{\mathrm{tr}(A)}{2\|A\|}}$, respectively. Since $\|A\|$, $\|A\|_2$ and $\mathrm{tr}(A)$ all scale linearly with $\mathsf{SNR}$, both bounds saturate as $\mathsf{SNR}$ increases. The strong $\chi^2$~bound in (\ref{chi2_bound}) is trivial as it evaluates to $1$ due to $\lambda^*_{\max}(A)>0$ and $t+\mathrm{tr}(A)=0$. The modified HW and modified $m_\infty$ bounds also saturate as $\mathsf{SNR}$ increases. Let us explain this behaviour for the $m_\infty$~bound of order $k$ in Proposition~\ref{prop_17}. The order $k$ is even due to the requirement that $\mathrm{rank}(A)-k=4M-k$ must be even. The matrix $A$ has a total of $2M$ negative eigenvalues $\psi^-_1, \psi_1^-, \psi_2^-, \psi_2^-,\dots, \psi_{M}^-, \psi_{M}^-$ , as given in~(\ref{poolgh}). Therefore, the possibilities for the order $k$ are $2,4,\dots,2M$. Recall the parameters $a_k, b_k$ in (\ref{a_k_b_k}). If $k=2M$, then $b_{2M}=2\sum_{i=1}^{M}\psi_{i}^-$ and hence,
\begin{align}
t+b_{2M}=-\mathrm{tr}(A)+2\sum_{i=1}^{M}\psi_{i}^-=-2\sum_{i=1}^{M}\psi_i^+-2\sum_{i=1}^{M}\psi_i^-+2\sum_{i=1}^{M}\psi_{i}^-=-2\sum_{i=1}^{M}\psi_i^+<0.
\end{align}     
This violates the condition specified in (\ref{poloi_2}). Therefore, $k\leq 2M-2$. Since $0\leq d_1\leq d_2\leq \cdots\le d_{M}$, then $\psi^-_1\leq \psi_2^-\leq \cdots \leq \psi_{M}^-<0$ and we get 
\begin{align}
\label{qcv_1}
t+b_{k} &=-\mathrm{tr}(A)+2\sum_{i=1}^{\frac{k}{2}}\psi_i^-\notag\\&=-2\sum_{i=1}^{M}\psi_i^+-2\sum_{i=1}^{M}\psi_i^-+2\sum_{i=1}^{\frac{k}{2}}\psi_i^-\notag\\&=-2\sum_{i=1}^{M}\psi_i^+-2\sum_{i=\frac{k}{2}+1}^{M}\psi_i^-\notag\\
&=\Theta(\mathsf{SNR})\,\,\, \text{as $\mathsf{SNR}\to\infty$},
\end{align}
where the first step follows from the fact that $b_k$ is the sum of the $k$ largest negative eigenvalues of $A$ and the last step is due to (\ref{pio1}), (\ref{pio2}) and the fact that the sum $\sum_{i=\frac{k}{2}+1}^{M}\psi_i^-$ is not void due to $k\leq 2M-2$. Similarly, It is evident that 
\begin{align}
\label{qcv_2}
a_k=\max\{\psi_1^+,\cdots, \psi_M^+, -\psi^-_{\frac{k}{2}+1},\cdots, -\psi^-_{M} \}=\Theta(\mathsf{SNR})\,\,\, \text{as $\mathsf{SNR}\to\infty$}.
\end{align}
By (\ref{qcv_1}) and (\ref{qcv_2}), both $t+b_k$ and $a_k$ scale linearly with $\mathsf{SNR}$. Hence, the ratio $\frac{t+b_k}{a_k}$ (and consequently, the $m_\infty$~bound of order $k\leq 2M-2$) saturates as $\mathsf{SNR}$ increases.  

We observe that the~HM~bound in (\ref{H_M}) scales like $\frac{1}{\mathsf{SNR}^M}$ when $\mathsf{SNR}$ is large. As such, none of HW, $m_\infty$ and strong~$\chi^2$ bounds and the modified versions of HW and $m_\infty$ bounds will be sharper than the~HM~bound in the large signal-to-noise ratio regime. Numerical results indicate that these bounds are also not tighter than the~HM~bound for small or moderate values of $\mathsf{SNR}$. The only remaining bounds to consider are the modified strong $\chi^2$ bounds in (\ref{modified_chi2_final}). Since $t=-\mathrm{tr}(A)$, the condition in (\ref{poloi_4}) of Proposition~\ref{prop_9} always holds. Let us look at the strong $\chi^2$~bound of maximum possible order $k=2M$. Defining 
\begin{align}
\psi_{\max}^+=\max_{1\leq i\leq M}\psi_i^+,\quad \tilde{d}_{m}=2\sum_{i=1}^M\left(\frac{\frac{-2\psi_i^-}{\psi^+_{\max}}}{\frac{-\psi_i^-}{\psi^+_{\max}}+1}\right)^m,\,\,\,m\ge1,
\end{align}
the bound is given by   
\begin{align}
\label{full_chi2}
p_e\leq \sum_{m=0}^{M-1}\frac{\tilde{c}_{m}}{2^m},
\end{align}
where the coefficients $\tilde{c}_{m}$ are computed via the recursion 
\begin{align}
\label{}
 \tilde{c}_0=\frac{1}{\prod_{i=1}^{M}(\frac{-\psi_i^-}{\psi^+_{\max}}+1)},\quad \tilde{c}_{m}=\frac{1}{2m}\sum_{j=0}^{m-1}  \tilde{d}_{m-j}\tilde{c}_{j},\,\,\,m\geq1.
\end{align}
As observed earlier, each $-\psi^{-}_i$ scales like $\mathsf{SNR}$, whereas $\psi_{\max}^+<1$ does not scale with $\mathsf{SNR}$. Thus, $\tilde{c}_0$ scales like $\frac{1}{\mathsf{SNR}^M}$ and $\tilde{d}_1,\dots, \tilde{d}_{M-1}$ do not scale with $\mathsf{SNR}$. Consequently, all coefficients $\tilde{c}_1,\dots, \tilde{c}_M$ and therefore the entire bound scale like~$\frac{1}{\mathsf{SNR}^M}$, which matches the scaling of the~HM~bound. Additionally,~(\ref{full_chi2}) can offer substantial improvements over the~HM~bound for a wide range of values of $\mathsf{SNR}$, as demonstrated by the following example.
\begin{figure*}[t]
\centering
\subfigure[]{
\includegraphics[scale=0.37]{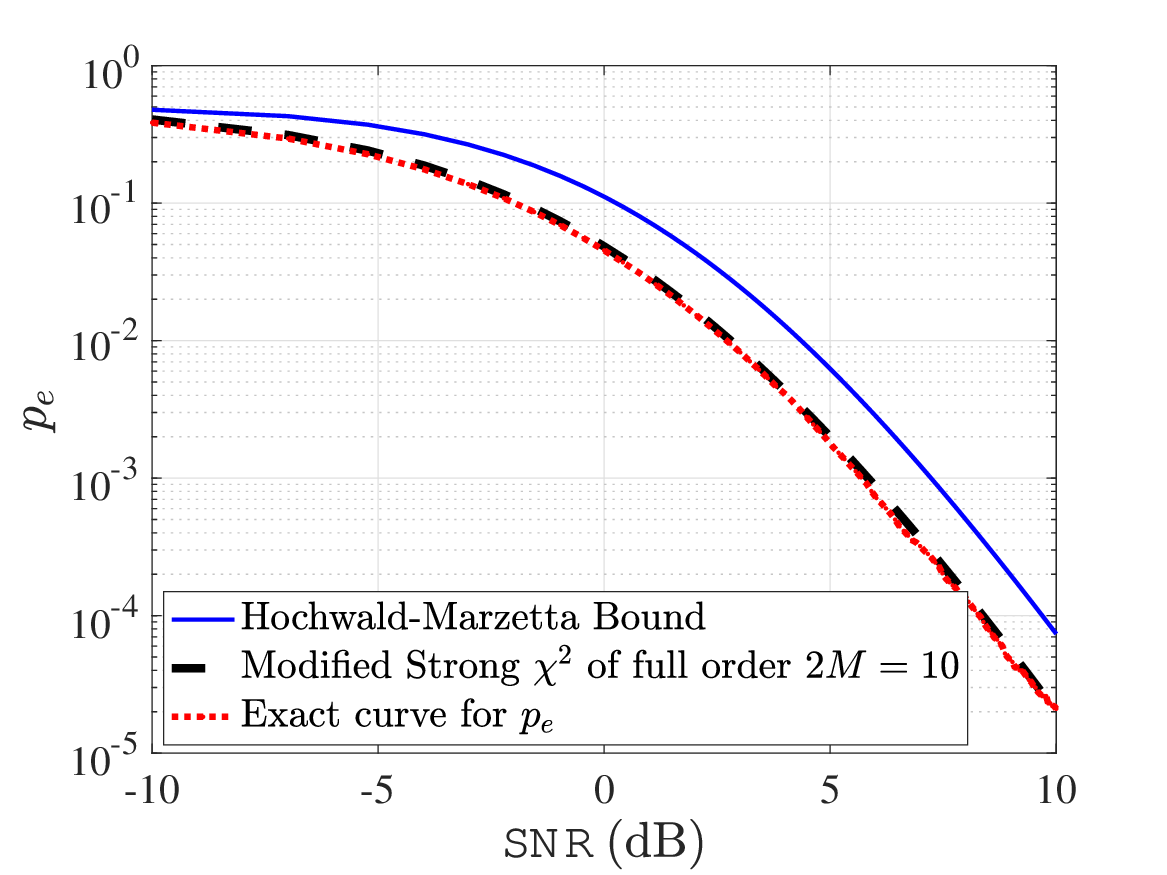}
\label{}
}
\subfigure[]{
\includegraphics[scale=0.37]{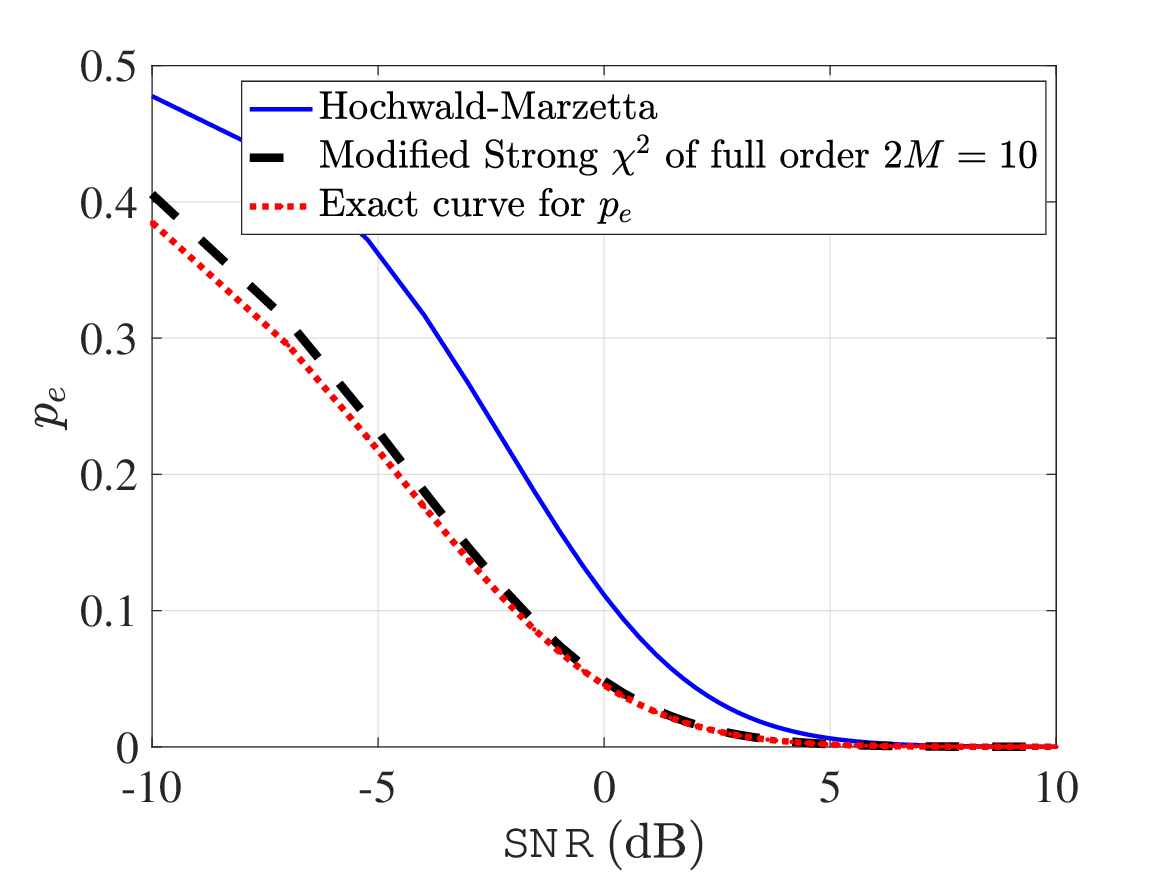}
\label{}
}
\label{fig:subfigureExample}
\caption[Optional caption for list of figures]{Let $M=5, T=11$. Consider the unitary space-time code construction based on the DFT matrix in \cite{Hochwald2}. Panel~(a) presents plots of the~HM~bound in (\ref{H_M}) and the strong $\chi^2$~bound of full order $k=2M=10$ in (\ref{full_chi2}) with respect to $\mathsf{SNR}$ between $-10\,\mathrm{dB}$ to $10\,\mathrm{dB}$. For comparison, the exact curve for $p_e$, evaluated via Monte~Carlo simulation with $10^6$ trials, is also included. Panel~(b) displays the same curves on a linear vertical axis rather than a logarithmic one.}
\label{fig22}
\end{figure*}

\textbf{Example}- Let $M=5$ and $T=11$. We utilize the unitary space-time code construction based on the DFT matrix as proposed in~\cite{Hochwald2}. Let $\Phi_0$ be the $11\times 5$ unitary matrix whose columns are the first $5$ columns of the DFT matrix of size $11\times 11$. Let $\Upsilon$ be an $11\times 11$ diagonal matrix whose diagonal entries are $-1$ or $1$ and set $\Phi_1=\Upsilon\Phi_0$. Among all $2^{11}$ possibilities for $\Upsilon$, we select the one for which $\max_{1\leq i\leq M}d_i$ is the smallest. This results in the singular values $d_1=0.0001, d_2=0.0018, d_3= 0.0207, d_4=0.1433$ and $d_5=0.5789$. Fig.~\ref{fig22} presents plots of the~HM~bound in (\ref{H_M}) and the strong $\chi^2$~bound of full order $2M=10$ in (\ref{full_chi2}) with respect to $\mathsf{SNR}$ between $-10\,\mathrm{dB}$ to $10\,\mathrm{dB}$. For comparison, the exact curve for $p_e$, evaluated using Monte~Carlo simulation with $10^6$ trials, is also displayed.~$\diamond$ 

We end this subsection by observing that the exact value for $p_e$ in (\ref{WCom}) can be computed, though the computation is cumbersome. Since the nonzero eigenvalues for $A$ appear in pairs, then $\underline{\boldsymbol{x}}^{\mathsf{T}}A\underline{\boldsymbol{x}}=\sum_{i=1}^M (\psi_i^-\boldsymbol{\eta}_i^-+\psi_i^+\boldsymbol{\eta}_i^+)$, where $\boldsymbol{\eta}_1^-, \boldsymbol{\eta}_1^+,\dots, \boldsymbol{\eta}_M^-,\boldsymbol{\eta}_M^+$ are independent $\chi^2_2$ or equivalently, exponential random variables with parameter $\frac{1}{2}$. The PDF for $\psi_i^-\boldsymbol{\eta}_i^-+\psi_i^+\boldsymbol{\eta}_i^+$ is given by  
\begin{align}
\label{no_pdf}
p_{\psi_i^-\boldsymbol{\eta}_i^-+\psi_i^+\boldsymbol{\eta}_i^+}(z)=\frac{1}{2(\psi_i^+-\psi_i^-)}(e^{-\frac{z}{2\psi_i^-}}\mathds{1}_{z<0}+e^{-\frac{z}{2\psi_i^+}}\mathds{1}_{z>0}).
\end{align}
The PDF for $\underline{\boldsymbol{x}}^{\mathsf{T}}A\underline{\boldsymbol{x}}$ is the convolution of the PDFs in (\ref{no_pdf}) over $i=1,\dots, M$, which takes the form 
\begin{align}
\label{density_QF}
p_{\underline{\boldsymbol{x}}^{\mathsf{T}}A\underline{\boldsymbol{x}}}(z)=\sum_{i=1}^M(P_i^-(z)e^{-\frac{z}{2\psi_i^-}}\mathds{1}_{z<0}+P_i^+(z)e^{-\frac{z}{2\psi_i^+}}\mathds{1}_{z>0}),
\end{align}
where $P_i^{\pm}(z)$ are polynomials in $z$. As such, $p_e=\Pr(\underline{\boldsymbol{x}}^{\mathsf{T}}A\underline{\boldsymbol{x}}>0)=\int_{0}^\infty\sum_{i=1}^MP_i^+(z)e^{-\frac{z}{2\psi_i^+}}dz$, which can be evaluated analytically by applying the identity $\int_0^\infty z^ke^{-az}dz=\frac{k!}{a^{k+1}}$ for integer $k\ge0$ and $a>0$. A similar observation was made in \cite{Hochwald1} using an alternative approach based on characteristic functions, where $p_e$ is expressed as the sum of residues of a rational function in the complex plane. 

\section{Proof of Proposition~\ref{prop_1}}
\label{sec4}
The indispensable starting point in proving HWI is the exponential Markov inequality which gives 
\begin{align}
\Pr(\boldsymbol{\Delta}>t)\leq e^{-st}\mathbb{E}[e^{s\boldsymbol{\Delta}}],
\end{align} 
where $s\ge0$ is arbitrary. Let $A=UD U^{\mathsf{T}}$ be the eigenvalue decomposition of $A$ where $U$ is an orthogonal matrix and the diagonal matrix $D$ carries the (real) eigenvalues $\lambda_1, \dots, \lambda_n$ of~$A$. Then 
\begin{align}
\label{innooo}
\boldsymbol{\Delta}=\sum_{i=1}^n \lambda_i(\boldsymbol{\xi}^2_i-1),
\end{align}
where $\boldsymbol{\xi}_1,\dots, \boldsymbol{\xi}_n$ are independent standard normal random variables that are the entries of the random vector $\underline{\boldsymbol{\xi}}=U^{\mathsf{T}}\underline{\boldsymbol{x}}$. By independence,  
\begin{align}
\mathbb{E}[e^{s\boldsymbol{\Delta}}]=\prod_{i=1}^n e^{-s\lambda_i}\mathbb{E}[e^{s\lambda_i\boldsymbol{\xi}^2_i}].
\end{align} 
Since $\mathbb{E}[e^{s\lambda_i\boldsymbol{\xi}^2_i}]=\frac{1}{\sqrt{1-2s\lambda_i}}$ for $2s\lambda_i<1$, we get 
\begin{align}
\label{psd_1}
\Pr(\boldsymbol{\Delta}>t)\leq e^{-st}\prod_{i=1}^n e^{-s\lambda_i}e^{-\frac{1}{2}\ln(1-2s\lambda_i)},
\end{align}
for every $s\ge0$ such that $2s\lambda_i<1$ for all $1\leq i\le n$. Let $a>0$ and $0<b<1$ be such that (\ref{theone}) holds. The bound in~(\ref{psd_1}) can be relaxed as 
\begin{align}
\Pr(\boldsymbol{\Delta}>t)\leq e^{-st}\prod_{i=1}^ne^{-s\lambda_i}e^{\frac{1}{2}(2s\lambda_i+a(2s\lambda_i)^2)}=\exp\Big(-st+2as^2\sum_{i=1}^n\lambda_i^2\Big)=\exp\big(-st+2a\|A\|_2^2s^2\big),
\end{align} 
for every $s\ge0$ such that $|2s\lambda_i|\leq b$ for all $1\leq i\leq n$, or equivalently,  $0\leq s\leq \frac{b}{2\|A\|}$. Since $b<1$, the earlier condition $2s\lambda_i<1$ for all $1\leq i\le n$ is also met. The minimum value for the function $$f(s)=-st+2a\|A\|_2^2s^2$$ over $0\le s\le\frac{b}{2\|A\|}$ occurs at $s_{opt}=\min\{s_0, \frac{b}{2\|A\|}\}$ where $s_0=\frac{t}{4a\|A\|_2^2}$. We consider two cases:
\begin{enumerate}
  \item Let $s_0\leq \frac{b}{2\|A\|}$. Then 
  \begin{align}
  \label{fa1}
f(s_{opt})=f(s_0)=-\frac{t^2}{4a\|A\|_2^2}+2a\|A\|_2^2\Big(\frac{t}{4a\|A\|_2^2}\Big)^2=-\frac{t^2}{8a\|A\|_2^2}.
\end{align}
  \item Let $s_0>\frac{b}{2\|A\|}$. Then 
  \begin{align}
  \label{fa2}
f(s_{opt})=f\Big(\frac{b}{2\|A\|}\Big)=-\frac{bt}{2\|A\|}+2a\|A\|_2^2\Big(\frac{b}{2\|A\|}\Big)^2&\stackrel{}{<} -\frac{bt}{2\|A\|}+2a\|A\|_2^2 \times s_0\times\frac{b}{2\|A\|}\notag\\&=-\frac{bt}{2\|A\|}+\frac{t}{2}\times\frac{b}{2\|A\|}=-\frac{bt}{4\|A\|}.
\end{align}
\end{enumerate} 
By (\ref{fa1}) and (\ref{fa2}), 
\begin{align}
f(s_{opt})\leq  \max\Big\{-\frac{t^2}{8a\|A\|_2^2}, -\frac{bt}{4\|A\|}\Big\}=-\min\Big\{\frac{t^2}{8a\|A\|_2^2}, \frac{bt}{4\|A\|}\Big\}\leq -\kappa(a,b)\min\Big\{\frac{t^2}{\|A\|_2^2}, \frac{t}{\|A\|}\Big\},
\end{align}
where
\begin{align}
\kappa(a,b)=\min\Big\{\frac{1}{8a}, \frac{b}{4}\Big\}.
\end{align}
Next, let us explore the constants $a,b>0$. Choosing $a=1$ and $b=\frac{1}{2}$ results in $\kappa=\kappa(1,0.5)=0.125$ and recovers HWI as given in Reference~\cite{Giraud}. But, one can do better.  We need the special case of the following lemma for~$m=1$. 
\begin{lemma}
\label{lem_1}
Let $m\geq1$ be an integer. For every $0<b<1$, the smallest $a$ (tightest upper bound) for which~(\ref{gen_ineq1}) holds is $a=\theta_m(b)$ as defined in~(\ref{theta_func}). Equivalently, for every $a>\frac{1}{m+1}$, the largest $b$ for which the inequality in (\ref{theone}) holds is $\theta_m^{-1}(a)$ where $\theta_m^{-1}(\cdot)$ is the inverse function to $\theta_m(\cdot)$.  Conversely, if $a\leq \frac{1}{m+1}$, there exists no $b>0$ such that (\ref{theone}) holds.  
\end{lemma}
\begin{proof}
See Appendix~\ref{A1}.
\end{proof}
 Applying Lemma~\ref{lem_1} with $m=1$ results in 
\begin{align}
\kappa=\kappa(\theta_1(b),b)=\min\Big\{\frac{1}{8\theta_1(b)}, \frac{b}{4}\Big\}.
\end{align}
 for arbitrary $0<b<1$. Since $\frac{1}{8\theta_1(b)}=\frac{1}{-\frac{8}{b^2}(\ln(1-b)+b)}$ is decreasing in $b$ and $\frac{b}{4}$ is increasing in $b$,  the largest value for their minimum is achieved when $\frac{1}{8\theta_1(b)}=\frac{b}{4}$, i.e., $2b\theta_1(b)=1$. This equation is solved for $b\approx 0.583$ and the achieved maximum value is slightly larger than $0.1457$.

\section{Proof of Proposition~\ref{prop_2}}
\label{sec5}
Since $A$ is positive-semidefinite, its eigenvalues are all nonnegative. We begin with the inequality in~(\ref{psd_1}). Assume (\ref{thetwo}) holds for some $0<a< 1$ and $0<b<1$. Then 
\begin{align}
\label{shadow}
\Pr(\boldsymbol{\Delta}>t)\leq  e^{-st}\prod_{i=1}^n e^{-s\lambda_i}e^{\frac{1}{2}(2s\lambda_i+\frac{(2s\lambda_i)^2}{2(1-2as\lambda_i)})}= e^{-st}\prod_{i=1}^n e^{\frac{s^2\lambda^2_i}{1-2as\lambda_i}}=\exp\Big(-st+\sum_{i=1}^n\frac{s^2\lambda_i^2}{1-2as\lambda_i}\Big),
\end{align}
for every $s\geq 0$ such that $0\leq 2s\lambda_i\le b$ for all $1\le i\le n$, or equivalently, $0\leq s\leq \frac{b}{2\|A\|}$. Since $ab<1$ by assumption, then 
\begin{align}
\label{boul}
0<1-ab\leq 1-2as\lambda_i,
\end{align}
for every $1\leq i\leq n$. Therefore,  $0<1-2as\|A\|\leq 1-2as\lambda_i$ for every $1\leq i\leq n$ and we can loosen the upper bound in (\ref{shadow}) as 
\begin{align}
\label{mary}
\Pr(\boldsymbol{\Delta}>t)\leq  \exp\Big(-st+\sum_{i=1}^n\frac{s^2\lambda_i^2}{1-2as\|A\|}\Big)=\exp\Big(-st+\frac{\|A\|_2^2s^2}{1-2as\|A\|}\Big),\quad 0\leq s\leq \frac{b}{2\|A\|}.
\end{align}
Let $\alpha$ and $\beta$ be as in (\ref{conv}) and call the exponent on the right-hand side of (\ref{mary}) by $f(s)$, i.e.,  
\begin{align}
f(s)=-st+\frac{\beta s^2}{1-2a\alpha s}.
\end{align}
It is straightforward to check that the function $f$ is convex~(concave upward) over the entire interval $(-\infty,\frac{1}{2a\alpha})$ on the left side of its vertical asymptote $s=\frac{1}{2a\alpha}$. Moreover, it achieves an absolute minimum value over that interval at 
\begin{align}
s_0=\frac{1}{2a\alpha}\bigg(1-\frac{1}{\big(1+\frac{2a\alpha t}{\beta}\big)^{\frac{1}{2}}}\bigg).
\end{align}
Note that $\frac{b}{2\alpha}<\frac{1}{2a\alpha}$ due to the assumption $ab<1$. It follows that the expression on the right-hand side of (\ref{mary}) achieves its minimum value over $0\leq s\leq \frac{b}{2\alpha}$ at $s_{opt}=\min\{s_0, \frac{b}{2\alpha}\}$.
We consider two cases:
\begin{enumerate}
  \item Let $s_0\leq \frac{b}{2\alpha}$, or equivalently, 
  \begin{align}
  \label{c}
\frac{\alpha t}{\beta}\leq  \frac{1}{2a}\Big(\frac{1}{(1-ab)^2}-1\Big).
\end{align}
Then $f(s_{opt})=f(s_0)=-\frac{t}{2a\alpha}+\frac{\beta}{2a^2\alpha^2}\big((1+\frac{2a\alpha t}{\beta})^{\frac{1}{2}}-1\big)$,
where we omit the underlying algebra.  
  \item Let $s_0>\frac{b}{2\alpha}$. Then $f(s_{opt})=f(\frac{b}{2\alpha})=-\frac{bt}{2\alpha}+\frac{b^2\beta}{4(1-ab)\alpha^2}$.
\end{enumerate}
Let us denote the expression on the right-hand side of (\ref{c}) by $c$, i.e., 
\begin{align}
\label{tau}
c=\frac{1}{2a}\Big(\frac{1}{(1-ab)^2}-1\Big).
\end{align}
We have shown that 
\begin{align}
\label{imp_lm}
\Pr(\boldsymbol{\Delta}>t)\leq e^{f(s_{opt})}
\end{align}
and 
\begin{align}
\label{astar}
f(s_{opt})=\left\{\begin{array}{cc}
    -\frac{t}{2a\alpha}+\frac{\beta}{2a^2\alpha^2}\Big(\big(1+\frac{2a\alpha t}{\beta}\big)^{\frac{1}{2}}-1\Big)  & \frac{\alpha t}{\beta}\leq c   \\
    -\frac{bt}{2\alpha}+\frac{b^2\beta}{4(1-ab)\alpha^2}  &  \frac{\alpha t}{\beta}> c
\end{array}\right..
\end{align}
The next lemma explores possibilities for $a, b$. It addresses an inequality for which (\ref{thetwo}) is a special~case.  
\begin{lemma}
\label{lem_2}
Let $m\geq 1$ be an integer. For every $\frac{m+1}{m+2}< a<1$, the inequality 
 \begin{align}
  \label{gen_ineq2}
-\ln(1-x)\leq x+\frac{x^2}{2}+\frac{x^3}{3}+\cdots+\frac{x^m}{m}+\frac{x^{m+1}}{(m+1)(1-ax)},\quad 0\leq x\leq b,
\end{align}
 holds with $b=\frac{m+1-(m+2)a}{a(m-(m+1)a)}$.
\end{lemma} 
\begin{proof}
See Appendix~\ref{A2}.
\end{proof}
Applying Lemma~\ref{lem_2} with $m=1$, we obtain $b=\frac{3a-2}{a(2a-1)}$ for $\frac{2}{3}<a< 1$. Note that $b<1$ as promised earlier. Moreover, $1-ab=\frac{1-a}{2a-1}$ and the threshold $c$ in (\ref{tau}) is given~by $c=\frac{3a-2}{2(1-a)^2}$.
\section{Proof of Proposition~\ref{prop_3}}
\label{sec6}
To see how (\ref{en_LM}) implies HWI, let us write HWI in (\ref{HW}) as 
\begin{align}
\label{hw_diff}
\Pr(\boldsymbol{\Delta}>t)\leq e^{-\kappa\min\left\{\frac{t^2}{\beta} ,  \frac{t}{\alpha}\right\}},
\end{align}
where $\alpha, \beta$ are given in (\ref{conv}). Note 
\begin{align}
\label{emanuelle}
\min\left\{\frac{t^2}{\beta} ,  \frac{t}{\alpha}\right\}=\left\{\begin{array}{cc}
   \frac{t^2}{\beta}   & \frac{\alpha t}{\beta}\leq 1   \\
      \frac{t}{\alpha}   & \frac{\alpha t}{\beta} >1   
\end{array}\right..
\end{align}
We aim to show that there exists $\kappa>0$ such that the upper bound in (\ref{en_LM}) is less than or equal to the upper bound in (\ref{hw_diff}), or equivalently, $\Lambda(t,a)\leq -\kappa\min\{\frac{t^2}{\beta}, \frac{t}{\alpha}\}$. We want this to hold regardless of $t, \alpha, \beta$ and $\kappa$ to be as large as possible. Thus, we are looking to compute   
\begin{align}
\label{will}
\kappa= \inf_{t, \alpha, \beta> 0 }-\frac{\Lambda(t,a)}{\min\{\frac{t^2}{\beta}, \frac{t}{\alpha}\}}.
\end{align}
Interestingly, $-\frac{\Lambda(t,a)}{\min\{\frac{t^2}{\beta}, \frac{t}{\alpha}\}}$  depends on $t, \alpha, \beta$ only through the ratio $\rho=\frac{\alpha t}{\beta}$ in (\ref{rho_first}). For convenience in comparing the right-hand sides in (\ref{en_exp}) and (\ref{emanuelle}), we fix $a$ such that $c=1$, i.e., $\frac{3a-2}{2(1-a)^2}=1$. This quadratic equation has only one root in the admissible interval $\frac{2}{3}<a<1$. We denote this value of $a$ by $a_0$ given~by 
\begin{align}
\label{a0}
a_0=\frac{7-\sqrt{17}}{4}\approx 0.7192.
\end{align}
 One can easily check that the corresponding value for $b$ is $b_0=\frac{3a_0-2}{a_0(2a_0-1)}=\frac{1}{2}$. We have
\begin{align}
\label{ani}
-\frac{\Lambda(t,a_0)}{\min\{\frac{t^2}{\beta}, \frac{t}{\alpha}\}}=\left\{\begin{array}{cc}
     \frac{1}{2a_0\rho}-\frac{1}{2a_0^{2}\rho^2}\big(\sqrt{1+2a_0\rho}-1\big) &   \rho\leq 1 \\
   \frac{1}{4}-\frac{1}{16(1-\frac{a_0}{2})\rho}   &   \rho>1
\end{array}\right..
\end{align}
The continuous function of $\rho$ that sits on the right-hand side of (\ref{ani}) achieves its absolute minimum value at $\rho=1$. It follows that the infimum in (\ref{will}) is given by
\begin{align}
\frac{1}{4}-\frac{1}{16(1-\frac{a_0}{2})}=\frac{9-\sqrt{17}}{32}\approx 0.152.
\end{align}
It is shown in Appendix~\ref{A3} that the choice of $a=a_0$ in (\ref{a0}) indeed results in the largest value for $\kappa$ and that it was not just a choice of convenience. 

Finally, we address LMI itself. To see why (\ref{LM}) and (\ref{LM2}) are equivalent,  call $t'=2\sqrt{\beta t}+2\alpha t$. Solving for $t$, we get $t=(\frac{\sqrt{\beta+2\alpha t'}-\sqrt{\beta}}{2\alpha})^2$. Relabeling $t'$ as $t$, LMI in (\ref{LM}) becomes 
\begin{align}
\Pr(\boldsymbol{\Delta}>t)\leq \exp\left(-\Big(\frac{\sqrt{\beta+2\alpha t}-\sqrt{\beta}}{2\alpha}\Big)^2\right).
\end{align} 
 This inequality implies HWI in (\ref{hw_diff}) with 
\begin{align}
\label{bird}
\kappa= \inf_{t,\alpha, \beta>0}-\frac{(\frac{\sqrt{\beta+2\alpha t}-\sqrt{\beta}}{2\alpha})^2}{\min\{\frac{t^2}{\beta}, \frac{t}{\alpha}\}}.
\end{align}
Once again, the ratio on the right-hand side of (\ref{bird}) depends entirely on $\rho=\frac{\alpha t}{\beta}$. It is given by 
\begin{align}
-\frac{(\frac{\sqrt{\beta+2\alpha t}-\sqrt{\beta}}{2\alpha})^2}{\min\{\frac{t^2}{\beta}, \frac{t}{\alpha}\}}=\left\{\begin{array}{cc}
    (\frac{\sqrt{1+2\rho}-1}{2\rho})^2  & \rho\leq 1   \\
         \frac{(\sqrt{1+2\rho}-1)^2}{4\rho}  \rho>1
\end{array}\right..
\end{align}  
This function of $\rho$ achieves it absolute minimum value of $1-\frac{\sqrt{3}}{2}$ at $\rho=1$. 
\section{Proof of Proposition~\ref{prop_4}}
\label{sec7}
We begin with the bound in (\ref{psd_1}). By Lemma~\ref{lem_1} in Section~\ref{sec4}, for every integer $m\ge1$ and $0<b<1$, the inequality in (\ref{gen_ineq1}) is the tightest for $a=\theta_m(b)$.  Then 
\begin{align}
\label{arms}
\Pr(\boldsymbol{\Delta}>t)&\leq e^{-st}\prod_{i=1}^ne^{-s\lambda_i}e^{\frac{1}{2}(2s\lambda_i+a|2s\lambda_i|^{m+1}+\sum_{k=2}^m\frac{(2s\lambda_i)^k}{k})}\notag\\
&=\exp\Big(-st+a2^{m}s^{m+1}\sum_{i=1}^n|\lambda_i|^{m+1}+\sum_{k=2}^m\frac{2^{k-1}s^{k}}{k}\sum_{i=1}^n\lambda_i^{k}\Big)\notag\\
&=\exp\Big(-st+a2^{m}s^{m+1}\|A\|_{m+1}^{m+1}+\sum_{k=2}^m\frac{2^{k-1}s^{k}\mathrm{tr}(A^k)}{k}\Big),
\end{align}  
for every $s\geq0$ such that $|2s\lambda_i|\leq b$ holds for all $1\leq i\leq n$, or equivalently, $0\leq s\leq \frac{b}{2\|A\|}$. The value of $s$ that minimizes the right-hand side of (\ref{arms}) is not tractable for arbitrary integer $m$. Instead, let $s_0$ be the value of $s$ such that the sum of the first two terms in the exponent, i.e., $-st+a2^{m}s^{m+1}\|A\|_{m+1}^{m+1}$ is minimized. It is given~by  
\begin{align}
s_0=\frac{1}{2}\left(\frac{t}{(m+1) a\|A\|_{m+1}^{m+1}}\right)^{\frac{1}{m}}.
\end{align}  
We compute the right-hand side of (\ref{arms}) at $s=\hat{s}_{opt}$ defined by 
\begin{align}
\hat{s}_{opt}=\min\Big\{s_0, \frac{b}{2\|A\|}\Big\}.
\end{align}
Then the third term in the exponent in (\ref{arms}) is bounded by 
\begin{align}
\sum_{k=2}^m\frac{2^{k-1}\hat{s}_{opt}^{k}\mathrm{tr}(A^k)}{k} \leq  \sum_{k=2}^m\frac{2^{k-1}(\frac{b}{2\|A\|})^{k}|\mathrm{tr}(A^k)|}{k}=\sum_{k=2}^m\frac{b^k}{2k}\frac{|\mathrm{tr}(A^k)|}{\|A\|^k}\leq \frac{r(A)}{2}\sum_{k=2}^m \frac{b^k}{k},
\end{align} 
where the first step is due to $\hat{s}_{opt}^{k}\mathrm{tr}(A^k)\leq |\hat{s}_{opt}^{k}\mathrm{tr}(A^k)|$ and $|\hat{s}_{opt}|=\hat{s}_{opt}\leq \frac{b}{2\|A\|}$ and the last step is due~to
\begin{align}
|\mathrm{tr}(A^k)|=\big|\sum_{i=1}^n\lambda_i^k\big|\leq \sum_{i=1}^n|\lambda_i|^k\leq r(A)\max_{1\leq i\leq n}|\lambda_i|^k= r(A)\|A\|^k.
\end{align}
Here, we noted that $r(A)=|\{i: \lambda_i\neq 0\}|$.  Next, we bound the sum of the first two terms in the exponent, i.e., $-\hat{s}_{opt}t+a2^{m}\hat{s}_{opt}^{m+1}\|A\|_{m+1}^{m+1}$. We consider two cases:
\begin{enumerate}
  \item Let $s_0\leq \frac{b}{2\|A\|}$. Then $\hat{s}_{opt}=s_0$ and simple algebra shows that
  \begin{align}
  \label{yani1}
-\hat{s}_{opt}t+a2^{m}\hat{s}_{opt}^{m+1}\|A\|_{m+1}^{m+1}=-\kappa_{m,1}\frac{t^{1+\frac{1}{m}}}{\|A\|_{m+1}^{1+\frac{1}{m}}},
\end{align}
where $\kappa_{m,1}$ is defined by 
\begin{align}
\kappa_{m,1}=\frac{1-\frac{1}{m+1}}{2(m+1)^{\frac{1}{m}}a^{\frac{1}{m}}}.
\end{align}
  \item Let $s_0> \frac{b}{2\|A\|}$. Then $\hat{s}_{opt}=\frac{b}{2\|A\|}$ and we have 
  \begin{align}
  \label{yani2}
-\hat{s}_{opt}t+a2^{m}\hat{s}_{opt}^{m+1}\|A\|_{m+1}^{m+1} &=-\hat{s}_{opt}t+a2^m \hat{s}_{opt}\times \hat{s}_{opt}^m\|A\|_{m+1}^{m+1}\notag\\
&\leq -\frac{bt}{2\|A\|}+a2^m \frac{b}{2\|A\|}\times s_0^m\|A\|_{m+1}^{m+1}\notag\\
&=-\frac{bt}{2\|A\|}+a2^m \frac{b}{2\|A\|}\times \frac{t}{2^m(m+1)a\|A\|_{m+1}^{m+1}} \|A\|_{m+1}^{m+1}\notag\\
&=-\kappa_{m,2}\frac{t}{\|A\|},
\end{align} 
where the second step is due to $\hat{s}_{opt}\leq s_0$ and $\kappa_{m,2}$ is defined by 
\begin{align}
\kappa_{m,2}=\frac{b}{2}\Big(1-\frac{1}{m+1}\Big).
\end{align}
\end{enumerate}
By (\ref{yani1}) and (\ref{yani2}), 
\begin{align}
-\hat{s}_{opt}t+a2^{m}\hat{s}_{opt}^{m+1}\|A\|_{m+1}^{m+1} &\leq  \max\Big\{-\kappa_{m,1}\frac{t^{1+\frac{1}{m}}}{\|A\|_{m+1}^{1+\frac{1}{m}}}, -\kappa_{m,2}\frac{t}{\|A\|} \Big\}\notag\\
&=-\min\Big\{\kappa_{m,1}\frac{t^{1+\frac{1}{m}}}{\|A\|_{m+1}^{1+\frac{1}{m}}}, \kappa_{m,2}\frac{t}{\|A\|} \Big\}\notag\\
&\leq -\min\{\kappa_{m,1}, \kappa_{m,2}\}\min\bigg\{\frac{t^{1+\frac{1}{m}}}{\|A\|_{m+1}^{1+\frac{1}{m}}}, \frac{t}{\|A\|} \bigg\}.
\end{align}
The proof of Proposition~\ref{prop_4} is now complete.

\section{Proof of Proposition~\ref{prop_6}}
\label{sec9}
Let us write the $m_\infty$~bound in (\ref{infty_m}) and the relaxed HW bound in (\ref{relaxed_HW}) in terms of  $r=\frac{t}{\|A\|}$. The $m_\infty$~bound becomes 
\begin{align}
\label{whistle}
 \Pr(\boldsymbol{\Delta}>t)\leq \Big(1+\frac{r}{n}\Big)^{\frac{n}{2}}e^{-\frac{r}{2}}=e^{-(\frac{r}{2}-\frac{n}{2}\ln(1+\frac{r}{n}))}.
\end{align}  
and the relaxed HWI becomes 
\begin{align}
\Pr(\boldsymbol{\Delta}>t)\leq e^{-\kappa\min\{\frac{r^2}{n}, r\}},\,\,\, \kappa=\frac{9-\sqrt{17}}{32}.
\end{align}
The $m_\infty$~bound is tighter than the relaxed HWI if and only if 
\begin{align}
\label{arezoome}
\frac{r}{2}-\frac{n}{2}\ln\Big(1+\frac{r}{n}\Big)>\kappa\min\Big\{\frac{r^2}{n}, r\Big\},\,\,\, r>0.
\end{align}
We consider two cases:
\begin{enumerate}
  \item Let $0<r\leq n$. Then $\min\{\frac{r^2}{n}, r\}=\frac{r^2}{n}$ and (\ref{arezoome}) becomes 
  \begin{align}
\frac{r}{2}-\frac{n}{2}\ln\Big(1+\frac{r}{n}\Big)>\frac{\kappa r^2}{n}.
\end{align}
Dividing both sides by $\frac{n}{2}$ and rearranging terms,
\begin{align}
\label{leave}
\ln\Big(1+\frac{r}{n}\Big)<\frac{r}{n}\Big(1-\frac{2\kappa r}{n}\Big).
\end{align} 
Denoting $x=\frac{r}{n}$, we need to show that $\ln(1+x)<x(1-2\kappa x)$ for all $0<x\leq 1$. Define $f(x)=x(1-2\kappa x)-\ln(1+x)$. Then $f'(x)=1-4\kappa x-\frac{1}{1+x}=\frac{(1-4\kappa)x-4\kappa x^2}{1+x}$. The numerator is the quadratic function $q(x)=(1-4\kappa)x-4\kappa x^2$. Solving $q(x)=0$, we see that $f(x)$ has only one critical number inside the interval $(0,1)$ given by $x_0=\frac{1-4\kappa}{4\kappa}=\frac{1+\sqrt{17}}{8}\approx 0.6404$. Moreover, $f'(x)>0$ for $0<x<x_0$ and $f'(x)<0$ for $x_0<x<1$. By the first derivative test, the function $f$ increases over $(0, x_0)$, achieves a local maximum value at $x=x_0$ and then decreases over $(x_0,1)$. As such, $f$ achieves it absolute minimum value over the interval $[0,1]$ either at $x=0$ or at $x=1$. But, $f(0)=0$ and $f(1)=1-2\kappa-\ln2\approx 0.002>0$. It follows that $f(x)>0$ over $(0,1]$ as desired. 
  \item Let $r>n$. Then $\min\{\frac{r^2}{n}, r\}=r$ and (\ref{arezoome}) becomes $\frac{r}{2}-\frac{n}{2}\ln(1+\frac{r}{n})>\kappa r$. Dividing both sides by $\frac{n}{2}$ and rearranging terms, we get the inequality
  \begin{align}
\ln\Big(1+\frac{r}{n}\Big)<(1-2\kappa)\frac{r}{n}.
\end{align}
Denoting $x=\frac{r}{n}$, we show that $\ln(1+x)<(1-2\kappa)x$ for every $x>1$. Define the function $f(x)=(1-2\kappa) x-\ln(1+x)$. Then $f'(x)=1-2\kappa-\frac{1}{1+x}>1-2\kappa-\frac{1}{2}=\frac{\sqrt{17}-1}{16}>0$ where we used $-\frac{1}{1+x}>-\frac{1}{2}$ for $x>1$. Since $f(1)=1-2\kappa-\ln2\approx 0.002>0$ and $f$ is increasing over $[1,\infty)$, it follows that $f(x)>0$ for every~$x>1$ as desired. 
\end{enumerate}
Next, we look at the relaxed LMI in (\ref{relaxed_LM}). Writing the bound in terms of $r=\frac{t}{\|A\|}$ and expanding the complete square in the exponent, we have 
\begin{align}
\Pr(\boldsymbol{\Delta}>t)\leq e^{-\frac{r}{2}+\frac{n}{2}(\sqrt{1+\frac{2r}{n}}-1)}.
\end{align}
Therefore, the upper bound in the $m_\infty$~inequality in (\ref{whistle}) is less than the upper bound in the relaxed LMI if and only if 
\begin{align}
\ln\Big(1+\frac{r}{n}\Big)<\Big(1+\frac{2r}{n}\Big)^{1/2}-1.
\end{align} 
However, this holds true for all $r>0$ due to the inequality $\ln(1+x)<\sqrt{1+2x}-1$ for every $x>0$. To see this, let $f(x)=\sqrt{1+2x}-1-\ln(1+x)$. Then $f'(x)=\frac{1}{\sqrt{1+2x}}-\frac{1}{1+x}$. But, $\sqrt{1+2x}<1+x$ for every $x>0$. As such, $f'(x)>0$ for every $x>0$. Since $f(0)=0$ and $f$ is increasing over $[0,\infty)$, then $f(x)>f(0)=0$ for every $x>0$.

Finally, we verify that the $m_\infty$~bound always beats the large-deviations bound. Call $x=\frac{1+r}{n}\geq 1$. The bound in (\ref{ld_chi2}) is looser than that in (\ref{minfty_pol}) if and only if 
\begin{align}
\label{paq}
\Big(1+\frac{r}{n}\Big)^{\frac{n}{2}}<\frac{1}{\sqrt{e}}\Big(\frac{e(1+r)}{n}\Big)^{\frac{n}{2}}\iff \Big(1+\frac{nx-1}{n}\Big)^{\frac{n}{2}}<\frac{1}{\sqrt{e}}(ex)^{\frac{n}{2}}&\stackrel{}{\iff} 1+\frac{nx-1}{n}<\frac{1}{e^{\frac{1}{n}}}ex\notag\\
&\iff 1-\frac{1}{n}+x<e^{1-\frac{1}{n}}x\notag\\
&\iff 1-\frac{1}{n}<(e^{1-\frac{1}{n}}-1)x.
\end{align} 
Using the inequality $z<e^z-1$ for $z\neq 0$,  we have $1-\frac{1}{n}<e^{1-\frac{1}{n}}-1$ for every $n\geq 2$. Therefore, (\ref{paq}) is certainly true due to $x\geq 1$.

\section{Proof of Proposition~\ref{prop_7} }
\label{sec10}
Recall the polynomials $P_n^{(1)}(r)$ and $P_n^{(2)}(r)$ in (\ref{poly_1}) and (\ref{poly_2}), respectively. It is straightforward to check that $P^{(2)}_n(r)< P_n^{(1)}(r)$ for all $r\geq 0$ if $n=2,4,6$. In fact, 
\begin{align}
\label{n_2}
P_2^{(1)}(r)-P_2^{(2)}(r)=\frac{r}{2}+1-\frac{1}{\sqrt{e}},
\end{align}
\begin{align}
\label{n_4}
P_4^{(1)}(r)-P_4^{(2)}(r)=\frac{r^2}{16}+\Big(\frac{1}{2}-\frac{1}{2\sqrt{e}}\Big)r+1-\frac{3}{2\sqrt{e}}
\end{align}
and 
\begin{align}
\label{n_6}
P_6^{(1)}(r)-P_6^{(2)}(r)=\frac{r^3}{216}+\Big(\frac{1}{12}-\frac{1}{8\sqrt{e}}\Big)r^2+\Big(\frac{1}{2}-\frac{3}{4\sqrt{e}}\Big)r+1-\frac{13}{8\sqrt{e}}.
\end{align}
The polynomials on the right-hand sides in (\ref{n_2}), (\ref{n_4}) and (\ref{n_6}) all have positive coefficients. This proves supremacy of the weak $\chi^2$~bound over the $m_\infty$~bound when $n=2,4,6$. If $n$ is an even integer greater than or equal to $8$, the proof strategy is to first show that the polynomial $P_n^{(1)}(r)-P_n^{(2)}(r)$ has at least two positive roots and then proceed to prove that $P_n^{(1)}(r)-P_n^{(2)}(r)$ can not have more than two positive roots. We observe that $P_n^{(1)}(0)>P_n^{(2)}(0)$ regardless of the value of $n$. In fact, 
\begin{align}
\label{at_0}
P_n^{(1)}(0)=1,\quad P_n^{(2)}(0)=\frac{1}{\sqrt{e}}\sum_{i=0}^{\frac{n}{2}-1}\frac{1}{2^ii!}<\frac{1}{\sqrt{e}}\sum_{i=0}^{\infty}\frac{(\frac{1}{2})^i}{i!}=1.
\end{align}
Since $\deg (P_n^{(1)})=\frac{n}{2}$ is larger than $\deg (P_n^{(2)})=\frac{n}{2}-1$, we also have $P_n^{(1)}(r)>P_n^{(2)}(r)$ for all sufficiently large $r$. If we can find a positive number $r^*$ such that $P_n^{(1)}(r^*)<P_n^{(2)}(r^*)$, then the Intermediate Value Theorem~(IVT) implies that there must exist numbers $0<r_n<r^*<r'_n$ such that $P_n^{(1)}(r_n)=P_n^{(2)}(r_n)$ and $P_n^{(1)}(r'_n)=P_n^{(2)}(r'_n)$, i.e., $r_n, r'_n$ are zeros of $P_n^{(1)}(r)-P_n^{(2)}(r)$. We verify that the choice $r^*=1$ works for every even integer $n\geq 8$. Let us state this as a lemma. 
\begin{lemma}
\label{lem_lem}
Let the polynomials $P_n^{(1)}(r)$ and $P_n^{(2)}(r)$ be as in (\ref{poly_1}) and (\ref{poly_2}), respectively. Then $P_n^{(1)}(1)<P_n^{(2)}(1)$ for every even integer $n\ge8$, i.e., 
\begin{align}
\label{love_1}
\Big(1+\frac{1}{n}\Big)^{\frac{n}{2}}<\frac{1}{\sqrt{e}}\sum_{i=0}^{\frac{n}{2}-1}\frac{1}{i!},\,\,\, n\geq 8.
\end{align}
\end{lemma}
\begin{proof}
Denote $m=\frac{n}{2}$. We need to show that 
\begin{align}
\label{love_00}
\Big(1+\frac{1}{2m}\Big)^{m}<\frac{1}{\sqrt{e}}\sum_{i=0}^{m-1}\frac{1}{i!},\,\,\, m\ge4.
\end{align}
The cases $4\leq m\leq 49$ and $m\geq 50$ are addressed separately. The reason for this split becomes clear by the end of the proof. One can easily check with the help of a computer that (\ref{love_00}) holds for $4\leq m\leq 49$. 
For the remainder of the proof, we assume $m\geq 50$. Let $\boldsymbol{x}\sim \mathrm{Poisson}(1)$ be a Poisson random variable with parameter $1$. Then 
\begin{align}
\label{love_2}
\frac{1}{\sqrt{e}}\sum_{i=0}^{m-1}\frac{1}{i!} &=\sqrt{e}\sum_{i=0}^{m-1}\frac{e^{-1}1^i}{i!}=\sqrt{e}\Pr(\boldsymbol{x}\leq m-1)=\sqrt{e}\big(1-\Pr(\boldsymbol{x}\geq m)\big)\notag\\
&=\sqrt{e}\big(1-\Pr(\boldsymbol{x}^2\geq m^2)\big)\ge \sqrt{e}\Big(1-\frac{\mathbb{E}[\boldsymbol{x}^2]}{m^2}\Big)=\sqrt{e}\Big(1-\frac{2}{m^2}\Big),
\end{align}
where the penultimate step is due to Markov's inequality and the last step is due to $\mathbb{E}[\boldsymbol{x}^2]=\mathrm{Var}(\boldsymbol{x})+(\mathbb{E}[\boldsymbol{x}])^2=1+1^2=2$. By (\ref{love_1}) and (\ref{love_2}), it suffices to show that 
\begin{align}
\label{love_3}
\Big(1+\frac{1}{2m}\Big)^{m}<\sqrt{e}\Big(1-\frac{2}{m^2}\Big),\,\,\, m\geq 50.
\end{align}
Define 
\begin{align}
f(x)=\frac{1}{2}+\ln\Big(1-\frac{2}{x^2}\Big)-x\ln\Big(1+\frac{1}{2x}\Big),\,\,\, x\geq 50.
\end{align}
The inequality in (\ref{love_3}) can be written as $e^{f(m)}>1$ for $m\geq 50$. Thus, we aim to show that 
\begin{align}
\label{stumble_oo}
f(x)>0,\,\,\, x\ge50. 
\end{align}
To achieve this goal, we prove that
\begin{align}
\label{love_4_1}
\lim_{x\to\infty}f(x)=0 
\end{align}
and
\begin{align}
\label{love_4_2}
f'(x)<0,\,\,\, x\geq50.
\end{align}
 Verifying (\ref{love_4_1}) is an easy task. To prove~(\ref{love_4_2}), we use a similar method. We show that 
\begin{align}
\label{love_5}
\lim_{x\to\infty}f'(x)=0 
\end{align}
and 
\begin{align}
\label{love_6}
f''(x)>0,\,\,\, x\geq50.
\end{align}
We have 
\begin{align}
f'(x)=\frac{4}{x(x^2-2)}+\frac{1}{2x+1}-\ln\Big(1+\frac{1}{2x}\Big).
\end{align}
Then (\ref{love_5}) clearly holds. To verify (\ref{love_6}), we compute 
\begin{align}
\label{love_7}
f''(x)=\frac{x^5-48x^4-52x^3+20x^2+36x+8}{x^2(x^2-2)^2(2x+1)^2},
\end{align}
where we omit the tedious algebra. The denominator in (\ref{love_7}) is clearly positive for every $x\geq 50$. Let us denote the numerator by $p(x)$, i.e.,
\begin{align}
p(x)=x^5-48x^4-52x^3+20x^2+36x+8.
\end{align} 
To show that $p(x)>0$ for every $x\geq 50$, let us look at the shifted polynomial $p(x+50)$. It is given by 
\begin{align}
p(x+50)=x^5 + 202x^4 + 15348x^3 + 522220x^2 + 6862036x + 6051808.
\end{align}
Since all coefficients are positive in $p(x+50)$, we conclude that $p(x+50)\geq6051808>0$ for all~$x\geq 0$ and hence, $p(x)> 0$ for all $x\geq 50$ as desired. One can check that $p(49)=-303155<0$. This is why we addressed $4\leq m\leq 49$ and $m\ge50$ separately.
 \end{proof}
 Next, we prove that $r_n, r'_n$ are the \textit{only} two positive roots for the polynomial $P_n^{(1)}(r)-P_n^{(2)}(r)$. Let us recall \emph{Descartes' rule of signs}. It says that for a polynomial $a_0+a_1x+a_2x^2+\cdots+a_dx^d$ of degree~$d$, the number $N$ of positive roots and the number $N'$ of sign changes in the coefficients have the same parity and $N\leq N'$. The number $N'$ is the number of indices $0\leq j\leq d-1$ such that $a_ja_{j+1}<0$. The rest of proof is devoted to showing that there are only two sign changes in the coefficients of $P_n^{(1)}(r)-P_n^{(2)}(r)$ for every even integer $n\ge8$. Since we already showed the existence of at least two positive roots $r_n, r'_n$, it follows by Descartes' rule of signs that $r_n, r'_n$ are the only two positive roots for $P_n^{(1)}(r)-P_n^{(2)}(r)$. Let us first determine the coefficients of $P_n^{(1)}(r)$ and $P_n^{(2)}(r)$. In case of $P_n^{(1)}(r)$, they are given by 
\begin{align}
P_n^{(1)}(r)=\sum_{j=0}^{\frac{n}{2}}a^{(1)}_jr^j,\quad a^{(1)}_j=\frac{{\frac{n}{2}\choose j}}{n^j}.
\end{align}   
In case of $P_n^{(2)}(r)$, we write  
\begin{align}
P_n^{(2)}(r)=\frac{1}{\sqrt{e}}\sum_{i=0}^{\frac{n}{2}-1}\frac{1}{2^ii!}\sum_{j=0}^i{i\choose j}r^j=\frac{1}{\sqrt{e}}\sum_{j=0}^{\frac{n}{2}-1}\sum_{i=j}^{\frac{n}{2}-1}\frac{{i\choose j}}{2^ii!}r^j,
\end{align}
where in the second step we have changed the order of summations. This shows that 
\begin{align}
P_n^{(2)}(r)=\sum_{j=0}^{\frac{n}{2}}a^{(2)}_jr^j,\quad a^{(2)}_j=\left\{\begin{array}{cc}
    \frac{1}{\sqrt{e}}\sum_{i=j}^{\frac{n}{2}-1}\frac{{i\choose j}}{2^ii!}  &  0\leq j\leq \frac{n}{2}-1   \\
     0 &  j=\frac{n}{2} 
\end{array}\right..
\end{align}
The coefficients of $P_n^{(1)}(r)-P_n^{(2)}(r)$ are $a_j^{(1)}-a_j^{(2)}$ for $j=0,1,\dots, \frac{n}{2}$. These coefficients incur a total of two sign changes. More precisely, we show
\begin{align}
\label{t1}
a^{(1)}_0-a^{(2)}_0>0,\quad a^{(1)}_1-a^{(2)}_1>0,
\end{align} 
\begin{align}
\label{t2}
a^{(1)}_j-a^{(2)}_j<0,\,\,\, j=2,3,\dots, \frac{n}{2}-1
\end{align}
and 
\begin{align}
\label{t3}
a^{(1)}_{\frac{n}{2}}-a^{(2)}_{\frac{n}{2}}>0.
\end{align} 
Verifying (\ref{t1}) is trivial. In fact, $a_0^{(1)}=P_n^{(1)}(0), a_0^{(2)}=P_n^{(2)}(0)$ and we checked in (\ref{at_0}) that $P_n^{(1)}(0)>P_n^{(2)}(0)$. Also, 
\begin{align}
a^{(1)}_1=\frac{\frac{n}{2}}{n}=\frac{1}{2}\end{align}
and 
\begin{align}
a^{(2)}_1=\frac{1}{\sqrt{e}}\sum_{i=1}^{\frac{n}{2}-1}\frac{i}{2^ii!}=\frac{1}{\sqrt{e}}\sum_{i=1}^{\frac{n}{2}-1}\frac{1}{2^i(i-1)!}=\frac{1}{\sqrt{e}}\sum_{i=0}^{\frac{n}{2}-2}\frac{1}{2^{i+1}i!}=\frac{1}{2\sqrt{e}}\sum_{i=0}^{\frac{n}{2}-2}\frac{1}{2^{i}i!}<\frac{1}{2\sqrt{e}}\sum_{i=0}^{\infty}\frac{(\frac{1}{2})^i}{i!}=\frac{1}{2}.
\end{align}
The inequality in (\ref{t3}) is also trivial as $a_{\frac{n}{2}}^{(2)}=0<a_{\frac{n}{2}}^{(1)}=\frac{1}{n^{\frac{n}{2}}}$. To prove (\ref{t2}), first  we write it in a different form.  The condition $a^{(1)}_j<a^{(2)}_j$ for $j=2,3,\dots, \frac{n}{2}-1$ can be written as 
\begin{align}
\frac{{\frac{n}{2}\choose j}}{n^j}<\frac{1}{\sqrt{e}}\sum_{i=j}^{\frac{n}{2}-1}\frac{{i\choose j}}{2^i i!}&\stackrel{(*)}{\iff} \frac{{\frac{n}{2}\choose j}}{n^j}<\frac{1}{\sqrt{e}}\sum_{i=0}^{\frac{n}{2}-j-1}\frac{{i+j\choose j}}{2^{i+j} (i+j)!}\notag\\
&\iff \frac{\frac{\frac{n}{2}!}{j! (\frac{n}{2}-j)!}}{n^j}<\frac{1}{\sqrt{e}}\sum_{i=0}^{\frac{n}{2}-j-1}\frac{\frac{(i+j)!}{j! i!}}{2^{i+j}(i+j)!}\notag\\
&\iff \frac{\frac{\frac{n}{2}!}{ (\frac{n}{2}-j)!}}{n^j}<\frac{1}{\sqrt{e}}\sum_{i=0}^{\frac{n}{2}-j-1}\frac{1}{2^{i+j}i!}\notag\\
&\iff \frac{\frac{\frac{n}{2}!}{ (\frac{n}{2}-j)!}}{(\frac{n}{2})^j}<\frac{1}{\sqrt{e}}\sum_{i=0}^{\frac{n}{2}-j-1}\frac{1}{2^{i}i!}\notag\\
&\iff \frac{\frac{n}{2}(\frac{n}{2}-1)(\frac{n}{2}-2)\cdots(\frac{n}{2}-j+1)}{(\frac{n}{2})^j}<\frac{1}{\sqrt{e}}\sum_{i=0}^{\frac{n}{2}-j-1}\frac{1}{2^ii!},
\end{align}
where in $(*)$ we have changed the index $i$ to $i+j$. Call $m=\frac{n}{2}$. Then we need to show 
\begin{align}
\label{bell}
\frac{(m-1)(m-2)\cdots(m-j+1)}{m^{j-1}}<\frac{1}{\sqrt{e}}\sum_{i=0}^{m-j-1}\frac{1}{2^ii!},\,\,\, 2\leq j\leq m-1,
\end{align}
holds for all $m\geq 4$. We address $m=4$ and $m\geq 5$ separately: 
\begin{enumerate}
  \item Let $m=4$. The possibilities for $j$ are  $j=2,3$. If $j=2$, then~(\ref{bell}) becomes $\frac{3}{4}<\frac{1}{\sqrt{e}}(1+\frac{1}{2})$ which is true. If $j=3$, then~(\ref{bell}) becomes $\frac{3\times 2}{16}<\frac{1}{\sqrt{e}}\times 1$ which is also true.
  \item Assume $m\geq 5$. Let $\boldsymbol{x}\sim \mathrm{Poisson}(\frac{1}{2})$ be a Poisson random variable with parameter $\frac{1}{2}$. Then the right-hand side of (\ref{bell}) can be written~as 
\begin{align}
\label{bell_2}
\frac{1}{\sqrt{e}}\sum_{i=0}^{m-j-1}\frac{1}{2^ii!} &=\sum_{i=0}^{m-j-1}\frac{e^{-\frac{1}{2}}(\frac{1}{2})^i}{i!}=\Pr(\boldsymbol{x}\leq m-j-1)=1-\Pr(\boldsymbol{x}\geq m-j)\notag\\&\stackrel{(a)}{\geq}  1-\frac{\mathbb{E}[\boldsymbol{x}]}{m-j}\stackrel{(b)}{=} 1-\frac{\frac{1}{2}}{m-j}=\frac{m-j-\frac{1}{2}}{m-j},
\end{align}
where~$(a)$ is due to Markov's inequality and $(b)$ is due to $\mathbb{E}[\boldsymbol{x}]=\frac{1}{2}$. By (\ref{bell}) and (\ref{bell_2}), it is enough to show that for every $m\geq 5$, 
\begin{align}
\label{bell_3}
\frac{(m-1)(m-2)\cdots(m-j+1)}{m^{j-1}}<\frac{m-j-\frac{1}{2}}{m-j},\,\,\, j=2,3,\dots, m-1.
\end{align}
For $j=2$, (\ref{bell_3}) becomes $\frac{m-1}{m}<\frac{m-\frac{5}{2}}{m-2}$ which is equivalent to $(m-1)(m-2)<m(m-\frac{5}{2})$. This simplifies to $m>4$ which is true. This is why we separated the cases $m=4$ and $m>4$. If $j\geq 3$, we write~(\ref{bell_3}) as 
\begin{align}
\frac{(m-1)(m-2)\cdots(m-j+2)}{m^{j-2}}\times \frac{m-j+1}{m}<\frac{m-j-\frac{1}{2}}{m-j},
\end{align}
or equivalently, 
\begin{align}
\Big(1-\frac{1}{m}\Big)\Big(1-\frac{2}{m}\Big)\cdots\Big(1-\frac{j-2}{m}\Big)  \frac{m-j+1}{m}<\frac{m-j-\frac{1}{2}}{m-j}.
\end{align}
Since $(1-\frac{1}{m})(1-\frac{2}{m})\cdots(1-\frac{j-2}{m})<1$, it suffices to prove 
\begin{align}
 \frac{m-j+1}{m}<\frac{m-j-\frac{1}{2}}{m-j},\,\,\, j=3,\dots, m-1.
\end{align}
This is equivalent to $(m-j+1)(m-j)<m(m-j-\frac{1}{2})$ which simplifies to 
\begin{align}
\label{folk}
m>\frac{j(j-1)}{j-\frac{3}{2}},\,\,\, j=3,\dots, m-1.
\end{align}
The function $f(x)=\frac{x(x-1)}{x-\frac{3}{2}}$ is monotone increasing for $x\geq 3$. Therefore, (\ref{folk}) holds if and only if $m>f(m-1)=\frac{(m-1)(m-1-1)}{m-\frac{3}{2}}$. This simplifies to $m>\frac{4}{3}$ which is true. The proof of Proposition~\ref{prop_7} is now complete. 
\end{enumerate}

\section{Proof of Proposition~\ref{prop_17}}
\label{sec11}
\textcolor{Black}{Without loss of generality, we assume $A$ is full-rank, i.e., $r(A)=n$.  As explained in (\ref{innooo}), we can write $\boldsymbol{\Delta}=\sum_{i=1}^n\lambda_{i}\boldsymbol{\xi}_{i}^2-\sum_{i=1}^n\lambda_{i}$, where $\lambda_{i}$ are the (real and nonzero) eigenvalues of $A$ and $\boldsymbol{\xi}_{i}$ are independent standard normal random variables. Denote $\underline{\boldsymbol{\xi}}_k=(\boldsymbol{\xi}_{k+1},\cdots, \boldsymbol{\xi}_n)^{\mathsf{T}}$, $D_k=\mathrm{diag}(\lambda_{k+1},\cdots, \lambda_n)$ and $\boldsymbol{\Delta}_k=\underline{\boldsymbol{\xi}}^{\mathsf{T}}_kD_k\underline{\boldsymbol{\xi}}_k-\mathrm{tr}(D_k)$. We have 
\begin{align}
\label{koochil_1}
\Pr(\boldsymbol{\Delta}>t)&\stackrel{(a)}{=} \Pr\Big(\boldsymbol{\Delta}_k>t+b_k-\sum_{i=1}^{k}\lambda_{i}\boldsymbol{\xi}^2_{i}\Big)\notag\\
&\stackrel{(b)}{=} \mathbb{E}_{z_i\sim \boldsymbol{\xi}_i^2, 1\le i\le k}\Big[\Pr\Big(\boldsymbol{\Delta}_k>t+b_k-\sum_{i=1}^{k}\lambda_{i}\boldsymbol{\xi}^2_{i}\,\big|\, \boldsymbol{\xi}_1^2=z_1,\cdots, \boldsymbol{\xi}^2_k=z_k\Big)\Big]\notag\\
&\stackrel{(c)}{=} \mathbb{E}_{z_i\sim \boldsymbol{\xi}_i^2, 1\le i\le k}\Big[\Pr\Big(\boldsymbol{\Delta}_k>t+b_k-\sum_{i=1}^{k}\lambda_{i}z_{i}\Big)\Big],
\end{align}
where $(a)$ is due to $\boldsymbol{\Delta}=\sum_{i=1}^k \lambda_{i}\boldsymbol{\xi}_{i}^2-b_k+\boldsymbol{\Delta}_k$, $(b)$ is due to conditioning on the event $\{\boldsymbol{\xi}_1^2=z_1,\cdots, \boldsymbol{\xi}^2_k=z_k\}$ for $z_1,\dots, z_k\ge0$ and $(c)$ is due to independence of $\boldsymbol{\Delta}_k$ and $(\boldsymbol{\xi}_1,\cdots, \boldsymbol{\xi}_k)$. 
By~(\ref{poloi_1}) and (\ref{poloi_2}), $t+b_k-\sum_{i=1}^{k}\lambda_{i}z_{i}\ge0$ and hence, the original $m_\infty$~inequality gives 
\begin{align}
\label{koochil_3}
\Pr\big(\boldsymbol{\Delta}_k>t+b_k-\sum_{i=1}^{k}\lambda_{i}z_{i}\big)\leq \Big(1+\frac{ t+b_k-\sum_{i=1}^{k}\lambda_{i}z_{i}}{ (n-k)a_k}\Big)^{\frac{n-k}{2}}e^{-\frac{t+b_k-\sum_{i=1}^{k}\lambda_{i}z_{i}}{2a_k}}.
\end{align}
Since $n-k$ is even, one can use the multinomial expansion to write 
\begin{align}
\label{koochil_4}
\Big(1+\frac{ t+b_k-\sum_{i=1}^{k}\lambda_{i}z_{i}}{ (n-k)a_k}\Big)^{\frac{n-k}{2}}=\sum_{\substack{j_0,j_1,\dots, j_k\ge0\\ j_0+j_1+\cdots+j_k=\frac{n-k}{2}}}\frac{(\frac{n-k}{2})!}{j_0!j_1!\cdots j_{k}!}\Big(1+\frac{t+b_k}{(n-k)a_k}\Big)^{j_0}\prod_{i=1}^k\Big(\frac{-\lambda_iz_i}{(n-k)a_k}\Big)^{j_i}.
\end{align} 
 Putting (\ref{koochil_1}), (\ref{koochil_3}) and (\ref{koochil_4}) together,
\begin{align}
\label{monkey_1}
\Pr(\boldsymbol{\Delta}>t)&\leq\mathbb{E}\Big[\Big(1+\frac{ t+b_k-\sum_{i=1}^{k}\lambda_{i}\boldsymbol{\xi}^2_{i}}{ (n-k)a_k}\Big)^{\frac{n-k}{2}}e^{-\frac{t+b_k-\sum_{i=1}^{k}\lambda_{i}\boldsymbol{\xi}^2_{i}}{2a_k}}\Big]\notag\\
&\leq e^{-\frac{t+b_k}{2a_k}} \sum_{\substack{j_0,j_1,\dots, j_k\ge0\\ j_0+j_1+\cdots+j_k=\frac{n-k}{2}}}\frac{(\frac{n-k}{2})!}{j_0!j_1!\cdots j_{k}!}\Big(1+\frac{t+b_k}{(n-k)a_k}\Big)^{j_0}\prod_{i=1}^k\Big(\frac{-\lambda_i}{(n-k)a_k}\Big)^{j_i}\mathbb{E}\big[\boldsymbol{\xi}_i^{2j_i}e^{\frac{\lambda_i\boldsymbol{\xi}^2_{i}}{2a_k}}\big].\notag\\
\end{align}
The PDFs $p_{\boldsymbol{\xi}^2_{i}}(\cdot)$ are given by 
\begin{align}
\label{koochil_2}
p_{\boldsymbol{\xi}^2_{i}}(z)=\left\{\begin{array}{cc}
      \frac{1}{\sqrt{2\pi z}}e^{-\frac{z}{2}}   & z>0 \\ 
     0 &   z\le0
\end{array}\right.,\,\,\, i=1,\dots, k.
\end{align}
Then  
\begin{align}
\mathbb{E}\big[\boldsymbol{\xi}_i^{2j_i}e^{\frac{\lambda_i\boldsymbol{\xi}^2_{i}}{2a_k}}\big]&=\frac{1}{\sqrt{2\pi}}\int_0^\infty z^{j_i}e^{\frac{\lambda_i z}{2a_k}}p_{\boldsymbol{\xi_i^2}}(z)dz=\frac{1}{\sqrt{2\pi}}\int_0^\infty z^{j_i-\frac{1}{2}}e^{-\frac{1}{2}(1-\frac{\lambda_i}{a_k})z}dz\notag\\&=\frac{2^{j_i+\frac{1}{2}}}{\sqrt{2\pi}(1-\frac{\lambda_i}{a_k})^{j_i+\frac{1}{2}}}\int_0^\infty \tilde{z}^{j_i-\frac{1}{2}}e^{-\tilde{z}}d\tilde{z}=\frac{2^{j_i}\Gamma(j_i+\frac{1}{2})}{\sqrt{\pi}(1-\frac{\lambda_i}{a_k})^{j_i+\frac{1}{2}}},\,\,\,i=1,\dots,k,
\end{align}
where in the penultimate step we have applied the change of variable $\tilde{z}=\frac{1}{2}(1-\frac{\lambda_i}{a_k})z$ and $\Gamma(\cdot)$ is the Gamma function. It is well-known that $\Gamma(j+\frac{1}{2})=\frac{(2j-1)!!}{2^j}\sqrt{\pi}$ for an integer $j\geq 0$. Therefore, 
\begin{align}
\label{monkey_2}
\mathbb{E}\big[\boldsymbol{\xi}_i^{2j_i}e^{\frac{\lambda_i\boldsymbol{\xi}^2_{i}}{2a_k}}\big]=\frac{(2j_i-1)!!}{(1-\frac{\lambda_i}{a_k})^{j_i+\frac{1}{2}}},\,\,\,\,\,i=1,\dots,k.
\end{align}
Putting (\ref{monkey_1}) and (\ref{monkey_2}) together, we arrive at 
\begin{align}
\label{petr_petr}
\Pr(\boldsymbol{\Delta}>t)&\leq e^{-\frac{t+b_k}{2a_k}} \sum_{\substack{j_0,j_1,\dots, j_k\ge0\\ j_0+j_1+\cdots+j_k=\frac{n-k}{2}}}\frac{(\frac{n-k}{2})!}{j_0!}\Big(1+\frac{t+b_k}{(n-k)a_k}\Big)^{j_0}\prod_{i=1}^k\frac{(2j_i-1)!!}{j_i!}\frac{\big(\frac{-\lambda_i}{(n-k)a_k}\big)^{j_i}}{\big(1-\frac{\lambda_i}{a_k}\big)^{j_i+\frac{1}{2}}}\notag\\
&=e^{-\frac{t+b_k}{2a_k}} \sum_{j_0=0}^{\frac{n-k}{2}}\frac{(\frac{n-k}{2})!}{j_0!}\Big(1+\frac{t+b_k}{(n-k)a_k}\Big)^{j_0}\sum_{\substack{j_1,\dots, j_k\ge0\\ j_1+\cdots+j_k=\frac{n-k}{2}-j_0}}\prod_{i=1}^k\frac{(2j_i-1)!!}{j_i!}\frac{\big(\frac{-\lambda_i}{(n-k)a_k}\big)^{j_i}}{\big(1-\frac{\lambda_i}{a_k}\big)^{j_i+\frac{1}{2}}}.
\end{align}
Reindexing $q=\frac{n-k}{2}-j_0$ and defining $c_{k,q}$ by (\ref{def_c_l}) for $q=0,1,\cdots, \frac{n-k}{2}$, the proof is complete. }

\section{Proof of Proposition~\ref{prop6.1}}
\label{sec_123321}
We need the following result from \cite{Kotz1}.
\begin{lemma}
\label{lem4}
Let the functions $f$ and $\ln f$ admit Maclaurin series expansions with nonzero radii of convergence, given by $f(x)=\sum_{q=0}^\infty u_qx^q$ and $\ln f(x)=v_0+\sum_{q=1}^\infty \frac{v_q}{q}x^q$. Then $u_q=\frac{1}{q}\sum_{j=0}^{q-1}v_{q-j}u_j$ for $q\geq1$.
\end{lemma}
For given integer $k\geq1$,  $y_0>0$ and $y_1,\cdots, y_k\in \mathbb{R}\setminus\{0\}$, consider the function 
\begin{align}
f(x)=y_0\prod_{i=1}^k (1-y_ix)^{-\frac{1}{2}}.
\end{align}
It is straightforward to see that both $f$ and $\ln f$ admit Maclaurin series expansions over the interval $|x|<\frac{1}{\max_{1\leq i\leq k} |y_i|}$. Using the binomial series $(1+x)^{-\frac{1}{2}}=\sum_{j=0}^\infty{-\frac{1}{2}\choose j}x^j$ for  $|x|<1$, the Maclaurin series for $f$ is given by 
\begin{align}
\label{ML_f}
f(x)&\stackrel{}{=}y_0\prod_{i=1}^k \sum_{j=0}^\infty{-\frac{1}{2}\choose j}(-y_ix)^j\notag\\
&=y_0\sum_{j_1,\dots, j_k\ge0}\prod_{i=1}^k{-\frac{1}{2}\choose j_i}(-y_ix)^{j_i}\notag\\
&=y_0\sum_{q=0}^\infty \sum_{\substack{j_1,\dots, j_k\ge0\\ j_1+\cdots+j_k=q}}\prod_{i=1}^k{-\frac{1}{2}\choose j_i}(-y_ix)^{j_i}\notag\\
&=\sum_{q=0}^\infty\Big(y_0\sum_{\substack{j_1,\dots, j_k\ge0\\ j_1+\cdots+j_k=q}}\prod_{i=1}^k{-\frac{1}{2}\choose j_i}(-y_i)^{j_i}\Big)x^q.
\end{align}
But, ${-\frac{1}{2}\choose j}=\frac{-\frac{1}{2}(-\frac{1}{2}-1)\cdots(-\frac{1}{2}-j+1)}{j!}=\frac{(-1)^{j}(2j-1)!!}{2^jj!}$ for every $j\geq 0$. Using this in (\ref{ML_f}), we find $f(x)=\sum_{q=0}^\infty u_qx^q$ where
\begin{eqnarray}
u_q=u_q(y_0,y_1,\cdots, y_k)=\frac{y_0}{2^q}\sum_{\substack{j_1,\dots, j_k\ge0\\ j_1+\cdots+j_k=q}}\prod_{i=1}^k\frac{(2j_i-1)!!}{j_i!}y_i^{j_i},\,\,\,q\ge0
\end{eqnarray}
Using the expansion $\ln(1-x)=-\sum_{q=1}^\infty\frac{x^q}{q}$ for $|x|<1$, the Maclaurin series for $\ln f$ is given by 
\begin{align}
\label{}
  \ln f(x)  =\ln y_0-\frac{1}{2}\sum_{i=1}^k \ln(1-y_ix) =\ln y_0+\frac{1}{2}\sum_{i=1}^k  \sum_{q=1}^\infty\frac{(y_ix)^q}{q}=\ln y_0+\sum_{q=1}^\infty\Big(\frac{1}{2q}\sum_{i=1}^ky_i^q \Big)x^q.
\end{align}
Hence, $\ln f(x)=v_0+\sum_{q=1}^\infty \frac{v_q}{q}x^q$ where 
\begin{align}
\label{}
  v_0=\ln y_0,\quad v_q=v_q(y_1,\cdots, y_k)=\frac{1}{2}\sum_{i=1}^ky_i^q,\,\,\, q\geq 1.
\end{align}
By Lemma~\ref{lem4}, 
\begin{align}
\label{Pro_1}
u_q=\frac{1}{q}\sum_{j=0}^{q-1}v_{q-j}u_j,\,\,\, q\geq1.
\end{align}
Recall $\phi_1,\dots,\phi_k$ in (\ref{def_phi_i}). Setting $y_i=\phi_i$ for $i=1,\dots, k$ and $y_0=\phi_0:=\frac{(\frac{r(A)-k}{2})!}{\prod_{i=1}^k(1-\frac{\lambda_i}{a_k})^{\frac{1}{2}}}$, we observe that $c_{k,q}$ in (\ref{def_c_l}) can be written~as  
\begin{align}
\label{Pro_2}
c_{k,0}=\phi_0,\quad c_{k,q}=2^qu_q(\phi_0,\phi_1,\cdots, \phi_k),\,\,\, q\geq1.
\end{align}  
Let 
\begin{align}
\label{Pro_3}
   d_{k,q}=2^{q+1}v_q(\phi_1,\cdots, \phi_k)=\sum_{i=1}^k(2\phi_i)^q,\,\,\, q\geq1.
\end{align}
By (\ref{Pro_1}), (\ref{Pro_2}) and (\ref{Pro_3}), $\frac{c_{k,q}}{2^q}=\frac{1}{q}\sum_{j=0}^{q-1}\frac{d_{k,q-j}}{2^{q-j+1}}\frac{c_{k,j}}{2^j}$, which can be written as $c_{k,q}=\frac{1}{2q}\sum_{j=0}^{q-1}d_{k,q-j}c_{k,j}$. The proof of Proposition~\ref{prop6.1} is complete. 

 \section{Proof of Proposition~\ref{prop_11_last}}
 \label{secXII}
\textcolor{Black}{ We write $\boldsymbol{\Delta}=\boldsymbol{Q}^+-\boldsymbol{Q}^--\mathrm{tr}(A)$ where $\boldsymbol{Q}^-$ and $\boldsymbol{Q}^+$ are positive-definite quadratic forms given by $\boldsymbol{Q}^-=-\sum_{i=1}^{u}\lambda_{i}\boldsymbol{\xi}^{2}_{i}$ and $\boldsymbol{Q}^+=\sum_{i=u+1}^{r(A)}\lambda_i\boldsymbol{\xi}^{2}_{i}$ and  $\boldsymbol{\xi}_{1},\dots, \boldsymbol{\xi}_{r(A)}$ are independent standard normal random variables. Let us denote $t_0=t+\mathrm{tr}(A)$. We consider two cases:
 \begin{enumerate}
  \item Let $t_0\ge0$. Then 
  \allowdisplaybreaks
  \begin{align}
\label{RK_plus}
  \Pr(\boldsymbol{\Delta}>t)&=\Pr(\boldsymbol{Q}^+>t_0+\boldsymbol{Q}^-)\notag\\&\stackrel{(a)}{=}\mathbb{E}_{Q\sim \boldsymbol{Q}^-}\big[\Pr\big(\boldsymbol{Q}^+>t_0+Q \big)\big] \notag\\
    &\stackrel{}{=}1-\mathbb{E}\big[F_{\boldsymbol{Q}^+}(t_0+\boldsymbol{Q}^-)\big]\notag\\&\stackrel{(b)}{=}1-\mathbb{E}\Big[\sum_{l=0}^\infty c_l^{\text{RK}^+}F_{v+2l}\Big(\frac{t_0+\boldsymbol{Q}^-}{s^{+}} \Big)\Big] \notag\\
      &\stackrel{(c)}{=}1-\sum_{l=0}^\infty c_l^{\text{RK}^+}\mathbb{E}\Big[F_{v+2l}\Big(\frac{t_0+\boldsymbol{Q}^-}{s^{+}} \Big)\Big]\notag\\& \stackrel{(d)}{\le}1-\sum_{l=0}^L c_l^{\text{RK}^+}\mathbb{E}\Big[F_{v+2l}\Big(\frac{t_0+\boldsymbol{Q}^-}{s^{+}} \Big)\Big]\notag\\
       &\stackrel{(e)}{\le}1-\sum_{l=0}^L c_l^{\text{RK}^+}\mathbb{E}\Big[F_{\hat{v}+2l}\Big(\frac{t_0+\boldsymbol{Q}^-}{s^{+}} \Big)\Big]\notag\\
       &\stackrel{(f)}{=}1-\sum_{l=0}^L c_l^{\text{RK}^+}\Big(1-e^{-\frac{t_0}{2s^{+}}}\sum_{m=0}^{\frac{\hat{v}+2l}{2}-1}\frac{1}{2^mm!}\mathbb{E}\Big[\Big(\frac{t_0+\boldsymbol{Q}^-}{s^{+}}\Big)^me^{-\frac{\boldsymbol{Q}^-}{2s^{+}}} \Big]\Big)\notag\\
       &\stackrel{}{=}1-\sum_{l=0}^Lc_l^{\text{RK}^+}+e^{-\frac{t_0}{2s^{+}}}\sum_{l=0}^L c_l^{\text{RK}^+}\sum_{m=0}^{\frac{\hat{v}}{2}+l-1}\frac{1}{2^mm!}\mathbb{E}\Big[\Big(\frac{t_0-\sum_{i=1}^{u}\lambda_i\boldsymbol{\xi}^2_i }{s^{+}}\Big)^me^{\frac{\sum_{i=1}^{u}\lambda_i\boldsymbol{\xi}^2_i }{2s^{+}}} \Big],
\end{align}
where $(a)$ follows from conditioning and independence of $\boldsymbol{Q}^-$ and $\boldsymbol{Q}^+$, $(b)$ is due to the Ruben-Kotz expansion in (\ref{RK_expansion}), $(c)$ is due to the Monotone Convergence Theorem~\cite{Axler}, $(d)$  is due to truncation at $l=L$, $(e)$ is due to $v\leq \hat{v}$ and the fact that the CDF for the chi-squared distribution is non-increasing with respect to the degrees of freedom and $(f)$ is due to (\ref{even_n}) and the fact that $\hat{v}+2l$ is even.   
  \item Let $t_0< 0$. Following similar lines of reasoning in (\ref{RK_plus}), we find  
  \allowdisplaybreaks
  \begin{align}
\label{RK_minus}
  \Pr(\boldsymbol{\Delta}>t)&=\Pr(\boldsymbol{Q}^-<\boldsymbol{Q}^+-t_0 \big)=\mathbb{E}\big[F_{\boldsymbol{Q}^-}(\boldsymbol{Q}^+-t_0 )\big]= \sum_{l=0}^\infty c_l^{\text{RK}^-}\mathbb{E}\Big[F_{u+2l}\Big(\frac{\boldsymbol{Q}^+-t_0}{s^-} \Big)\Big]\notag\\
 &\stackrel{(a)}{\le} \sum_{l=L+1}^\infty c_l^{\text{RK}^-}+\sum_{l=0}^L c_l^{\text{RK}^-}\mathbb{E}\Big[F_{u+2l}\Big(\frac{\boldsymbol{Q}^+-t_0}{s^-} \Big)\Big]\notag\\
  &\stackrel{(b)}{\le} \sum_{l=L+1}^\infty c_l^{\text{RK}^-}+\sum_{l=0}^L c_l^{\text{RK}^-}\mathbb{E}\Big[F_{\check{u}+2l}\Big(\frac{\boldsymbol{Q}^+-t_0}{s^-} \Big)\Big]\notag\\
  &\stackrel{(c)}{=} 1-\sum_{l=0}^L c_l^{\text{RK}^{-}}+\sum_{l=0}^L c_l^{\text{RK}^-}\Big(1-e^{\frac{t_0}{2s^-}}\sum_{m=0}^{\frac{\check{u}+2l}{2}-1}\frac{1}{2^mm!}\mathbb{E}\Big[\Big(\frac{\boldsymbol{Q}^+-t_0}{s^-}\Big)^me^{-\frac{\boldsymbol{Q}^+}{2s^-}} \Big]\Big)\notag\\
 &= 1-e^{\frac{t_0}{2s^-}}\sum_{l=0}^L c_l^{\text{RK}^-}\sum_{m=0}^{\frac{\check{u}}{2}+l-1}\frac{1}{2^mm!}\mathbb{E}\bigg[\Big(\frac{\sum_{i=u+1}^{r(A)}\lambda_i\boldsymbol{\xi}^{2}_{i}-t_0}{s^-}\Big)^me^{-\frac{\sum_{i=u+1}^{r(A)}\lambda_i\boldsymbol{\xi}^{2}_{i}}{2s^-}} \bigg],
\end{align}
where $(a)$ is due to $\mathbb{E}[F_{u+2l}(\frac{\boldsymbol{Q}^+-t_0}{s^-} )]\leq1$ for all $l\ge L+1$, $(b)$ is due to $u\geq \check{u}$ and the fact that the CDF for the chi-squared distribution is non-increasing with respect to the degrees of freedom and $(c)$ is due to $ \sum_{l=L+1}^\infty c_l^{\text{RK}^-}=1- \sum_{l=0}^L c_l^{\text{RK}^-}$.
\end{enumerate}
 Computing the expectations in (\ref{RK_plus}) and (\ref{RK_minus}) is similar to computing the expectation in (\ref{expectation_App_E}) in Appendix~\ref{A5}. They are given by $\sum_{q=0}^{m}\frac{m!}{(m-q)!}(\frac{t_0}{s^+})^{m-q}c^+_{q}$ and $\sum_{q=0}^{m}\frac{m!}{(m-q)!}(-\frac{t_0}{s^-})^{m-q}c^-_{q}$, respectively, where $c^-_q$ and $c^+_q$ are defined in (\ref{c_neg_q}) and (\ref{c_plus_q}). }
 
\textcolor{Black}{ To assess the tightness of the bound in (\ref{RK_plus}), we note that two inequalities are introduced only in steps $(d)$ and $(e)$. Denoting $\boldsymbol{R}=\frac{t_0+\boldsymbol{Q}^-}{s^{+}}$ for simplicity, the difference between the upper bound on $\Pr(\boldsymbol{\Delta}>t)$ in (\ref{RK_plus}) and $\Pr(\boldsymbol{\Delta}>t)$ itself is given by 
\begin{align}
\label{E_Bound}
    \textrm{gap in $(d)$}+\textrm{gap in $(e)$}&=\sum_{l=L+1}^\infty c_l^{\text{RK}^+}\mathbb{E}[F_{v+2l}(\boldsymbol{R})] +\sum_{l=0}^L c_l^{\text{RK}^+}\mathbb{E}[F_{v+2l}(\boldsymbol{R})-F_{\hat{v}+2l}(\boldsymbol{R})]\notag\\
    &\leq \sum_{l=L+1}^\infty c_l^{\text{RK}^+}+\sum_{l=0}^L c_l^{\text{RK}^+}\delta_{v+2l}\mathds{1}_{v\neq \hat{v}},
\end{align}
where the last step follows from the definition of the numbers $\delta_k$ in (\ref{KSD}). We also noted  that if $v\neq\hat{v}$ (i.e., $v$ is odd), then $\hat{v}=v+1$, meaning that $v+2l$ and $\hat{v}+2l$ are consecutive integers. The expression on the right hand side in (\ref{E_Bound}) is exactly $\mathsf{Gap}^+$ given in (\ref{RK_error_plus}). The proof of (\ref{RK_error_minus}) is similar and is omitted for brevity. }

\section{Concluding Remarks}
We studied concentration of measure for Gaussian quadratic chaos in the non-asymptotic regime, where existing bounds were improved and several new bounds were proposed.  We tightened HWI by increasing the absolute constant $\kappa$ in its formula from the best currently known value of $0.125$ to at least $0.1457$ in the symmetric case. This was achieved by determining the best quadratic upper bound on the function $-\ln(1-x)$ around $x=0$. In the positive-semidefinite case, it was shown that LMI implies HWI with $\kappa=1-\frac{\sqrt{3}}{2}$. We verified that a larger value $\kappa=\frac{9-\sqrt{17}}{32}$ is attainable by improving LMI itself. This improvement manifests in two forms, namely, the augmented LMI and the optimal LMI. The latter requires numerically computing the unique root of a quintic polynomial in the interval $(\frac{2}{3},1)$. The former approximates this unique root with a closed-form expression, where the underlying approximation error is less than $0.035$ regardless of the positive-semidefinite matrix~$A$ and the tail parameter $t$. 

We explored beyond HWI for a symmetric matrix $A$ by deriving the best polynomial upper bound of degree $m+1$ on the function $-\ln(1-x)$ around $x=0$, where $m\geq1$ is arbitrary. As a result, we developed a sequence of concentration bounds indexed by $m$, where the case $m=1$ recovers HWI. These bounds are expressed in terms of Schatten norms of $A$. It was conjectured that they undergo a phase transition in the sense that for a given matrix $A$, there exists a threshold $\tau_c\ge0$ such that if $t<\tau_c$, then $m=1$ (HWI) is the tightest and if $t>\tau_c$, then $m=\infty$ results in the sharpest bound. This led to a novel bound, called the $m_\infty$~bound, which depends on $A$ only through its rank $r(A)$ and its operator norm $\|A\|$. An upper bound was also established on $\tau_c$.  

All three HWI, LMI, and the $m_\infty$~inequality rely on Markov's inequality. Following a line of reasoning that avoids Markov's inequality, we introduced the strong $\chi^2$~inequality in the symmetric case and its relaxed version, the weak $\chi^2$~inequality, in the positive-semidefinite~case. To further investigate the $m_\infty$ and $\chi^2$ bounds, we restricted our attention to positive-definite matrices and explored all available concentration bounds that depend on $A$ only through its operator norm. Five candidates were examined, namely, the $m_\infty$~bound, the weak $\chi^2$~bound, relaxed versions of HWI and LMI, and the large-deviations bound. Among these, the tightest bound is always either the \(m_{\infty }\) bound or the weak \(\chi ^{2}\) bound. If $n=2,4,6$, the weak $\chi^2$~bound is tighter than the $m_\infty$~bound. For even~$n\ge8$, the $m_\infty$~bound is sharper than the weak $\chi^2$~bound if and only if the ratio $\frac{t}{\|A\|}$ lies inside an open interval which expands indefinitely as $n$ grows.

 Using a simple conditioning argument, modified versions of $m_\infty$, strong $\chi^2$, and HW inequalities of various orders were also proposed, offering significant improvements over their original counterparts. While any two different orders for the modified $m_\infty$ and strong $\chi^2$ inequalities must share the same parity, this condition is not required for the modified HW inequality. It was shown that:    
 \begin{itemize}
  \item The modified strong $\chi^2$~inequality of order $k+2$ is sharper than the one of order $k$ for all~$t$. 
  \item The modified $m_\infty$ inequality of order $k+2$ is sharper than the one of order $k$ for sufficiently large~$t$. 
  \item If $a_{k+1}=a_k$, then the modified HWI of order $k+1$ is \textit{looser} than the one of order $k$ for all $t$, while for $a_{k+1}<a_k$, the higher-order inequality is \textit{sharper} for sufficiently large~$t$. 
\end{itemize} 
Although increasing the order typically yields tighter bounds, it comes at the expense of higher computational complexity. This trade-off is particularly critical for the modified \(m_{\infty }\) and \(\chi ^{2}\) bounds. Nevertheless, a practical advantage of the modified bounds is that they do not involve special functions or numerical integration. 

\textcolor{Black}{By extending the classical Ruben-Kotz expansion to general symmetric quadratic forms, we introduced RK inequalities of various truncation levels $L\geq0$. Although the computational complexity scales at the order of $\Omega(L^2)$, a key advantage of RKIs is that they provide an upper bound on the approximation gap. Specifically, when either $t+\mathrm{tr}(A)<0$ and $A$ has an even number $u$ of negative eigenvalues, or $t+\mathrm{tr}(A)\ge0$ and $A$ has an even number $v$ of positive eigenvalues, the approximation gap vanishes as $L$ approaches infinity. If these conditions are not satisfied, the approximation gap never exceeds $\max\{\delta_{u-1},\delta_v\}$, where $\delta_k$ is given in (\ref{closed_form_KSD}). }

These bounds proved highly effective in two applications, namely, change detection in the covariance matrix of a Gaussian vector and performance evaluation for unitary space-time coding over a block-fading MISO~channel.


\appendices

\section{Proof of Lemma~\ref{lem_1}}
\label{A1}

Define the function $f$ by 
\begin{align}
f(x)=x+\frac{x^2}{2}+\frac{x^3}{3}+\cdots+\frac{x^m}{m}+a|x|^{m+1}+\ln(1-x).
\end{align}
We study the cases $x>0$ and $x<0$ separately:
\begin{enumerate}
  \item Let $x>0$. Then $f(x)=x+\frac{x^2}{2}+\frac{x^3}{3}+\cdots+\frac{x^m}{m}+ax^{m+1}+\ln(1-x)$ and we get 
\begin{align}
\label{push}
f'(x)&=1+x+x^2+\cdots+x^{m-1}+(m+1)ax^m-\frac{1}{1-x}\notag\\
&=\frac{1-x^m+(m+1)ax^m(1-x)-1}{1-x}\notag\\
&=\frac{x^m((m+1)a-1-(m+1)ax)}{1-x}.
\end{align}
The only first-order critical numbers for $f$ are $0$ and $x_0=1-\frac{1}{(m+1)a}$. If $a\leq \frac{1}{m+1}$, then one easily checks that $f'(x)<0$ over the open interval $(0,1)$. Since $f(0)=0$, the Mean Value Theorem (MVT) implies that $f(x)<0$ over $(0,1)$ and hence, (\ref{theone}) can not hold for any $b>0$. Conversely, let $a>\frac{1}{m+1}$. Then $0<x_0<1$, $f'(x)>0$ over the interval $(0, x_0)$ and $f'(x)<0$ over the interval $(x_0,1)$. Note that $\lim_{x\to1^-}f(x)=-\infty$. This tells us that $f$ rises above $0$ on right of $x=0$, reaches a maximum value at the critical number~$x_0$ and then goes down from there on and escapes to $-\infty$ on left of $x=1$. As such, $f$ has a unique zero (x-crossing) somewhere over the interval $(x_0,1)$. We let $b$ be this zero. Then~(\ref{gen_ineq1}) clearly holds for $0\leq x\leq b$. Moreover, the equation $f(b)=0$ gives $a=\theta_m(b)$ as defined in (\ref{theta_func}).   
  \item Let $x<0$. If $m$ is odd, then $f(x)=x+\frac{x^2}{2}+\frac{x^3}{3}+\cdots+\frac{x^m}{m}+ax^{m+1}+\ln(1-x)$ as in the previous case and one easily checks that for $a>\frac{1}{m+1}$, $f'(x)$ given in (\ref{push}) is negative for every $x<0$. Then MVT implies that $f(x)>0$ for every $x<0$ as desired. If $m$ is even, then $f(x)=x+\frac{x^2}{2}+\frac{x^3}{3}+\cdots+\frac{x^m}{m}-ax^{m+1}+\ln(1-x)$. A similar computation as in~(\ref{push}) gives that 
  \begin{align}
f'(x)=\frac{-x^m((m+1)a+1-(m+1)ax)}{1-x}.
\end{align}
Since $m$ is even, we have $f'(x)<0$ for every $x<0$ and MVT applies once again to confirm that $f(x)>0$ for $x<0$.
\end{enumerate}

\section{Proof of Lemma~\ref{lem_2}}
\label{A2}
Assume $m\geq1$ is an integer and fix $\frac{m+1}{m+2}<a<1$. Define the function $f$ by 
\begin{align}
\hskip1cmf(x)=x+\frac{x^2}{2}+\frac{x^3}{3}+\cdots+\frac{x^m}{m}+\frac{x^{m+1}}{(m+1)(1-ax)}+\ln(1-x),\,\,\, 0\le x<1.
\end{align}
One can easily check that  
\begin{align}
f'(x)=\frac{x^{m+1}\big((m+2)a-m-1-a((m+1)a-m)x\big)}{(m+1)(1-x)(1-ax)^2}.
\end{align}
Thanks to the condition $\frac{m+1}{m+2}<a<1$, the only positive first-order critical number for $f$ over the interval $(0,1)$ is given by $x_0=\frac{m+1-(m+2)a}{a(m-(m+1)a)}$. Moreover, $f'(x)>0$ for $0<x<x_0$ and $f'(x)<0$ for $x_0<x<1$. Since $f(0)=0$, then MVT implies that $f(x)>0$ for every $x$ in the interval\footnote{In fact, $f$ has a unique x-crossing at some number $x_1$ inside the interval $(x_0,1)$ and $f(x)\ge 0$ for every $0\le x\le x_1$. Hence, (\ref{gen_ineq2}) also holds with $b=x_1$. However, unlike $x_0$, the number $x_1$ does not admit a closed-form expression. } $(0,x_0]$ and hence, (\ref{gen_ineq2}) holds with $b=x_0$.

\section{The choice of $a=a_0$  in (\ref{a0}) yields the largest $\kappa$}
\label{A3}
\begin{figure}
   \centering
    \includegraphics[width=0.45\textwidth]{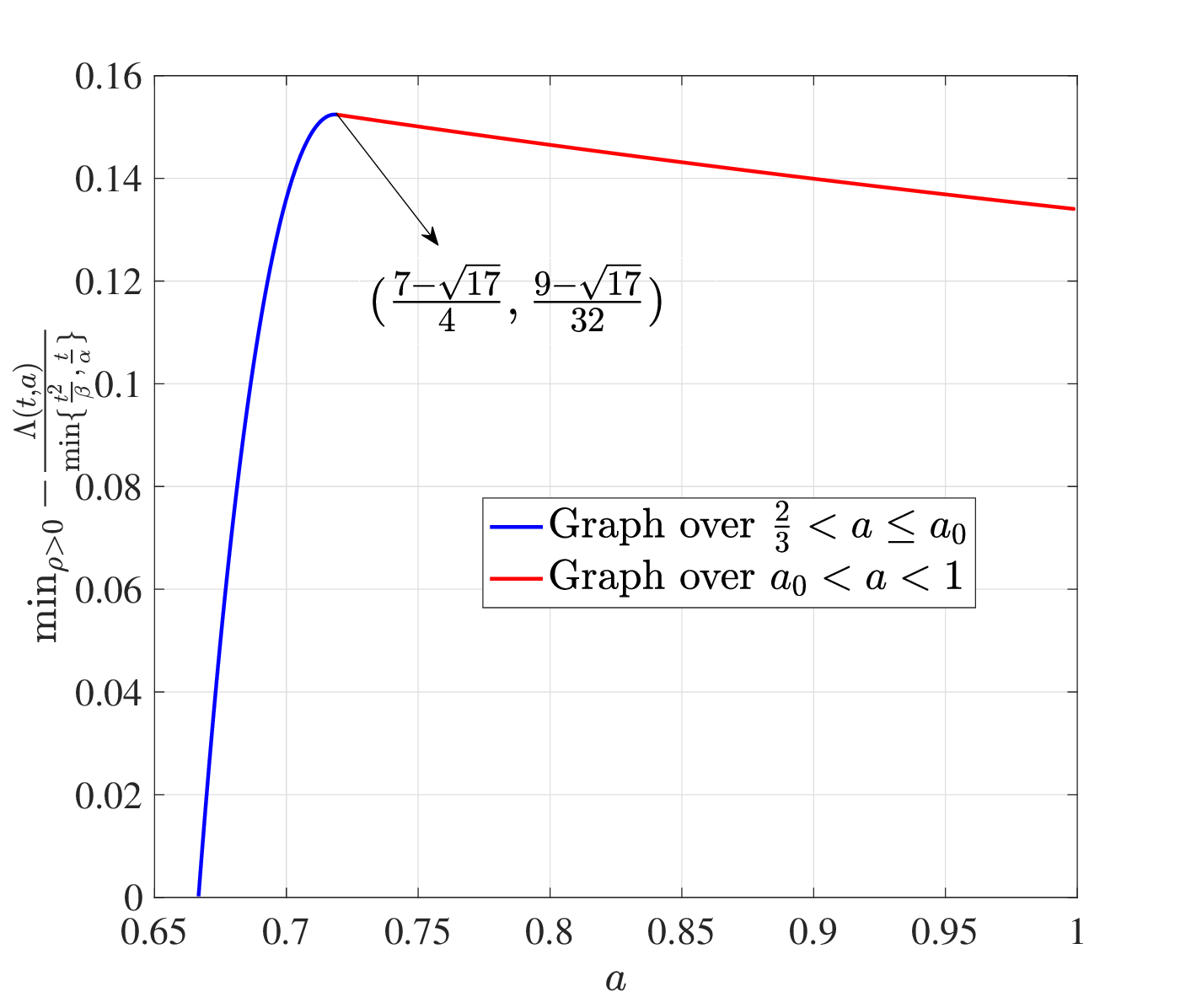}
     \caption{ The graph of the function of $a$ that sits on the right-hand side of (\ref{day}) over the interval $(\frac{2}{3},1)$. It achieves its absolute maximum value of $\frac{9-\sqrt{17}}{32}$ at $a=a_0=\frac{7-\sqrt{17}}{4}$.}
    \label{pict_11111}
\end{figure}

In order to compare $\Lambda(t,a)$ in (\ref{en_exp}) with $\min\{\frac{t^2}{\beta}, \frac{t}{\alpha}\}$, we address the cases $c\leq 1$ and $c>1$ separately: 
\begin{enumerate}
  \item Let $c\leq 1$, or equivalently, $\frac{2}{3}<a\leq a_0$. Then 
\begin{align}
-\frac{\Lambda(t,a)}{\min\{\frac{t^2}{\beta}, \frac{t}{\alpha}\}}=\left\{\begin{array}{cc}
    \frac{1}{2a\rho}-\frac{1}{2a^{2}\rho^2}\big(\sqrt{1+2a\rho}-1\big)  & 0<\rho\leq c   \\
     \frac{b}{2\rho}-\frac{b^2}{4(1-ab)\rho^2} &   c<\rho\leq 1\\
    \frac{b}{2}-\frac{b^2}{4(1-ab)\rho} & \rho>1
\end{array}\right..
\end{align}  
One can easily check that this is a continuous function of $\rho$ and it achieves its absolute minimum value of $\frac{b}{2}-\frac{b^2}{4(1-ab)}=\frac{6a^3-a^2-8a+4}{4a^2(2a^2-3a+1)}$  at $\rho=1$.
  \item Let $c>1$, or equivalently, $ a_0<a<1$.  Then 
  \begin{align}
-\frac{\Lambda(t,a)}{\min\{\frac{t^2}{\beta}, \frac{t}{\alpha}\}}=\left\{\begin{array}{cc}
    \frac{1}{2a\rho}-\frac{1}{2a^{2}\rho^2}\big(\sqrt{1+2a\rho}-1\big)  & 0<\rho\leq 1   \\
     \frac{1}{2a}-\frac{1}{2a^2\rho}(\sqrt{1+2a\rho}-1) &   1<\rho\leq c\\
    \frac{b}{2}-\frac{b^2}{4(1-ab)\rho} & \rho>c
\end{array}\right..
\end{align}   
Once again, this is a continuous function of $\rho$ and it achieves its absolute minimum value of $\frac{1}{2a}-\frac{1}{2a^2}(\sqrt{1+2a}-1)$ at $\rho=1$.
\end{enumerate}
We have shown that 
\begin{align}
\label{day}
\min_{\rho>0}-\frac{\Lambda(t,a)}{\min\{\frac{t^2}{\beta^2}, \frac{t}{\alpha}\}}=\left\{\begin{array}{cc}
    \frac{6a^3-a^2-8a+4}{4a^2(2a^2-3a+1)}  & \frac{2}{3}<a\leq a_0   \\
     \frac{1}{2a}-\frac{1}{2a^2}(\sqrt{1+2a}-1) &  a_0<a<1 
\end{array}\right..
\end{align}
We need to choose $\frac{2}{3}<a<1$ such that the right hand side in (\ref{day}) is as large as possible. Fig.~\ref{pict_11111} shows the graph of this function of $a$ over the interval $(\frac{2}{3},1)$. We see that it achieves its absolute maximum value of $\frac{9-\sqrt{17}}{32}$ at $a=a_0=\frac{7-\sqrt{17}}{4}$.

\section{Proof of (\ref{final_claim})}
\label{A4}
For an even integer $n\geq 8$, let us denote $m=\frac{n}{2}\geq 4$. During the proof, we adopt the notation $\forall_\infty m$ to say that a statement indexed by the integer $m$ holds for all but finitely many values for $m$, i.e., it eventually holds for all sufficiently large~$m$. We begin by showing that 
\begin{align}
\label{App_e_1}
\forall_\infty m\geq 4,\,\,\, r'_{2m}>2m.
\end{align}
 This establishes the second claim in~(\ref{final_claim}). Recall from the proof of Proposition~\ref{prop_7} that $r_n, r'_n$ are the only two positive zeros for the polynomial $P_n^{(1)}(r)-P_n^{(2)}(r)$ and that this polynomial is negative if and only if $r_n<r<r'_n$. In order to verify (\ref{App_e_1}), it is enough to show that
\begin{align}
\label{this_that}
\forall_\infty m\geq 4,\,\,\, P_{2m}^{(1)}(2m)<P_{2m}^{(2)}(2m).
\end{align} 
We have $P_{2m}^{(2)}(2m)=\frac{1}{\sqrt{e}}\sum_{i=0}^{m-1}\frac{1}{2^ii!}(1+2m)^{i}$. We show that the last term in this sum is eventually larger than $P_{2m}^{(1)}(2m)=2^{m}$, i.e., 
\begin{align}
\label{fell0}
\forall_\infty m\geq 4,\,\,\, 2^m<\frac{1}{\sqrt{e}}\frac{(1+2m)^{m-1}}{2^{m-1}(m-1)!}=\frac{1}{\sqrt{e}}\frac{(m+\frac{1}{2})^{m-1}}{(m-1)!}.
\end{align}
This is easily done by invoking the Stirling's approximation formula. Let us recall the inequality\footnote{See Section 2.9 in \cite{Feller}.} 
\begin{align}
\label{fell1}
k!<\sqrt{2\pi}\,k^{k+\frac{1}{2}}e^{-k+\frac{1}{12k}},\,\,\, k\geq 1,
\end{align}
Plugging $k=m-1\geq 3$ in (\ref{fell1}) and observing that $\sqrt{2\pi}e^{\frac{1}{12(m-1)}}\leq \sqrt{2\pi}e^{\frac{1}{36}}\approx 2.577<e$, we get $(m-1)!<(m-1)^{m-\frac{1}{2}}e^{2-m}$ and hence, 
\begin{align}
\label{fell2}
\frac{e^{m-2}}{(m-1)^{m-\frac{1}{2}}}<\frac{1}{(m-1)!},\,\,\, m\geq 4.
\end{align}
By (\ref{fell2}), the proof of (\ref{fell0}) is complete if we can verify the stronger statement 
\begin{align}
\label{fell3}
\forall_\infty m\geq 4,\,\,\, 2^m<\frac{1}{\sqrt{e}}\frac{(m+\frac{1}{2})^{m-1}e^{m-2}}{(m-1)^{m-\frac{1}{2}}}.
\end{align}
Note that $\frac{(m+\frac{1}{2})^{m-1}}{(m-1)^{m-\frac{1}{2}}}>\frac{(m+\frac{1}{2})^{m-1}}{(m+\frac{1}{2})^{m-\frac{1}{2}}}=\frac{1}{(m+\frac{1}{2})^{\frac{1}{2}}}$. Thus, (\ref{fell3}) follows if 
\begin{align}
\forall_\infty m\geq 4,\,\,\, 2^m<\frac{1}{\sqrt{e}}\frac{e^{m-2}}{(m+\frac{1}{2})^{\frac{1}{2}}}.
\end{align} 
But, this is certainly true as the right hand side scales like $e^{m-O(\frac{1}{2}\ln m)}$. The proof of (\ref{App_e_1}) is now complete. 

To verify the first claim in (\ref{final_claim}), we show that 
\begin{align}
\label{scnm_1}
\forall \epsilon>0,\, \forall_\infty m\geq 4,\,\,\, r_{2m}<\epsilon.
\end{align}
Following the discussion before (\ref{this_that}), it suffices to prove that 
\begin{align}
\label{kaaff}
\forall_\infty m\geq 4,\,\,\, P_{2m}^{(1)}(\epsilon)<P_{2m}^{(2)}(\epsilon).
\end{align} 
The proof follows similar lines of reasoning to show $P_n^{(1)}(1)<P_n^{(2)}(1)$ in Lemma~\ref{lem_lem}. Fix $\epsilon>0$. Let $\boldsymbol{x}\sim \mathrm{Poisson}(\frac{1+\epsilon}{2})$ be a Poisson random variable with parameter $\frac{1+\epsilon}{2}$. Then 
\begin{align}
\label{kaaf}
P_{2m}^{(2)}(\epsilon)&=\frac{1}{\sqrt{e}}\sum_{i=0}^{m-1}\frac{(1+\epsilon)^i}{2^ii!}=e^{\frac{\epsilon}{2}}\sum_{i=0}^{m-1}\frac{e^{-\frac{1+\epsilon}{2}}(\frac{1+\epsilon}{2})^i}{i!}=e^{\frac{\epsilon}{2}}\Pr(\boldsymbol{x}\le m-1)\notag\\
&=e^{\frac{\epsilon}{2}}(1-\Pr(\boldsymbol{x}\ge m))\stackrel{(a)}{\geq} e^{\frac{\epsilon}{2}}\Big(1-\frac{\mathbb{E}[\boldsymbol{x}^2]}{m^2}\Big)\stackrel{(b)}{=} e^{\frac{\epsilon}{2}}\Big(1-\frac{\frac{1+\epsilon}{2}+(\frac{1+\epsilon}{2})^2}{m^2}\Big)=e^{\frac{\epsilon}{2}}\Big(1-\frac{\epsilon^2+4\epsilon+3}{4m^2}\Big),
\end{align}
where $(a)$ follows by Markov's inequality and $(b)$ is due to $\mathbb{E}[\boldsymbol{x}^2]=\frac{1+\epsilon}{2}+(\frac{1+\epsilon}{2})^2$. Also, $P_{2m}^{(1)}(\epsilon)=(1+\frac{\epsilon}{2m})^m$. By (\ref{kaaf}), if we can show that 
\begin{align}
\label{plk}
\forall_\infty m\geq 4,\,\,\, \Big(1+\frac{\epsilon}{2m}\Big)^m<e^{\frac{\epsilon}{2}}\Big(1-\frac{\epsilon^2+4\epsilon+3}{4m^2}\Big),
\end{align} 
then (\ref{kaaff}) follows. Define the function 
\begin{align}
f_{\epsilon}(x)=\frac{\epsilon}{2}+\ln\Big(1-\frac{\epsilon^2+4\epsilon+3}{4x^2}\Big)-x\ln\Big(1+\frac{\epsilon}{2x}\Big),\,\,\, x>\frac{1}{2}\sqrt{\epsilon^2+4\epsilon+3}.
\end{align}
The statement in (\ref{plk}) can be written as $\forall_\infty m\geq 4, e^{f_{\epsilon}(m)}>1$. Thus, it is enough to prove that $f_{\epsilon}(x)>0$ for sufficiently large $x$. This is accomplished if we can verify 
\begin{align}
\label{alps_1}
\lim_{x\to\infty}f_{\epsilon}(x)=0
\end{align}
and
\begin{align}
\label{alps_2}
f'_{\epsilon}(x)<0, 
\end{align}
for all sufficiently large $x$. Checking (\ref{alps_1}) is trivial. To verify (\ref{alps_2}), we use a similar argument and show that  
\begin{align}
\label{alps_3}
\lim_{x\to\infty}f'_{\epsilon}(x)=0
\end{align}
and
\begin{align}
\label{alps_4}
f''_{\epsilon}(x)>0, 
\end{align}
for all sufficiently large $x$. We have 
\begin{align}
f'_{\epsilon}(x)=\frac{2(\epsilon^2+4\epsilon+3)}{x(4x^2-\epsilon^2-4\epsilon-3)}+\frac{\epsilon}{2x+\epsilon}-\ln\Big(1+\frac{\epsilon}{2x}\Big).
\end{align}
Then (\ref{alps_3}) clearly holds. Also, 
\begin{align}
f''_{\epsilon}(x)=\frac{p_{\epsilon}(x)}{x^2(2x+\epsilon)^2(4x^2-\epsilon^2-4\epsilon-3)^2},
\end{align}
where $p_{\epsilon}(x)$ is a polynomial of degree $5$ in $x$ given by
\begin{align}
p_{\epsilon}(x)=16\epsilon^2 x^5+\textrm{terms of lower degrees in $x$}.
\end{align}
 The leading term in the polynomial $p_{\epsilon}(x)$ is $16\epsilon^2 x^5$. Since the coefficient $16\epsilon^2$ behind $x^5$ is positive, then~(\ref{alps_4}) must hold for all sufficiently large $x$ depending on $\epsilon$ of course. The proof of (\ref{scnm_1}) is complete. 
 
\section{Proof of Proposition~\ref{prop_9}}
\label{A5}
\textcolor{Black}{For simplicity of presentation,  we assume $A$ is full-rank, i.e., $r(A)=n$. Let $\boldsymbol{\xi}_1,\dots, \boldsymbol{\xi}_n$ be as in Section~\ref{sec11}. Then
\begin{align}
\label{chi2_proof}
\Pr(\boldsymbol{\Delta}>t)=\Pr\big(\sum_{i=1}^{n}\lambda_i\boldsymbol{\xi}_i^2>t+\mathrm{tr}(A)\big)&=\Pr\big(\sum_{i=k+1}^{n}\lambda_i\boldsymbol{\xi}_i^2>t+\mathrm{tr}(A)-\sum_{i=1}^{k}\lambda_i\boldsymbol{\xi}_i^2\big)\notag\\
&\leq\Pr\big(\lambda_{\max}(A)\chi^2_{n-k}>t+\mathrm{tr}(A)-\sum_{i=1}^{k}\lambda_i\boldsymbol{\xi}_i^2\big)\notag\\
&=1-\mathbb{E}\Big[F_{n-k}\Big(\frac{t+\mathrm{tr}(A)-\sum_{i=1}^{k}\lambda_i\boldsymbol{\xi}_i^2}{\lambda_{\max}(A)}\Big)\Big],
\end{align}  
where the penultimate step is due to $\sum_{i=k+1}^{n}\lambda_i\boldsymbol{\xi}_i^2\leq \lambda_{\max}(A)\sum_{i=k+1}^{n}\boldsymbol{\xi}_i^2=\lambda_{\max}(A)\chi^2_{n-k}$ and the last step follows from conditioning and the assumption that $\lambda_{\max}(A)>0$. By (\ref{poloi_4}), $t+\mathrm{tr}(A)-\sum_{i=1}^{k}\lambda_i\boldsymbol{\xi}^2_i\ge 0$. Since $n-k$ is even, it follows from (\ref{even_n}) and (\ref{chi2_proof}) that
\begin{align}
\label{modified_chi2}
\Pr(\boldsymbol{\Delta}>t)\leq e^{-\frac{t+\mathrm{tr}(A)}{2\lambda_{\max}(A)}}\sum_{m=0}^{\frac{n-k}{2}-1}\frac{1}{2^mm!}\mathbb{E}\bigg[\Big(\frac{t+\mathrm{tr}(A)-\sum_{i=1}^{k}\lambda_i\boldsymbol{\xi}^2_i}{\lambda_{\max}(A)}\Big)^me^{\frac{\sum_{i=1}^{k}\lambda_i\boldsymbol{\xi}^2_i}{2\lambda_{\max}(A)}} \bigg].
\end{align} 
Following similar lines of reasoning that led from (\ref{monkey_1}) to (\ref{petr_petr}) in Section~\ref{sec11}, the expectation in (\ref{modified_chi2}) can be computed in closed form after one expands $(\frac{t+\mathrm{tr}(A)-\sum_{i=1}^{k}\lambda_i\boldsymbol{\xi}^2_i}{\lambda_{\max}(A)})^m$ according to multinomial expansion. Omitting the details for brevity, we arrive at 
\begin{align}
\label{expectation_App_E}
 \mathbb{E}\bigg[\Big(\frac{t+\mathrm{tr}(A)-\sum_{i=1}^{k}\lambda_i\boldsymbol{\xi}^2_i}{\lambda_{\max}(A)}\Big)^me^{\frac{\sum_{i=1}^{k}\lambda_i\boldsymbol{\xi}^2_i}{2\lambda_{\max}(A)}} \bigg]=\sum_{q=0}^{m}\frac{m!}{(m-q)!}\Big(\frac{t+\mathrm{tr}(A)}{\lambda_{\max}(A)}\Big)^{m-q}\tilde{c}_{k,q},
\end{align}
 where the coefficients $\tilde{c}_{k,q}$ are given in (\ref{tilde_c_kl}). The proof of the recursion in (\ref{pill_111}) is identical to the proof of (\ref{Pill_000}), replacing $\phi_i$ with $\tilde{\phi}_i$ and $c_{k,0}$ with $\tilde{c}_{k,0}$.  }

 \section{Proof of Proposition~\ref{prop_8}}
 \label{A6}
\textcolor{Black}{ Assume (\ref{poloi_3}) holds. Then $t+b_k-\sum_{i=1}^k\lambda_iz_i\ge0$ for all $z_1,\dots,z_k\ge0$. Applying the original HWI to~(\ref{koochil_1}), 
 \begin{align}
 \label{rtbn}
\Pr(\boldsymbol{\Delta}>t)\leq \mathbb{E}_{z_i\sim \boldsymbol{\xi}_i^2}\Big[\exp\Big(-\kappa\min\Big\{\frac{(t+b_k-\sum_{i=1}^k \lambda_iz_i)^2}{\sum_{i=k+1}^n\lambda_i^2}, \frac{t+b_k-\sum_{i=1}^k \lambda_iz_i}{a_k}\Big\}\Big)\Big].
\end{align}
Using (\ref{poloi_3}) once more, $\frac{(t+b_k-\sum_{i=1}^k \lambda_iz_i)^2}{\sum_{i=k+1}^n\lambda_i^2}\geq  \frac{t+b_k-\sum_{i=1}^k \lambda_iz_i}{a_k}$ for all $z_1,\dots, z_k\ge0$. Then (\ref{rtbn}) becomes 
\begin{align}
\Pr(\boldsymbol{\Delta}>t)\leq \mathbb{E} \big[e^{- \frac{\kappa(t+b_k-\sum_{i=1}^k \lambda_i\boldsymbol{\xi}_i^2)}{a_k}}\big]=e^{-\frac{\kappa(t+b_k)}{a_k}}\prod_{i=1}^k \mathbb{E}\big[e^{\frac{\kappa \lambda_i \boldsymbol{\xi}_i^2}{a_k}}\big]=e^{-\frac{\kappa(t+b_k)}{a_k}}\prod_{i=1}^k \Big(1-2\kappa\frac{ \lambda_i}{a_k}\Big)^{-\frac{1}{2}},
\end{align}
where the last step uses the fact that $\boldsymbol{\xi}_i^2$ is a $\chi^2_1$ random variable.}

\end{document}